\documentclass[11pt]{article}
 
\usepackage{amsmath}
\usepackage{amssymb}
\usepackage{amsthm}
\usepackage{siunitx}
\usepackage{natbib}
\usepackage{booktabs}
\usepackage[margin=1in]{geometry}
\usepackage{graphicx}
\usepackage{hyperref}
\usepackage{url}
\usepackage{multirow}
\usepackage{siunitx}
\usepackage{caption}
\usepackage{subcaption}
\usepackage{algorithm}
\usepackage{algpseudocode}
\numberwithin{equation}{section}

\newcommand\comment[1]{{}}

\def\E{\mathbb E}
\newcommand\Fi[1]{F_{#1}^{\leftarrow}}
\def\ind{\perp\!\!\!\perp}

\def\I{\mathbb I}
\def\P{\mathbb P}
\def\Q{\mathbb Q}
\def\R{\mathbb R}
\def\RV{\text{RV}}
\def\N{\mathbb N}
\def\Z{\mathbb Z}
\def\seqSet{\mathcal{S}_{\alpha}} 
\newcommand{\multiplier}[2]{\kappa_{#1}(#2)}
\newcommand{\mmultiplier}[4]{\kappa_{#1, #2}(#3, #4)}
\newcommand{\normConst}[3]{\eta_{#1}({#2}, {#3})}
\newcommand{\pred}[1]{\widehat{Y}_{t + h}^{\text{(#1)}}}
\newcommand{\AROptPred}[3]{\widehat{Y}_{#1 + #2}(#3)}
\newcommand{\approxAROptPred}[3]{\widehat{Y}_{#1 + #2}(\widehat{#3})}

\def\Cov{\text{Cov}}
\def\series{\xi}

\def\spn{{\text{span}}}
\def\AR{\text{AR}}
\def\TP{\text{TP}}
\def\FP{\text{FP}}
\def\FN{\text{FN}}
\def\TN{\text{TN}}
\def\FPR{\text{FPR}}
\def\FNR{\text{FNR}}
\def\TPR{\text{TPR}}
\def\TSS{\text{TSS}}
\def\HSS{\text{HSS}}
\def\EDI{\text{EDI}}
\def\Z{{\mathbb Z}}
\def\Prec{\text{prec}}
\def\Range{\text{Range}}

\def\wh{\widehat }

\theoremstyle{plain}
\newtheorem{assumption}{Assumption}[section]
\newtheorem{example}{Example}[section]
\newtheorem{definition}{Definition}[section]
\newtheorem{problem}{Problem}[section]
\newtheorem{remark}{Remark}[section]
\theoremstyle{plain}
\newtheorem{lemma}{Lemma}[section]
\newtheorem{theorem}{Theorem}[section]
\newtheorem{proposition}{Proposition}[section]
\newtheorem{corollary}{Corollary}[section]

\DeclareMathOperator*{\argmin}{arg\,min}

\usepackage{chngcntr}
\counterwithin{table}{section}
\counterwithin{figure}{section}

\allowdisplaybreaks

\begin{document}
\title{On the optimal prediction of extreme events in heavy-tailed time series 
with applications to solar flare forecasting}

\author{Victor Verma,\ 
Stilian Stoev\thanks{Corresponding author,  {\tt sstoev@umich.edu};  Department of 
Statistics, 272 West Hall, 1085 South University Ave, Ann Arbor, Michigan 48109-1107},\ 
and Yang Chen\\   Department of Statistics, University of Michigan, Ann Arbor \\
}

\date{}

\maketitle

\thispagestyle{empty}

\begin{abstract}
The prediction of extreme events in time series is a fundamental problem arising in 
many financial, scientific, engineering, and other applications. We begin by establishing a 
general Neyman-Pearson-type characterization of optimal extreme event predictors in terms of density 
ratios.  This yields new insights and several closed-form optimal extreme event predictors 
for additive models.  These results naturally extend to time series, where we study optimal 
extreme event prediction for both light- and heavy-tailed autoregressive and moving average models.  Using a uniform law of large numbers for ergodic time series, we establish the asymptotic optimality 
of an empirical version of the optimal predictor for autoregressive models.  Using multivariate regular variation, we obtain an expression for the optimal extremal precision in heavy-tailed infinite moving averages, which provides theoretical bounds on the ability to predict extremes in this general class of models.  
We address the important problem of predicting solar flares by applying our theory and methodology to a state-of-the-art time series consisting of solar soft X-ray flux measurements. Our results demonstrate the success and limitations in solar flare forecasting of long-memory autoregressive models and long-range-dependent, heavy-tailed FARIMA models. 
\end{abstract}


\section{Introduction}\label{sect:intro}

Let $Y=\{Y_t,\ t\in \Z\}$ be a stationary time series.  The problem of predicting or 
forecasting a future value $Y_{t+h},\ h\ge 1$, of the time series based on 
its past observations $Y_s,\ s\le t$, is of fundamental theoretical and practical 
interest.  This problem arises naturally in many financial, scientific, and engineering 
applications, and it has been studied extensively in many articles and monographs \citep[][]{brockwell1991time, pourahmadi:2001, koul-book, box2015time, shumway2017time}.  Broadly speaking, 
the prediction problem amounts to finding a statistic 
$\hat Y_t(h) = g(Y_s,\ s\le t)$ that minimizes
a particular risk functional ${\cal R}(Y_{t+h},\hat Y_t(h))$. 
For finite-variance time series perhaps the most commonly used criterion 
has been the mean squared risk ${\cal R}(\xi,\hat\xi):= \E[ (\xi - \hat \xi)^2 ]$
\citep[][]{box2015time}.   In the case when $Y$ is a heavy-tailed infinite-variance time series, optimal predictors in the sense of least absolute deviation (LAD), i.e., 
${\cal R}(\xi,\hat \xi) = \E | \xi - \hat \xi|$, have been considered
\citep[][]{davis:resnick:1993}.
Risk functionals that are robust to outliers and extremes have also been employed 
for fitting models to time series data \citep[][]{koul-book,koul2010acla}.

Notably, Professor Hira Koul has pioneered and developed powerful weighted 
empirical process approaches applicable for the analysis of long-range dependent 
and/or heavy-tailed time series 
\citep[see, e.g., the monograph][and the references therein]{koul-book}.

In many scientific applications, however, one is interested in a particular form of binary prediction.  Namely, predicting {\em extreme events} $\{Y_{t+h} > y_0\}$ associated 
with the time series $Y$, where $y_0$ is a relatively large threshold.  An important example motivating our current work is the problem of forecasting {\em extreme solar flares}, which has long been a problem of interest in the scientific community \citep[see, e.g.,][and the references therein]{chen2023edit}. Solar flares are sudden, massive eruptions of electromagnetic radiation from the Sun's atmosphere 
\citep[][]{NASA_web_page,fletcher2011anob}. 
Predicting extreme solar flares is of great interest since they can have adverse effects like radio blackouts \citep[][]{swpc2023sola}, and they can be accompanied by coronal mass ejections (CMEs) and solar energetic particle events that themselves can have severe effects, such as electrical blackouts \citep[][]{bolduc2002gico} and irradiation of astronauts \citep[][]{mertens2019char}, respectively. Notably, in February 2022, a CME associated with a strong solar flare resulted in the destruction of 38 Starlink satellites \citep[][]{kataoka2022unex}.

\begin{figure}[t!]
    \centering
    \includegraphics[width=0.7\linewidth]{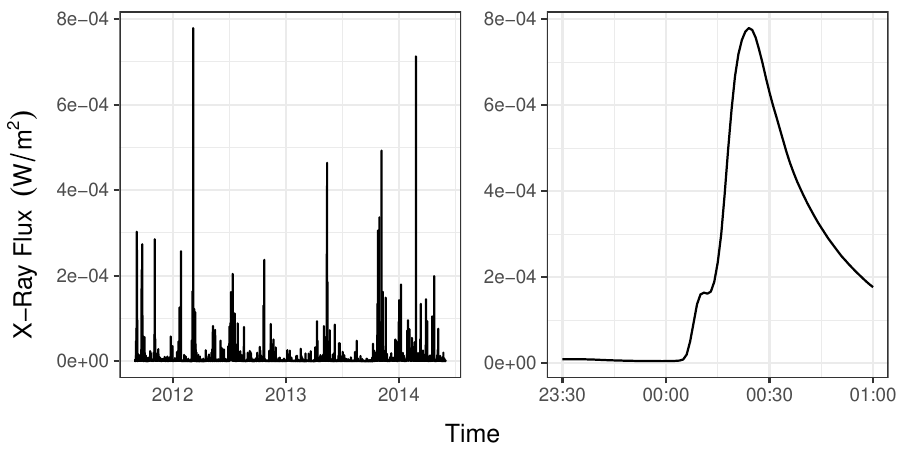}
    \caption{{\small \textbf{Left:} the GOES soft X-ray flux time series over the window 1 September 2011 to 31 May 2014. \textbf{Right:} the time series between 11:30 PM on 6 March 2014 and 1:00 AM on 7 March 2014, during which a strong solar flare occurred.}}
    \label{fig:goes_flux_time_series_nonlogged}
\end{figure}

The strength of a solar flare is defined as the peak value of the soft solar X-ray flux while it is ongoing \citep[][]{fletcher2011anob}.
The left panel of Figure~\ref{fig:goes_flux_time_series_nonlogged} displays a time series of flux measurements from the period 1 September 2011 to 31 May 2014 obtained from the Geostationary Operational Environmental Satellites (GOES) (cf. Section 
\ref{sect:data_analysis}); the right panel shows the same time series, but around the time of a strong solar flare that occurred on the night of 6–7 March 2014. The key features of the flux time series are heavy tails and complex temporal dependence \citep[][]{stanislavsky2009fari}.
The main objective is to provide timely forecasts of extreme flaring events 
$\{Y_{t+h}>y_0\}$, where $Y_{t + h}$ is the flux at a future time $t + h$ and $y_0$ is a high threshold.  This is one formulation of the so-called {\em operational flare forecasting} problem.  A plethora of modern machine learning tools have been used for operational flare forecasting, but their accuracy has been far from satisfactory; though they have had tremendous success in other prediction tasks, they have not managed to surpass more basic approaches, and the problem is still largely open \citep[][]{barnes2016acom, leka2019acomII, leka2019acomIII}. Note that our framing of the flare forecasting problem can be applied to other time series forecasting problems; for example, it can be applied to the problem of forecasting the solar radio flux, an important proxy for solar activity that has been widely studied \citep[][]{swpc2024f107, gromenko2017eval}.

In the flare forecasting context and similar contexts, a natural question 
arises as to what are the {\em optimal predictors} and the ultimate
{\em theoretical limits} of accuracy in predicting extreme events in time series, 
regardless of the methodology employed.
In this paper, we address the fundamental problem of {\em optimal prediction}
of extreme events $\{Y_{t+h}>y_0\}$ in time series.  Motivated by the solar flare
prediction problem, we focus on heavy-tailed time series models though many of 
our results apply in greater generality.  Somewhat 
surprisingly, this important binary prediction problem has been studied much less than other prediction problems in the literature.
Therefore, in Section \ref{sect:foundations}, we begin by outlining a framework for 
the optimal prediction of extreme events.  Briefly, 
suppose that the target event  occurs at a rate  $1-p = \P[Y_{t+h} > y_0]$, 
where $p\approx 1$.  Then, all predictors of these events can be written in the form
$\{\wh Y_t(h) > \wh y_0\}$, where $\wh Y_t(h) = g(Y_s,\ s\le t)$ is some statistic 
of the observed data.  We say that the predictor is {\em calibrated} if it raises alarms at a specified rate, i.e., 
$\P[\wh Y_t(h) > \wh y_0] =1-q$ for a specified $q$; typically, $q = p$. The {\em precision} 
of a calibrated predictor $\{\wh Y_t(h) > \wh y_0\}$ is
\begin{equation}\label{e:prec-intro}
    \Prec_{Y_{t+h} > y_0}(\wh Y_t(h) > \wh y_0) :=  \P[ Y_{t+h} > y_0 | \wh Y_t(h) > \wh y_0].
\end{equation}
Note that if the current time is $t$, then the event $\{Y_{t + h} > y_0\}$ is uncertain; it may or may not occur. If we raise an alarm at time $t$, we quantify our uncertainty about whether $Y_{t + h}$ will be extreme using the precision. We give a more general definition of our notion of precision in Definition~\ref{def:pred_properties}. It is the population analogue of the sample precision, a metric that is widely used in the evaluation of classifiers and is defined as the proportion of alarms that were followed by the event of interest \citep[][]{fawcett2006anin}. See Remark~\ref{rem:sample_metrics} for more detail on the motivation for the definition. A predictor will be said to be optimal if it {\em maximizes} the precision among all calibrated predictors. Thus, calibrated optimal predictors are ones that maximize the precision while controlling the alarm 
rate  (see also Sections  \ref{subsect:general_opt_preds} 
and \ref{subsect:roc_domin_prop} for optimality under other modes of calibration). 

Theorem \ref{thm:base_thm} establishes a general characterization of optimal
predictors by showing that $\wh Y_t(h)$ may be obtained as a ratio of
two probability densities.  This is akin to the celebrated Neyman-Pearson 
characterization of the most powerful tests in terms of likelihood ratios
and relates to a recent line of work on optimal binary classification 
\citep[][]{tong2016asur}. Building on this general result, in Section
\ref{subsect:explicit_opt_preds} we obtain closed-form characterizations of the optimal
event predictors in linear heteroskedastic models, leading to complete characterizations
for the Gaussian copula, linear heavy-tailed models, and a number of other contexts.

In Section \ref{subsect:opt_pred_asymp_aspects}, motivated by our main objective to
predict {\em extreme events}, we study the asymptotic aspects of optimal prediction
as the threshold $y_0$ grows, or equivalently, as $p\uparrow 1$. 
It turns out that the limiting precision in \eqref{e:prec-intro} is in fact the 
{\em upper tail dependence} coefficient between 
$Y_{t+h}$ and $\wh Y_t(h)$.  Tail dependence has been widely used 
and studied in extreme value theory 
\citep[see, e.g.,][and the references therein]{russell2016data,jansen:neblung:stoev:2023}.  However, to the best of our knowledge, our work is the first to explicitly 
characterize the {\em optimal extreme event} predictors, namely, the ones that lead to the 
strongest tail dependence (see Theorem \ref{thm:lambda-opt}, below).  This naturally leads us to establish new results on the optimal extremal precision for prediction of extreme events in heavy-tailed linear time series (cf. Section
\ref{subsect:MA_infty_mod_opt_pred_asymp_aspects}, particularly Theorem~\ref{thm:MA_infty_opt_asymp_precision}). These results
characterize the best possible asymptotic precision that can be achieved by any 
calibrated predictor of the target extreme events. Interestingly, for a large class of causal infinite moving average models with 
absolutely monotone sequence of coefficients, optimal extreme event predictors are simply achieved by using the most recent observation (Proposition \ref{prop:oracle}). In general, however, the optimal predictors for moving average models involve the complete history of the time series.  In the important special case of autoregressive models, these predictors can be computed with a finite sample (see Section \ref{subsect:explicit_lin_ts_mod_opt_preds}).  In Section \ref{subsect:AR_practical_pred}, we propose a practical methodology for the estimation of the optimal predictors in autoregressive time series.  Using results on uniform large numbers for empirical processes of ergodic time series \citep[][]{adams2012unif}, we establish the asymptotic calibration and optimality of the predictors. Table~\ref{tab:mod_summary} summarizes the results we present for the kinds of time series models we consider in this paper, infinite moving average and autoregressive models.

\begin{table}[!htb]
    \centering
    \begin{tabular}{|>{\centering\arraybackslash}p{2.7cm}|>{\centering\arraybackslash}p{2.2cm}|>{\centering\arraybackslash}p{2.2cm}|>{\centering\arraybackslash}p{2.5cm}|>{\centering\arraybackslash}p{3.2cm}|}
        \hline
        \textbf{Model} & \textbf{Oracle} & \textbf{Estimator} & \textbf{Inference} & \textbf{Extremal Oracle Precision} \\
        \hline
        AR($d$) & Theorem~\ref{thm:AR_opt_pred} & Section~\ref{subsubsect:AR_approx_opt_pred} & Theorem~\ref{thm:fixed-p} & Numerical \\
        \hline
        MA($\infty$) & Theorem~\ref{thm:MA_infty_opt_pred} & Not Included & Open Problem & Theorem~\ref{thm:MA_infty_opt_asymp_precision} \\
        \hline
        FARIMA$(0, d, 0)$ & Theorem~\ref{thm:MA_infty_opt_pred} & Not Included & Open Problem & Theorem~\ref{thm:MA_infty_opt_asymp_precision} \\
        \hline
    \end{tabular}
    \caption{A summary of the results we present for the models we consider in this paper. Since FARIMA$(0, d, 0)$ models are MA($\infty)$ models, the results given for the former are the same as those given for the latter. For those two kinds of models, estimators that approximate the optimal or oracle predictor can be obtained by truncating the oracle predictor; they are not included in this paper. The inferential properties of these estimators remain to be studied. The extremal precision of the oracle estimator for AR($d$) models can be obtained by computing the MA($\infty)$ representation and then applying Theorem~\ref{thm:MA_infty_opt_asymp_precision}.}
    \label{tab:mod_summary}
\end{table}

In Section \ref{sect:data_analysis}, we present an extensive 
application of the optimal extreme event prediction methodology to the 
state-of-the-art GOES flux time series data.  We model the data with both autoregressive and heavy-tailed long-range dependent (FARIMA) time series 
models, as advocated for in works such as \cite{burnecki2014algo}. 
The optimal extreme event predictors under these models provide new insights 
on the challenging solar flare prediction problem.  Specifically, they indicate 
that both heavy tails and long-range dependence should be accounted for when 
predicting extreme solar flares.  The resulting empirical precisions and true 
skill statistics (Table \ref{tab:data_analysis_results}) provide lower bounds on 
the optimal prediction accuracy one can hope to achieve when predicting solar 
flares using time series methods applied to historical GOES flux data. Operational 
flare forecasting, however, may be significantly improved 
through the use of covariate information such as images of the Sun at various wavelengths. 
Inference for optimal extreme event predictors in this context is a future 
direction of research, which we plan to pursue with the guidance of the general characterization 
results in Theorems \ref{thm:base_thm} and \ref{thm:lambda-opt}. \\

{\em The rest of the paper is structured as follows.} In Section 
\ref{subsect:general_opt_preds}, we give a framework and general 
characterization results for optimal event prediction.  
We study the optimal extremal precision in Section \ref{subsect:opt_pred_asymp_aspects}. Section \ref{subsect:explicit_lin_ts_mod_opt_preds} applies and extends these results to the important class of linear heavy-tailed time series models. Section \ref{subsect:AR_practical_pred} develops statistical inference for the optimal predictors in autoregressive time series, while Section \ref{subsect:MA_infty_mod_opt_pred_asymp_aspects} establishes general theoretical results
on the optimal extremal precision in heavy-tailed moving average models. Section \ref{sect:data_analysis} applies prediction methods to the solar flare prediction problem using GOES flux data from the most active periods of two recent solar cycles. Technical proofs and further details are given in the Appendix.

\section{Foundations of Optimal Extreme Event Prediction}\label{sect:foundations}

In this section, we present a framework for extreme 
event prediction.  We introduce the notion of calibrated predictors 
and establish a general characterization result for the optimal 
predictors in terms of density ratios (Section \ref{subsect:general_opt_preds}).  Examples of optimal predictors 
are presented in Section \ref{subsect:explicit_opt_preds}. 
The notion of {\em extremally optimal prediction} is discussed
in Section \ref{subsect:opt_pred_asymp_aspects}.

\subsection{A General Characterization of Optimal Predictors}\label{subsect:general_opt_preds}

Problems like those mentioned in Section~\ref{sect:intro} can be viewed as instances of a certain general statistical problem. Let $X$ be a random vector in $\R^d$ and $Y$ a random variable. We shall make the following assumption about the joint distribution of $X$ and $Y$ to simplify the proofs of our results:
\begin{assumption}\label{assump:joint_dens}
    $X$ and $Y$ have a joint density $f$ with respect to $\mu \times \nu$, where $\mu$ and $\nu$ are Borel measures on $\R^d$ and $\R$, respectively.
\end{assumption}
The general problem motivating our work in this paper can be stated as follows:
\begin{problem}[Optimal extreme event prediction]\label{prob:prob}
    For $p \in (0, 1)$ and $y_0 \in \R$ satisfying $\P(Y > y_0) = 1 - p$, predict the event $\{Y > y_0\}$ using $X$, 
    either for fixed $p$ or in the extreme regime $p \uparrow 1$.
\end{problem}

Predicting whether $Y > y_0$ is equivalent to predicting whether the indicator $\I(Y > y_0)$ equals one. Any predictor of $\I(Y > y_0)$ via $X$ can be expressed in the form $\I(h(X) > \tau)$ for a measurable function $h$ and a real number $\tau$; we predict that $Y > y_0$ if and only if this indicator equals one. Theorem~\ref{thm:base_thm} below characterizes the optimal $h$ for a class of predictors that are calibrated in a sense described in Definition~\ref{def:pred_properties}. In Section~\ref{subsect:roc_domin_prop}, we also show that this characterization is stronger in that it produces an optimal receiver operating characteristic (ROC) curve, and we discuss its connections to optimal Neyman-Pearson binary classification.

In the flare forecasting problem, $Y$ would represent the flux at a future time and $X$ would represent a vector of features engineered from current and past values of the flux as well as other relevant quantities that reflect the state of the Sun. The probability $p$ would correspond to a flare class threshold. The flares of interest would then occur with probability $1 - p$; we might want our predictor to be calibrated so that it raises an alarm with the same probability. Assuming that our predictor is calibrated in this way, or some other way, we would aim to optimize its precision, the probability that an alarm is raised given that a flare is imminent. Calibration, precision, and optimality are formally defined next.

\begin{definition}\label{def:pred_properties}
    Let $\I(h(X) > \tau)$ be a predictor of $\I(Y > y_0)$.

    (i) (\textit{Calibration}) The predictor $\I(h(X) > \tau)$ is said to be calibrated at level $q$ if 
    \begin{equation}\label{e:def:pred_properties:calibration}
    \P(h(X) > \tau) = 1 - q.
    \end{equation}
    The class of functions $h$ for which \eqref{e:def:pred_properties:calibration} holds for some
    $\tau$ will be denoted by ${\cal C}_q(X)$.

    (ii) (\textit{Precision}) The precision of the predictor $\I(h(X) > \tau)$ is defined as
     \begin{equation}\label{e:def:pred_properties:prec}
    \Prec_{Y > y_0}(h(X) > \tau) = \P(Y > y_0 \mid h(X) > \tau).
    \end{equation}
    
    (iii) (\textit{Optimality}) A predictor $\I(h(X) > \tau)$ that is calibrated at level $q$ is said to be optimal if for any other predictor $\I(k(X) > \rho)$ of $\I(Y > y_0)$ that is calibrated at level $q$, we have

    \begin{equation}\label{eq:optimality_cond}
        \P(Y > y_0 \mid h(X) > \tau) \ge \P(Y > y_0 \mid k(X) > \rho).
    \end{equation}
\end{definition}

For any random variable $\xi$, we denote by $F_{\xi}$ the distribution function of $\xi$,
and consider its left-continuous generalized inverse, defined as:
\begin{equation}\label{e:gen-inverse}
\Fi{\xi}(p) := \inf \{y : F_{\xi}(y) \ge p\},\ \ \ p\in (0,1).
\end{equation}
We will refer to $\Fi{\xi}(p)$ as to the $p$th quantile of $\xi$, for $p \in (0, 1)$.
We shall need the following result.
\begin{lemma}\label{lem:calibration_levels}
    For any random variable $\xi$, we have:
    
    {\em (i)} $F_\xi^\leftarrow(F_\xi(\xi)) = \xi$, almost surely.
    
    {\em (ii)} $F_\xi(F_\xi^\leftarrow(q))= q$, for all $q\in \Range(F_\xi)$, i.e., 
    \begin{equation}\label{e:exact-calibration}
    \P[ \xi > F_\xi^\leftarrow(q)] = 1-q,\ \ \mbox{ for all }q\in \Range(F_\xi).
    \end{equation}
\end{lemma}
\begin{proof}
    See Section~\ref{subsubsect:calibration_levels_lem_pf}.
\end{proof}

\begin{remark}  If the predictor $\I(h(X) > \tau)$ is calibrated 
at level $q$, i.e., \eqref{e:def:pred_properties:calibration} holds, then the
threshold $\tau$ therein can be chosen as $\tau:= F_{h(X)}^\leftarrow(q)$.
Indeed, having \eqref{e:def:pred_properties:calibration} for some $\tau$ is equivalent 
to having $q \in {\rm Range}(F_{h(X)})$ and hence by part (ii) of Lemma \ref{lem:calibration_levels}, we have
$
\P[ h(X) > F_{h(X)}^\leftarrow(q) ] = 1-q.
$
\end{remark}

\begin{remark}
    The probability of the event of interest, $\{Y > y_0\}$, is fixed at $1 - p$. It is natural to calibrate the predictor $\I(h(X) > \tau)$ at level $p$, assuming that $p \in \Range(F_{h(X)})$. Otherwise, an alarm for the event of interest would be raised more (or less) frequently than the event occurs. However, a practitioner may want to calibrate the predictor at a level different from $p$ to improve some performance measure, like the area under the ROC curve, at the expense of an alarm rate that is bigger or smaller than the rate of occurrence of the extreme event.
    For more details, see Section~\ref{subsect:roc_domin_prop}.
\end{remark}
\begin{remark}
    The notion of precision in Definition~\ref{def:pred_properties} is closely related to the notion of precision arising in the evaluation of classifiers. See \cite{fawcett2006anin} and Remark~\ref{rem:sample_metrics}.
\end{remark}

The following theorem characterizes the optimal predictor.
\begin{theorem}\label{thm:base_thm}
    Consider a threshold $y_0$ such that $p = \P(Y \le y_0) \in (0, 1)$. Let $r(x) = f_1(x) / f_0(x)$, where
    \[
    f_0(x) = \frac{1}{p}\int_{(-\infty, y_0]} f(x, y)\,\nu(dy)
    \quad \text{and} \quad
    f_1(x) = \frac{1}{1 - p}\int_{(y_0, \infty)} f(x, y)\,\nu(dy),
    \]
    and where $r(x)$ is interpreted as $\infty$ when $f_0(x) = 0$.
    For every $q \in \Range(F_{r(X)}) \cap (0, 1)$, the predictor $\I(r(X) > \tau)$ of $\I(Y > y_0)$ that is calibrated at level $q$ is optimal at that level, in the sense of Definition~\ref{def:pred_properties}, i.e., Relation~\eqref{eq:optimality_cond} holds.
\end{theorem}

\begin{remark}
    This theorem is similar to the characterization of the optimal predictor in the Neyman-Pearson paradigm for classification \citep[][]{tong2016asur}. The theorem characterizes optimal predictors as density ratios. Indeed, the functions $f_0$ and $f_1$ are the conditional densities of $X \mid Y \le y_0$ and $X \mid Y > y_0$, respectively. More insights and connections to the literature on binary classification are given in Section~\ref{subsect:roc_domin_prop}.
\end{remark}

\begin{proof}[Proof of Theorem~\ref{thm:base_thm}]
    Note that the density of $X$ with respect to $\mu$ can be written as
    \[
    f_0(x)\P(Y \le y_0) + f_1(x)\P(Y > y_0) = p f_0(x) + (1 - p)f_1(x).
    \]
    Let $S = r^{-1}((\tau, \infty))$. Then
    \[
    1 - q
    = \P(r(X) > \tau)
    = p\int_S f_0 + (1 - p)\int_S f_1,
    \]
    where here and in the rest of the proof we omit $\mu(dx)$ in the expressions of integrals for the sake of brevity. It follows that
    \begin{equation}\label{eq:int_f0_S}
    \int_S f_0 = \frac{1}{p}\left[1 - q - (1 - p)\int_S f_1\right].
    \end{equation}
    Let $\I(k(X) > \rho)$ be a predictor of $\I(Y > y_0)$ that is calibrated at level $q$ and let $T = k^{-1}((\rho, \infty))$. We must have
    \begin{equation}\label{eq:int_f0_T}
    \int_T f_0 = \frac{1}{p}\left[1 - q - (1 - p)\int_T f_1\right].
    \end{equation}
    We will show that it is not possible that:
    \begin{equation}\label{e:NP-counter-positive}
     \int_T f_0 < \int_S f_0.
    \end{equation} Suppose that this is true. It follows from \eqref{eq:int_f0_S} and \eqref{eq:int_f0_T} that $\int_T f_1 > \int_S f_1$. Also note that for $x \in T \setminus S$, $f_1(x) \le \tau f_0(x)$, and for
    $x\in S\setminus T$, $f_1(x) > \tau f_0(x)$. We have
    \begin{align}
        \int_T f_1
        & = \int_{T \cap S} f_1 + \int_{T \setminus S} f_1 \nonumber \\
        & = \int_{T \cap S} f_1 + \int_{T \setminus S} f_1 + \int_{S \setminus T} f_1 - \int_{S \setminus T} f_1 \nonumber \\
        & = \int_S f_1 +  \int_{T \setminus S} f_1 - \int_{S \setminus T} f_1  \le \int_S f_1 + \tau\int_{T \setminus S} f_0 - \tau\int_{S \setminus T} f_0 \nonumber \\
        & \le \int_S f_1 + \tau\int_{T \setminus S} f_0 - \tau\int_{S \setminus T} f_0 + \tau\left(\int_S f_0 - \int_T f_0\right) \label{e:by-assumption} \\
        & \le \int_S f_1 + \tau\int_{T \cup S} f_0 - \tau\int_{T \cup S} f_0 
        = \int_S f_1, \label{e:NP-impossible}
    \end{align}
    where Relation \eqref{e:by-assumption} follows from the assumed inequality \eqref{e:NP-counter-positive}.  The upper bound
    in Relation \eqref{e:NP-impossible}, however, contradicts \eqref{e:NP-counter-positive} in view of \eqref{eq:int_f0_S} and \eqref{eq:int_f0_T}.
    This shows that \eqref{e:NP-counter-positive} is impossible and hence, $\int_T f_0 \ge \int_S f_0$. 
    It follows that $\int_T f_1 \le \int_S f_1$. Both predictors are assumed to be calibrated at level $q$, so 
    $P(r(X) > \tau) = \P(k(X) > \rho) = 1 - q$. We have
    \begin{align*}
        \P(Y > y_0 \mid k(X) > \rho)
        &= \frac{\P(Y > y_0, k(X) > \rho)}{\P(k(X) > \rho)} 
        = \frac{\P(k(X) > \rho \mid Y > y_0)\P(Y > y_0)}{1 - q} \\
        &= \frac{1 - p}{1 - q}\int_T f_1 
        \le \frac{1 - p}{1 - q}\int_S f_1 
        \le \P(Y > y_0 \mid r(X) > \tau).
    \end{align*}
    This shows that $\I(r(X) > \tau)$ is optimal.
\end{proof}

\begin{remark}
    Observe that in general, the density ratio $r$ as well as $p$
    in Theorem~\ref{thm:base_thm} need to be estimated in practice. Also, 
    though several approaches for estimating density ratios have been developed \citep[see, e.g.,][]{sugiyama2012dens}, in applications, it may be difficult to estimate $r$, particularly when $p \approx 1$, when there would be few observations from $f_1$. In Section~\ref{subsect:explicit_opt_preds}, we discuss several broadly-applicable models for which optimal prediction can be performed without having to estimate $p$ or $r$.
    \comment{For example, when $Y = g(X) + \epsilon$, with $X \ind \epsilon$, as in 
    Theorem~\ref{thm:add_err_mod_thm}, then $\I(g(X) > \tau)$ is an optimal predictor of $\I(Y > y_0)$ for some $\tau \in \R$; this predictor depends on neither $p$ nor $r$.}
\end{remark}

\comment{\begin{remark}
    Note that $r(x)$ in Theorem~\ref{thm:base_thm} is only defined for $x$ belonging to the support of $f_0$.
\end{remark}}

\comment{
\begin{remark}
    When $q = p$, the precision is related to the notion of tail dependence. Suppose that we have a predictor $\I(h(X) > \tau)$ of $\I(Y > y_0)$, and assume that $h$ does not depend on $p$. If $p$ is in both $\Range(F_Y)$ and $\Range(F_{h(X)})$, then it follows from Lemma~\ref{lem:calibration_levels} that the predictor and response can be represented as $\I(h(X) > \Fi{h(X)}(p))$ and $\I(Y > \Fi{Y}(p))$, respectively. As $p \uparrow 1$, the precision $\Prec_{Y > \Fi{Y}(p)}(h(X) > \Fi{h(X)}(p))$ tends to the tail dependence coefficient between $Y$ and $h(X)$. For $p \approx 1$, the precision is very close to the tail dependence coefficient, which motivates the study of asymptotic optimality, where a predictor is optimal if it maximizes the tail dependence coefficient. In Section~\ref{subsect:opt_pred_asymp_aspects}, we discuss asymptotically optimal predictors of extreme events.
\end{remark}}

\begin{remark}
Consider a target event $A = \{Y > y_0\}$ that occurs with probability $1 - p$ and a
predictor based on $\wh A = \{h(X)>\tau\}$ that is calibrated at level $q$, where
$p,q\in (0,1)$. It is easy to see that if $p$ and $q$ are fixed, then maximizing the 
precision 
$\P[ A | \wh A]$ of the predictor is equivalent to minimizing the Hamming loss
$\E | \I(A) - \I({\wh A})|$.  Indeed, this follows from the observation that 
$
\E| \I(A) - \I({\wh A})| = \E ( \I(A) + \I({\wh A}) - 2 \I(A)\I({\wh A}))
= (1-p) + (1-q) - 2 (1-q) \P[A | \wh A].
$
\end{remark}

\comment{\begin{remark}
    Given a function $h$ and a covariate vector $X$, we can construct a family of predictors $\{\I(h(X) > \tau)\}$ by varying the threshold $\tau$, and thus the calibration level $q$. It turns out that if $h$ is the density ratio $r$ in Theorem~\ref{thm:base_thm}, then the family $\{\I(h(X) > \tau)\}$ has an optimal ROC curve (cf Section~\ref{subsect:roc_domin_prop}).
\end{remark}}

\begin{remark}
In applications, especially when predicting extreme events $(p\approx 1)$, we often 
consider the {\em balanced case} where $p = q$. This requirement is natural as it mandates that alarms be raised only as frequently as events occur. However, a practitioner may want to construct predictors for varying levels of the alarm rate $1 - q$. This naturally leads to the notion of ROC curve optimality of families of predictors that differ in their alarm rate. We discuss this notion in Section~\ref{subsect:roc_domin_prop}.
\end{remark}

\subsection{Explicit Optimal Predictors}\label{subsect:explicit_opt_preds}

In this section, we consider several models where a closed-form 
expression for the optimal predictor exists that does not involve the 
density ratio in Theorem~\ref{thm:base_thm}. These results will be used 
in the sequel to characterize optimal extreme event predictors for 
moving average and autoregressive models (Section~\ref{sect:lin_ts_mod_opt_preds}). Extensions to state space models are briefly discussed Section \ref{subsect:more_example_opt_preds}).  They are also of independent interest.

In the sequel all optimal predictors are understood in the sense of Definition \ref{def:pred_properties}.  Also, 
a function will be referred to as increasing (decreasing, respectively) if it is monotone non-decreasing 
(non-increasing, respectively).  When strict monotonicity is required, we will use the terms strictly 
increasing (decreasing, respectively). The lemma below will be helpful in the proofs of several characterization results.

\begin{lemma}\label{lem:opt_pred_cond_simplifier}
    Let $X$ be a random variable, let $h : \R \to \R$ be a (Borel) measurable function, and let $q \in (0, 1)$ and $\tau \in \R$ such that $\P(h(X) > \tau) = 1 - q$.
    
    (i) If $h$ is increasing and left-continuous, then $x_0 = \sup\{x \in \R : h(x) \le \tau\}$ is finite, and $h(X) > \tau$ if and only if $X > x_0$.

    (ii) If $h$ is decreasing and right-continuous, then $x_0 = \inf\{x \in \R : h(x) \le \tau\}$ is finite, and $h(X) > \tau$ if and only if $X < x_0$.
\end{lemma}
\begin{proof}
    (i) Let $S = \{x \in \R : h(x) \le \tau\}$. This set must be nonempty, as $S$ being empty would imply that $q = \P(h(X) \le \tau) = 0$. It must also be bounded above. If it were not, then $h(x) \le \tau$ for arbitrarily large $x$; since $h$ is increasing, it would follow that $h(x) \le \tau$ for all $x$, and we would get that $q = \P(h(X) \le \tau) = 1$. Since $S$ is nonempty and bounded above, $x_0 = \sup S \in \R$.

    Let $\{x_n\}$ be a sequence in $S$ such that $x_n \to x_0$ as $n \to \infty$. Since $h$ is left-continuous, letting $n$ tend to infinity in $h(x_n) \le \tau$ yields $h(x_0) \le \tau$. Hence, $h(X) > \tau$ if and only if $X > x_0$.

    (ii) The proof of this part is essentially the same as the proof of part (i).
\end{proof}

\begin{theorem}[Additive heteroskedastic models]\label{thm:add_err_mod_thm}
    Let $X = (X_i)_{i=1}^d$ and $\epsilon$ be independent and let
    \begin{equation}\label{e:thm:add_err_mod_thm}
    Y = g(X) + \sigma(X)\epsilon,
    \end{equation}
    for some deterministic measurable functions $g(\cdot)$ and $\sigma(\cdot)$, where $\sigma(\cdot)$ is positive on the range of $X$.
    Let $X$ have a density $f_X$ with respect to a $\sigma$-finite Borel measure $\mu$ on $\R^d$, and let $\epsilon$ have a density $f_{\epsilon}$ with respect to the Lebesgue measure on $\R$. Let $p \in (0, 1)$ and $y_0 \in \R$ 
    be such that $\P(Y \le y_0) = p$.  
    Then an optimal predictor of $\I(Y > y_0)$ is $\I(k(X) > \tau)$, where
    \[
    k(X) = \frac{g(X) - y_0}{\sigma(X)}
    \]
    and $\tau$ is some suitably chosen constant that calibrates the predictor.
\end{theorem}

\begin{remark}\label{rem:homoscedastic_noise}
    When the noise in \eqref{e:thm:add_err_mod_thm} is homoskedastic, i.e., when its variance is constant with respect to $X$, the optimal predictor can be written as
    $\I(g(X) > \tau)$
    for some constant $\tau$ that calibrates the predictor at the desired level.
\end{remark}

\comment{\begin{remark}
    Theorem~\ref{thm:add_err_mod_thm} remains valid if $X$ is a random element taking values in 
    a complete separable metric space.
\end{remark}}

\begin{proof}[Proof of Theorem \ref{thm:add_err_mod_thm}] Denote the Lebesgue measure on $\R$ by $\nu$ and 
    note that we can write the joint density $f$ of $X$ and $Y$ with respect to $(\mu\times\nu)(dx,dy) = 
    \mu(dx)\nu(dy)$ as
    \begin{equation}\label{e:thm:add_err_mod_thm-density}
    f(x, y)
    = f_X(x)f_{\epsilon}\left(\frac{y - g(x)}{\sigma(x)}\right)\frac{1}{\sigma(x)}.
    \end{equation}
    Indeed, for all Borel $A\subset \R^d$ and $B\subset \R$, using the independence of $X$ and 
    $\epsilon$, we obtain
    \begin{align}
    \P[ X\in A, Y\in B] &= \E\Big[ \I(X\in A) \E  \Big[ \I( \epsilon \in (B - g(X))/\sigma(X) | X \Big] \Big]
    \nonumber\\
    & = \int_A \Big(f_X(x) \int_{(B-g(x))/\sigma(x)} f_\epsilon(u) \nu(du)\Big)  \mu(dx) \nonumber\\
    & = \int_A \Big( \int_{B} f_X(x) f_\epsilon\Big(\frac{y-g(x)}{\sigma(x)} \Big) \frac{1}{\sigma(x)} 
    \nu(dy)\Big)\mu(dx) \label{e:thm:add_err_mod_thm-change}\\
    &= \int_{A\times B} f(x,y) (\mu\times \nu)(dx,dy), \nonumber
    \end{align}
    where in \eqref{e:thm:add_err_mod_thm-change} we used the change of variables $y = g(x) + \sigma(x)u$ and the homogeneity and translation-invariance of the Lebesgue measure $\nu$. The last relation follows from 
    Tonelli's theorem and since the sets $A\subset \R^d$ and $B\subset\R$ were arbitrary, it
    shows that $f(x,y)$ in \eqref{e:thm:add_err_mod_thm-density} is the joint density of $X$ and $Y$
    with respect to $(\mu \times \nu)(dx,dy)$.

    In light of the above derivation, the functions $f_0$ and $f_1$ defined below are the conditional densities of $X$ given $Y \le y_0$ and $Y > y_0$, respectively:
    \begin{align*}
        f_0(x)
        &:= \frac{f_X(x)}{p\sigma(x)} \int_{(-\infty, y_0]} f_{\epsilon}\left(\frac{y - g(x)}{\sigma(x)}\right)\,\nu(dy)
        = \frac{f_X(x)}{p} F_\epsilon\left(\frac{y_0 - g(x)}{\sigma(x)}\right) \\
        f_1(x)
        &:= \frac{f_X(x)}{(1 - p)\sigma(x)} \int_{(y_0, \infty)} f_{\epsilon}\left(\frac{y - g(x)}{\sigma(x)}\right)\,\nu(dy)
        = \frac{f_X(x)}{1 - p} \overline F_\epsilon\left(\frac{y_0 - g(x)}{\sigma(x)}\right),
    \end{align*}
    where $\overline F_\epsilon(u) = 1 - F_\epsilon(u) = \int_{(u, \infty)} f_\epsilon(v)\,\nu(dv)$ is the tail function of $\epsilon$.
    
    The ratio $r(x)$ of the conditional densities is thus
    \begin{equation}\label{eq:add_err_mod_thm_density_ratio}
        r(x)
        = \frac{f_1(x)}{f_0(x)}
        = \left(\frac{p}{1 - p}\right)\frac{\overline{ F}_{\epsilon}\left(\frac{y_0 - g(x)}{\sigma(x)}\right)}{F_{\epsilon}\left(\frac{y_0 - g(x)}{\sigma(x)}\right)}.
    \end{equation}
    For any $q \in \Range(F_{r(X)}) \cap (0, 1)$, an optimal predictor of $\I(Y > y_0)$ at level $q$ is $\I(r(X) > \Fi{r(X)}(q))$, by Theorem~\ref{thm:base_thm} and Lemma~\ref{lem:calibration_levels}. Define $\tilde{r} : \R \to \R$ by
    \[
    \tilde{r}(u) = \left(\frac{p}{1 - p}\right)\frac{\overline{F}_{\epsilon}(u)}{F_{\epsilon}(u)}.
    \]
    By the properties of distribution functions, $\tilde{r}$ is decreasing and right-continuous. If we let $U = (y_0 - g(X)) / \sigma(X)$, then $\P(\tilde{r}(U) > \Fi{r(X)}(q)) = 1 - q$. It follows from part (ii) of Lemma~\ref{lem:opt_pred_cond_simplifier} that $\tilde{r}(U) > \Fi{r(X)}(q)$ if and only if $U < \rho$ for some $\rho \in \R$. Hence, the optimal prediction condition $r(X) > \Fi{r(X)}(q)$, which can also be written as $\tilde{r}(U) > \Fi{r(X)}(q)$, can be reduced to the simpler condition
    \[
    \frac{y_0 - g(X)}{\sigma(X)} < \rho.
    \]
    Therefore, an optimal predictor of $\I(Y > y_0)$ at level $q$ is $\I(k(X) > \tau)$, where
    \[
    k(x) = \frac{g(x) - y_0}{\sigma(x)}\ \ \mbox{ and }\ \ \tau = -\rho.
    \]
\end{proof}

\begin{remark} The assumption that $\epsilon$ in Theorem~\ref{thm:add_err_mod_thm} has 
a density with respect to the Lebesgue measure cannot be easily relaxed.  Indeed, in 
\eqref{e:thm:add_err_mod_thm-change} we used the fact that $\nu$ is translation invariant and homogeneous.  
By the existence and uniqueness theorem for Haar measures, all locally finite translation-invariant 
Borel measures are necessarily proportional to the Lebesgue measure. See also Theorem 0.1 in 
\cite{elekes2006isle} for a more general result. For this reason, we just assume that $\epsilon$ has a 
density with respect to that measure.
\end{remark}

\begin{remark}\label{rem:q-calibration}
    If $r(X)$ in \eqref{eq:add_err_mod_thm_density_ratio} has a continuous distribution function, then $\Range(F_{r(X)}) \supseteq (0, 1)$, so an optimal predictor can be calibrated at any level $q \in (0, 1)$.
\end{remark}

\begin{remark}
    There are additional models for which the optimal extreme event predictor simplifies nicely. Another is the model $Y = \varphi_{g(X)}(\epsilon)$, where $\varphi$ is a function that satisfies certain smoothness and monotonicity conditions, $g : \R \to \R$ is some function, and $X$ and $\epsilon$ are independent random variables. See Proposition~\ref{prop:non-linear}. Yet another is the model $Y = g(X + \delta) + \epsilon$, where $g : \R \to \R$ is increasing and $X$, $\delta$, and $\epsilon$ are independent random variables; see Proposition~\ref{prop:monotone-link}.
\end{remark}


We conclude this section with several examples. 

\begin{example}\label{ex:Gauss} Let $(Y,X)$ be a Gaussian random vector 
in $\R^{d+1}$, with $X = (X_i)_{i = 1}^d$.  Then, for some constant vector $a = (a_i)_{i=1}^d$, one can always write
$$
 Y = a^\top X + \epsilon,
$$
where $X$ and $\epsilon$ are independent.  Theorem \ref{thm:add_err_mod_thm} therefore implies that for all $p,q\in (0,1)$ an optimal level-$q$ predictor of $\I(Y> \Fi{Y}(p))$ is of the form $\I( h(X) > \Fi{h(X)}(q))$, where $h$ is the linear function $h(x) = a^\top x = \sum_{i=1}^d a_i x_i$. Note that by Remark \ref{rem:q-calibration}, calibration is possible at level $q$ for all $q \in (0, 1)$.
\end{example}

The previous example naturally yields a characterization of the optimal extreme event
predictors in Gaussian copulas.

\begin{example}[Optimal event prediction in Gaussian copulas] 
\label{ex:Gaussian-copula}
Suppose now that $Y$ and $X_i,\ i=1,\cdots,d$ are random variables with 
{\em continuous marginal} distributions.  Assume that
\begin{equation}\label{e:Y_Xi_via_V_Ui}
Y \stackrel{a.s.}{=} g(V)\ \ \mbox{ and }\ \ X_i \stackrel{a.s.}{=}f_i(U_i),\ i=1,\cdots,d,
\end{equation}
where $g$ as well as the $f_i$'s are {\em strictly increasing functions}.
Since $V = g^{-1}(Y)$ and $U_i = f_i^{-1}(X_i)$, almost surely, and since $g$ is strictly
increasing, the prediction of $\{Y>F_Y^\leftarrow(p)\}$ in terms of
$X=(X_i)_{i=1}^d$ is equivalent to the prediction of $\{V>F_V^\leftarrow(p)\}$ in terms 
of $U=(U_i)_{i=1}^d$. Example \ref{ex:Gauss} then implies the following.

\begin{corollary} Assume that $(Y,X_1,\cdots,X_d)$ has continuous marginals and
a Gaussian copula. Then, for all $p,q\in (0,1)$, an optimal predictor of 
$\I(Y> F_Y^\leftarrow(p))$ calibrated at level $q$ is given by $\I(h(X)> F_{h(X)}^\leftarrow(q))$,
where
$
h(x) = \sum_{i=1}^d a_i \Phi^{-1}(F_{X_i}(x_i)),
$
and where $\Phi$ denotes the standard normal distribution function.
\end{corollary}
\begin{proof} Let $g(v) := F_Y^\leftarrow \circ \Phi(v)$ and 
$f_i(u):= F_{X_i}^\leftarrow \circ \Phi(u).$  Since the CDFs of $Y$ and the
$X_i$'s are continuous, $g$ and the $f_i$'s are strictly 
increasing.  Moreover, letting $V:= \Phi^{-1}(F_Y(Y))$ and 
$U_i:= \Phi^{-1}(F_{X_i}(X_i))$, by part (i) of Lemma \ref{lem:calibration_levels}, we
obtain that \eqref{e:Y_Xi_via_V_Ui} holds.  This implies that the optimal 
prediction of $\I(Y> F_Y^\leftarrow(p))$ in
terms of $X$ is equivalent to the optimal prediction of 
$\I(V> F_V^\leftarrow(p))$ in terms of $U$ and Example \ref{ex:Gauss} then yields the claim.
\end{proof}
\end{example}


\subsection{Asymptotic Aspects of Optimal Prediction: Extremal Precision}\label{subsect:opt_pred_asymp_aspects}

In this section, we consider the scenario where the target event 
$\{Y> F_Y^{\leftarrow}(p)\}$ is extreme, i.e., $p\approx 1$. 
To make this precise, we consider predictors $\I(h(X) >\tau)$ that are calibrated at 
level $q=p$ and investigate their {\em asymptotic precision} as $p\uparrow 1$.  In this case, it turns out 
that the asymptotic precision (when it exists) coincides with the notion of 
{\em tail dependence} coefficient, widely studied in extreme value theory 
\citep[see, e.g.,][and the references therein]{coles1999depe,jansen:neblung:stoev:2023}.

\begin{definition}[Tail dependence]\label{def:tdc}  For two random variables $\xi$ and $\eta$, we write
\begin{equation}\label{e:lambda-xi-eta}
 \lambda(\xi,\eta):= \lim_{p\uparrow 1} \P[ \xi > F_\xi^{\leftarrow}(p) | \eta > F_\eta^\leftarrow(p)],
\end{equation}
whenever the limit exists.  In this case, $\lambda(\xi,\eta)$ is referred to as the (upper) 
tail-dependence coefficient for $\xi$ and $\eta$. 
\end{definition}

One can interpret $\lambda(\xi,\eta)$ as the limiting probability that $\xi$ is extreme, given that $\eta$ is extreme.  Notice
that $\lambda(\xi,\eta)\in [0,1]$.  If $\lambda(\xi,\eta) = 0$, then $\xi$ and $\eta$ are said to be asymptotically or extremally 
independent.  

From the perspective of extreme event prediction, the tail-dependence coefficient may be interpreted as {\em asymptotic precision}.
Indeed, suppose for simplicity that $\xi$ and $\eta$ have continuous distributions, though this requirement can be relaxed substantially. 
Then, in view of Lemma \ref{lem:calibration_levels} $F_\xi(F_\xi^{\leftarrow}(p)) = p = F_\eta(F_\eta^\leftarrow(p))$, 
for all $p\in (0,1)$, and hence $\I(\eta> F_\eta^\leftarrow(p))$ is a calibrated predictor 
of $\I(\xi> F_\xi^{\leftarrow}(p))$, at level $p\in (0,1)$ (Definition \ref{def:pred_properties}).  This, in view of \eqref{e:def:pred_properties:prec}, implies 
that the tail-dependence coefficient $\lambda(\xi,\eta)$ in \eqref{e:lambda-xi-eta} is the limiting, asymptotic precision for predicting 
the extremes of $\xi$ in terms of the extreme values of $\eta$.  In general, however, the extremes of $\eta$ need not be optimal for
predicting the extreme values of $\xi$.  This observation motivates the following definition.

\begin{definition}[Extremal optimality] \label{def:asymp-opt}  Let $Y$ be a random variable and $X$ 
a random vector.\\

{\em (i)} The optimal {\em extremal precision} in predicting the extremes of $Y$ in terms of $X$, 
is defined as:
$$
\lambda^{\rm (opt)}(Y,X):= \limsup_{p\uparrow 1} \sup_{g \in {\cal C}_p(X) } 
    \Big\{ \Prec_{Y > F_Y^\leftarrow(p)}(g(X) > F_{g(X)}^\leftarrow(p))\Big\},
$$
where the supremum is over the class of all $p$-calibrated predictors of the form 
$\I(g(X)> F_{g(X)}^\leftarrow(p))$, i.e., over all predictors based on Borel functions $g$ such that 
$\P[ g(X) > F_{g(X)}^\leftarrow(p)]=1-p$, and $\Prec$ is the prediction 
precision defined in \eqref{e:def:pred_properties:prec}.\\

{\em (ii)} A random variable $h(X)$, for some (Borel) measurable $h$ is said to be an 
{\em  optimal extremal predictor}, if
\begin{equation}\label{e:def:asym-opt-ii}
\lambda^{\rm (opt)}(Y,X) = \limsup_{p\uparrow 1}\Big\{ \Prec_{Y>F_Y^\leftarrow(p)}( h(X)> F_{h(X)}^\leftarrow(p)) \Big\}.
\end{equation}
In particular, if the tail dependence coefficient $\lambda(Y,h(X))$ exists and $h$ is an optimal predictor, we have
$\lambda^{\rm (opt)}(Y,X) = \lambda(Y,h(X))$.
\end{definition}

\medskip
The following result characterizes the optimal extremal precision and the corresponding extremal 
predictor in a simple but important setting of additive error models.  Similar results can be obtained
for the other cases where closed-form expressions for the optimal predictors are available (cf Section \ref{subsect:explicit_opt_preds}).

\begin{theorem}\label{thm:lambda-opt} Suppose that $X$ and $\epsilon$ are independent and let
$$
Y =  h(X)+ \epsilon,
$$
for some Borel function $h$. Assume that $\epsilon$ has a density with respect to the Lebesgue measure in
$\R$.  Let $u_{h(X)} = \sup\{ u\, :\, F_{h(X)}(u)<1\}$ and suppose that 
the CDF of $h(X)$ is continuous on the interval $(u, u_{h(X)}]$, for some $u<u_{h(X)}$.

We have that
 $$
 \lambda^{\rm (opt)}(Y,X) = \limsup_{p\uparrow 1} \P[ Y>F_Y^\leftarrow(p)\, |\, h(X) > F_{h(X)}^\leftarrow(p) ].
 $$
 and thus $h(X)$ is an asymptotically optimal predictor for the extreme events of $Y$ 
(Definition \ref{def:asymp-opt}).  In particular, if the tail-dependence coefficient 
$\lambda(Y,h(X))$ exists, then 
$$
\lambda^{\rm (opt)}(Y,X) = \lambda(Y,h(X)).
$$
\end{theorem}

\begin{proof} 
By assumption, we have the CDF $F_{h(X)}$ is continuous on $(u,u_{h(X)}]$.  
This implies that $(p_0,1)\subset F_{h(X)}(\R)$, for some $0<p_0<1$.  Therefore,
the predictor $h(X)$ can be calibrated at all levels $p\in (p_0,1)$.  That is, 
$h\in {\cal C}_{p}(X),\ p\in (p_0,1)$. Also, by Theorem \ref{thm:add_err_mod_thm}, we have that for every 
other level $p$-calibrated predictor $g\in {\cal C}_{p}(X)$,
$$
\P[ Y> F_Y^\leftarrow(p) \, |\, h(X)> F_{h(X)}^\leftarrow(p)] \ge \P[ Y> F_Y^\leftarrow(p) \, |\, g(X)> F_{g(X)}^\leftarrow(p)] = \Prec_{Y > F_Y^\leftarrow(p)}(g(X) > F_{g(X)}^\leftarrow(p)).
$$
Therefore, for all $p\in (p_0,1)$,
$$
\Prec_{Y > F_Y^\leftarrow(p)}(h(X) > F_{h(X)}^\leftarrow(p)) \ge \sup_{g\in {\cal C}_p(X)} \Prec_{Y > F_Y^\leftarrow(p)}(g(X) > F_{h(X)}^\leftarrow(p)),
$$
which implies \eqref{e:def:asym-opt-ii} and proves the desired optimality.
\end{proof}

Theorem \ref{thm:lambda-opt} enables us to characterize the {\em optimal extremal precision}
$\lambda^{\rm (opt)}(Y,X)$ in a number of cases.  
The first is a well-known result in the folklore showing that Gaussian models have 
vanishing prediction precision at the extremes.

\begin{corollary}\label{cor:Gaussian-extr-precision}  Suppose that $(Y,X)$ is a random vector in $\R^{d+1}$ with continuous marginals
and Gaussian copula having a non-singular covariance matrix.  Then, $\lambda^{\rm (opt)} (Y,X) = 0$.
\end{corollary}
\begin{proof} In view of Example \ref{ex:Gaussian-copula}, without loss of generality we may assume that
$(Y,X)$ is Gaussian. Recalling Example \ref{ex:Gauss}, we can write 
$
Y = h(X) + \epsilon, 
$
where $h(X) := a^\top X$ and $\epsilon$ are independent Gaussian random variables and hence
$h(X)= a^\top X$ yields the optimal predictor for the events $\{Y>F_Y^\leftarrow(p)\}$
(Theorem \ref{thm:add_err_mod_thm}). Thus, by Theorem \ref{thm:lambda-opt}, $\lambda^{(opt)}(Y,X) = \lambda(Y, a^\top X)$. However, since $(Y,X)$ has a non-singular covariance matrix, ${\rm Corr}(Y, a^\top X) < 1$ and the tail-dependence coefficient $\lambda(Y,a^\top X) = 0$.
See, for example, Relation (5.1) in \cite{ledford:tawn:1996}
or Table 4.1 in \cite{joe:2015}.
\end{proof}

\begin{example}[Heavy-tailed linear models]\label{ex:linear-Pareto}
Suppose that
$$
Y = a^\top X + \epsilon,
$$
where $X_1,\cdots,X_d$ and $\epsilon$ are independent and identically distributed standard 
$\alpha$-Pareto random variables, where $\alpha>0$.  That is, 
$$
\P[X_i>x] = \P[\epsilon>x] = \frac{1}{x^\alpha},\ x>1.
$$
Then, by Theorems \ref{thm:add_err_mod_thm} and \ref{thm:lambda-opt}, we have
$$
\lambda^{\rm (opt)}(Y,X) = \lambda( a^\top X + \epsilon, a^\top X)
$$
Thus, assuming that $\sum_{i=1}^d (a_i)_+^\alpha >0$, i.e., $a_i>0$ for some $i$, 
Lemma \ref{lem:two_series_tdc} below implies that
$$
\lambda^{\rm (opt)}(Y,X) = 
\frac{\sum_{i=1}^d (a_i)_+^\alpha}{1+ \sum_{i=1}^d (a_i)_+^\alpha},
$$
where $(x)_+^\alpha = (\max\{0,x\})^\alpha$. Note that the $X_i$'s and $\epsilon$ have extremal skewness one (see Definition~\ref{def:rv_rand_var2}).
\end{example}

Corollary \ref{cor:Gaussian-extr-precision} and Example \ref{ex:linear-Pareto} illustrate two
fundamentally different regimes for linear models. Namely, linear models with 
light-tailed (Gaussian) errors can often have vanishing optimal extremal precision.  
On the other hand, heavy-tailed linear models can exhibit positive and often quite 
large optimal extremal precisions depending on the sizes of the coefficients $a_i$.  This naturally makes the class of linear heavy-tailed time series a natural and flexible 
class of models for predicting extreme events.

\section{Optimal Extreme Event Prediction in Linear Time Series}\label{sect:lin_ts_mod_opt_preds}

In this section, we consider optimal extreme event prediction for linear time series models that can be written as a causal moving average of iid innovations, i.e., MA$(\infty)$ models. We give the optimal predictor for these models in Theorem~\ref{thm:MA_infty_opt_pred}. When the innovations are heavy-tailed, we can derive the asymptotic precision of the optimal predictor as the event rate goes to one. This is given in Theorem~\ref{thm:MA_infty_opt_asymp_precision}, which relies on the concept of regular variation (reviewed in Section~\ref{subsect:rv_info}). Since, in the general case, the optimal predictor given by Theorem~\ref{thm:MA_infty_opt_pred} thresholds a moving average of infinitely many innovations, it cannot in general be estimated from data. However, in the special case of autoregressive models, the optimal predictor thresholds a linear combination of finitely many observations (Theorem~\ref{thm:AR_opt_pred}), so it can be estimated from data. In Section~\ref{subsubsect:AR_approx_opt_pred}, we describe a full methodology for predicting extreme events in autoregressive time series and prove that the predictor the methodology is based on has desirable asymptotic properties as the sample size goes to infinity (Theorem~\ref{thm:fixed-p}).

\subsection{Explicit Optimal Predictors}\label{subsect:explicit_lin_ts_mod_opt_preds}

We shall give a closed-form expression for the optimal predictor for both MA($\infty$) models and autoregressive models in this section. The simplification of Theorem~\ref{thm:add_err_mod_thm} described in Remark~\ref{rem:homoscedastic_noise} will be used to do this; in order to make that result applicable, we will need the following assumption.
\begin{assumption}\label{assump:innov_density}
    The innovations $\epsilon_t, t \in \Z$, are iid random variables that have a density with respect to the Lebesgue measure, and for some
    $\delta>0$, $\E[ |\epsilon_t|^\delta ] <\infty$ for all $t$.
\end{assumption}
We first handle MA($\infty$) models. Recall that a process $\{Y_t\}_{t = -\infty}^{\infty}$ is said
to be invertible if there exists a real sequence $\{b_j\}_{j = 0}^{\infty}$ such that for all $t$,
\[
\epsilon_t = \sum_{j = 0}^{\infty} b_j Y_{t - j},\ \ \mbox{ almost surely.}
\]
\begin{theorem}\label{thm:MA_infty_opt_pred}
    Let $\{a_j\}_{j = 0}^{\infty}$ be a real sequence, let $\{\epsilon_t\}_{t = -\infty}^{\infty}$ be a sequence of random variables satisfying Assumption~\ref{assump:innov_density}, and define a process $\{Y_t\}$ by
    \begin{equation}\label{eq:Y_t_def}
        Y_t := \sum_{j = 0}^{\infty} a_j\epsilon_{t - j}, \quad t \in \Z.
    \end{equation}
    Suppose that the series in \eqref{eq:Y_t_def} converges almost surely for all $t$ and that the process $\{Y_t\}$ is invertible. For $h \ge 1$, set
    \begin{equation}\label{eq:MA_infty_opt_pred}
    \pred{opt} := \sum_{j = 0}^{\infty} a_{j + h}\epsilon_{t - j}
    \end{equation}
    and assume that $\pred{opt}$ is not trivially zero. Let $p \in (0, 1)$ 
    and $y_0 \in \R$ such that $\P(Y \le y_0) = p$. The optimal predictor of $\I(Y_{t + h} > y_0)$ based on $Y_s$, $s \le t$, is $\I(\pred{opt} > \tau)$, $\tau$ being a constant that calibrates the 
    predictor.
\end{theorem}
\begin{proof}
    In view of \eqref{eq:Y_t_def} and  \eqref{eq:MA_infty_opt_pred}, we can write
    \begin{equation}\label{eq:future_Y_decomp}
        Y_{t + h}
        = \sum_{j = 0}^{\infty} a_j\epsilon_{t + h - j}
        = \sum_{j = 0}^{\infty} a_{j+h}\epsilon_{t  - j} + \sum_{j = 0}^{h - 1} a_j\epsilon_{t + h - j}
        = \pred{opt} + \sum_{j = 0}^{h - 1} a_j\epsilon_{t + h - j}.
    \end{equation}
    The random variable $\pred{opt}$ only depends on $\epsilon_s$ for $s \le t$; since $\{Y_t\}$ is invertible, this means that $\pred{opt}$ only depends on $Y_s$ for $s \le t$. It is independent of the 
    last sum as the last sum only depends on $\epsilon_{t + 1}, \ldots, \epsilon_{t + h}$. Furthermore, $\pred{opt}$ as well as the last sum have densities with respect to the Lebesgue measure by Assumption~\ref{assump:innov_density}. Taking $X$ to be $\pred{opt}$ and $g$ to be the identity function in Remark~\ref{rem:homoscedastic_noise}, we get that the optimal predictor of $Y_{t + h}$ based on $Y_s, s \le t$, thresholds $\pred{opt}$ appropriately.
\end{proof}

\begin{remark}
    It is assumed that the series in \eqref{eq:Y_t_def} converges almost surely for every $t \in \Z$. There are different sets of conditions on $\{a_j\}$ and $\{\epsilon_t\}$, appropriate for different scenarios, that guarantee this; see, e.g., Propositions 3.1.1 and 13.3.1 in \cite{brockwell1991time}. In Section~\ref{subsect:MA_infty_mod_opt_pred_asymp_aspects}, to calculate the asymptotic precision of the predictor in Theorem~\ref{thm:MA_infty_opt_pred}, we use a set of conditions from \cite{kulik2020heav}.
\end{remark}

\begin{remark}\label{rem:finite_history_pred}
    The optimal predictor $\pred{opt}$ checks whether the sum in \eqref{eq:MA_infty_opt_pred} is above a suitably chosen threshold. This sum can be written as
    \begin{equation}\label{eq:eq-opt-pred}
    \sum_{j = 0}^{\infty} \sum_{k = 0}^{\infty} a_{j + h}b_k Y_{t - j - k} = \sum_{r = 0}^{\infty}\left[\sum_{s = 0}^r a_{s + h}b_{r - s}\right]Y_{t - r}
    =: \sum_{r = 0}^{\infty} c_r Y_{t - r},
    \end{equation}
    which makes it clear that $\pred{opt}$ is based on infinitely many past observations. However, in practice, only finitely many past observations would be available at a given time $t$. Suppose that for some $\ell \ge 1$, the $\ell$ most recent observations, i.e., $Y_t, \ldots, Y_{t - \ell + 1}$, are available. Then a natural approximation to the optimal predictor is obtained by truncating \eqref{eq:eq-opt-pred} as follows:
    \[
    \widehat{Y}_{t + h}^{(\ell)} = \sum_{r = 0}^{\ell - 1} c_r Y_{t - r},
    \]
    where the $c_r$'s can be obtained using the formula $c_r = \sum_{s = 0}^r a_{s + h}b_{r - s}$. Under 
    various models, one can study the accuracy of the approximation 
    $\widehat{Y}_{t + h}^{(\ell)}$ to $\pred{opt}$ and the respective precisions in terms of
    the rate of decay of the $c_r$'s. We do not pursue this here.
\end{remark}

In the rest of this section, we consider the general autoregressive model:
\begin{equation}\label{eq:AR_mod}
    Y_t  = \sum_{i=1}^d \phi_i Y_{t-i} + \epsilon_t,\ \ \ t\in \Z.
\end{equation}
We shall assume that the innovations satisfy Assumption~\ref{assump:innov_density} and that for each $t$, $\epsilon_t$ is independent from the history of the process up to time $t-1$, i.e., the $\sigma$-field ${\cal H}_{t-1} := \sigma\{Y_s,\ s\le t-1\}$. Under the usual conditions that the complex polynomial $\pi(z):= 1 -\sum_{i=1}^d \phi_i z^i$ has its roots outside the unit circle, the autoregressive equation \eqref{eq:AR_mod} has a unique stationary solution:
\begin{equation}\label{e:Y-as-a-moving-average}
 Y_t = \sum_{i=0}^\infty a_i \epsilon_{t-i},
\end{equation}
where $\sum_{i=0}^\infty a_i z^i = 1/\pi(z),\ |z|\le 1$ is the series representation of the function $1/\pi(z)$, which is analytic in the unit circle $\{|z| \le 1\}$. The coefficients $a_i$ in the moving average representation \eqref{e:Y-as-a-moving-average} decay exponentially fast and therefore, since $\E[|\epsilon_t|^\delta]<\infty$ (Assumption \ref{assump:innov_density}), 
the series in \eqref{e:Y-as-a-moving-average} is convergent almost surely 
\citep[see, e.g., Propositions 13.3.1 and 13.3.2 in][]{brockwell1991time}. 

The following lemma will help us derive the optimal predictor for the model in \eqref{eq:AR_mod}. For $t \ge u$, we use the notation $Y_{t:u}$ to succinctly represent the vector $(Y_t, \ldots, Y_u)^{\top}$.
\begin{lemma}\label{lem:AR_opt_pred_recursion1}
    Let $\Phi$ be the $d \times d$ matrix $(\begin{matrix} \phi & e_1 & \cdots & e_{d - 1} \end{matrix})$, with $e_1, \ldots, e_{d - 1}$ being the first $d - 1$ standard basis vectors of $\R^d$. For each $h \ge 1$, set $\phi(h) = \Phi^h e_1$. Then
    \begin{equation}\label{eq:y_t_plus_h_expr1}
        Y_{t + h} = \phi(h)^{\top}Y_{t:(t - d + 1)} + \eta,
    \end{equation}
    where $\eta$ is a random variable that is independent of $Y_{t:(t - d + 1)}$ and has a density with respect to the Lebesgue measure.
\end{lemma}
\begin{proof}
    See Section~\ref{subsect:AR_opt_pred_recursion_lem_pf}.
\end{proof}

In the following result that draws on the version of Theorem \ref{thm:add_err_mod_thm} in Remark~\ref{rem:homoscedastic_noise}, we characterize the optimal predictor for the model in \eqref{eq:AR_mod}.
\begin{theorem}\label{thm:AR_opt_pred}
    Suppose that $\{Y_t\}$ is the unique stationary solution given by 
    \eqref{e:Y-as-a-moving-average} for the AR$(d)$ equation in \eqref{eq:AR_mod}. Set
    \begin{equation}\label{eq:AR_opt_pred}
        \AROptPred{t}{h}{\phi} := \phi(h)^{\top}Y_{t:(t - d + 1)},
    \end{equation}
    with $\phi(h)$ being defined as in Lemma~\ref{lem:AR_opt_pred_recursion1}. Let $p \in (0, 1)$ and $y_0 \in \R$ such that $\P(Y \le y_0) = p$. Then the optimal predictor of $\I(Y_{t+h} > y_0)$ via $Y_{t:(t - d + 1)}$ is $\I(\AROptPred{t}{h}{\phi} > \tau)$, where $\tau$ is a constant that calibrates the predictor.
\end{theorem}
\begin{proof}
    Using Lemma~\ref{lem:AR_opt_pred_recursion1}, we can write $Y_{t + h} = \AROptPred{t}{h}{\phi} + \eta$, where $\AROptPred{t}{h}{\phi} \ind \eta$. Both $\AROptPred{t}{h}{\phi}$ and $\eta$ have densities with respect to the Lebesgue measure. It follows from the result in Remark~\ref{rem:homoscedastic_noise} that the optimal predictor is $\I(\AROptPred{t}{h}{\phi} > \tau)$, where $\tau$ is chosen to achieve the desired calibration level.
\end{proof}

\begin{remark}
    As established by \cite{slud1982acha}, when $h = 1$, the decomposition \eqref{eq:y_t_plus_h_expr1} where $\eta$ and $Y_{t:(t - d + 1)}$ are independent is only possible if $\{Y_t\}$ is either an 
    AR$(d)$ or a 
    Gaussian time series. This suggests that there is no simple closed-form expression for 
    the optimal finite-history predictor in general, non-Gaussian MA$(\infty)$ models.
\end{remark}

\begin{remark}\label{rem:opt_pred_finite_var_ar}
    When the process $\{Y_t\}$ has finite variances, the random variable $\AROptPred{t}{h}{\phi}$ given in \eqref{eq:AR_opt_pred} is the same as the linear combination of $Y_t, \ldots, Y_{t - d + 1}$ that approximates $Y_{t + h}$ with minimal mean squared error (MSE), i.e.,
    \begin{equation}\label{eq:AR_MSE}
        \phi(h) = \argmin_{\varphi \in \R^d} \E[(Y_{t + h} - \varphi^{\top}Y_{t:(t - d + 1)})^2].
    \end{equation}
    See Section~\ref{subsubsect:opt_pred_finite_var_ar} for a proof.
\end{remark}

\subsection{Statistical Inference for Optimal Predictors in AR($d$) Models}\label{subsect:AR_practical_pred}

The optimal predictor for AR$(d)$ models requires knowledge of $\phi$, as well as the distribution of $\pred{opt}$ in order to calibrate the predictor with the appropriate quantile. In practice, both $\phi$ and the quantile need to be estimated. In this section, we study a natural plug-in approximation to the optimal predictor $\pred{opt}$ and its calibration from training data. 
We establish the asymptotic consistency and optimality of the inferential procedure as the size of the training sample goes to infinity.  We do so under very general conditions using
a uniform law of large numbers for empirical processes established in \cite{adams2012unif} (see the proof of Proposition~\ref{prop:U_convergence}).

\subsubsection{Estimation of the Optimal Predictor}\label{subsubsect:AR_approx_opt_pred}

As mentioned above, in practice we use an empirical plug-in counterpart of the optimal
predictor in \eqref{eq:AR_opt_pred}.  To this end, we need to briefly 
discuss the estimation of AR$(d)$ time series with heavy-tailed innovations.  

There is a sizable literature on the consistency and asymptotic distribution of 
estimators of the coefficient vector of an AR($d$) model. For example, the seminal work 
of \cite{davis1992mest} studies a general family of M-estimators and in particular the 
least absolute deviation (LAD).  Suppose that we have $n$ consecutive observations 
$Y_s,\ s=t-n+1,\cdots,t$ from the AR$(d)$ time series. Then the LAD estimator of 
the vector $\phi$ of AR coefficients is defined as
\begin{equation}\label{e:LAD-estimator}
\wh \phi := \mathop{\rm argmin}_{\varphi\in\R^d}
\frac{1}{n-d} \sum_{s=t-n+d}^{t-1} | Y_{s+1} - \varphi^\top Y_{s:(s-d+1)}|.
\end{equation}
Under general regular variation assumptions on the innovations (see Section~\ref{subsect:MA_infty_mod_opt_pred_asymp_aspects}), \cite{davis1992mest} established the consistency of the LAD estimator, i.e.,
$\wh \phi \overset{\P}{\to} \phi$ as $n\to\infty$.  In fact, the rate of convergence 
was established as well as the asymptotic distribution using point process convergence
techniques.

Inference in heavy-tailed time series models is well-developed. 
Notably, the monograph of \cite{koul-book} develops a powerful approach using 
weighted empirical processes which applies to stationary as well as non-stationary
models. In particular, \cite{koul2010acla} proves asymptotic normality of a robust 
minimum distance estimator of the coefficient vector under assumptions that allow the 
innovations to have infinite variance or even infinite mean. Our aim is not to contribute 
to the established literature on coefficient estimation for AR$(d)$ models, but to study
the estimation of optimal extreme event predictors.  To this end, we shall assume that we
have a {\em consistent} estimator of $\wh \phi$ obtained either using the LAD or 
another robust M-estimator.

Recall that our goal is to predict exceedance of the $p$th quantile $\Fi{Y}(p)$ of the marginal distribution of $Y_t$, when $\{Y_t\}$ is a stationary AR($d$) process satisfying \eqref{eq:AR_mod}. Our main result will require the following assumption.
\begin{assumption}\label{assump:F_Y_cont_at_quantile}
    The distribution function $F_Y$ is continuous at $\Fi{Y}(p)$.
\end{assumption}
Under this assumption, the events $\{Y > \Fi{Y}(p)\}$ and $\{Y \ge \Fi{Y}(p)\}$ differ by an event of probability zero, so instead of trying to predict $\I(Y_{t + h} > \Fi{Y}(p))$, we can try to predict $\I(Y_{t + h} \ge \Fi{Y}(p))$. Theorem~\ref{thm:AR_opt_pred} suggests using $\I(\approxAROptPred{t}{h}{\phi} \ge \Fi{\approxAROptPred{t}{h}{\phi}}(p))$ to do this, where $\widehat{\phi}$ is an estimator of $\phi$. However, in practice, we would not be able to predict with this, as the quantile $\Fi{\approxAROptPred{t}{h}{\phi}}(p)$ would be unknown. We could instead predict with
\begin{equation}\label{eq:AR_practical_pred}
    \I(\approxAROptPred{t}{h}{\phi} \ge \widehat{F}_{\approxAROptPred{t}{h}{\phi}}^{\leftarrow}(p)),
\end{equation}
where $\widehat{F}_{\approxAROptPred{t}{h}{\phi}}$ is the empirical distribution function of $\approxAROptPred{t}{h}{\phi}$. Suppose that we have observed $Y_t, \ldots, Y_{t - n + 1}$, where $n \ge d$, the order of the time series. For $k \ge 0$, define $Z_{t - k} := Y_{(t - k):(t - k - d + 1)}$. We would be able to compute $\approxAROptPred{t}{h}{\phi}, \ldots, \approxAROptPred{t - n + d}{h}{\phi}$, and for any $u \in \R$, we could calculate $\widehat{F}_{\approxAROptPred{t}{h}{\phi}}(u)$ as
\begin{equation}\label{eq:ecdf}
    \widehat{F}_{\approxAROptPred{t}{h}{\phi}}(u)
    = \frac{1}{n - d + 1}\sum_{k = 0}^{n - d} \I(\approxAROptPred{t - k}{h}{\phi} \le u)
    = \frac{1}{n - d + 1}\sum_{k = 0}^{n - d} \I(\widehat{\phi}(h)^{\top}Z_{t - k} \le u).
\end{equation}
These considerations motivate the following methodology:
\begin{algorithm}
\caption{Extreme Event Prediction for AR($d$) Time Series}
\label{alg:AR_pred}
\begin{algorithmic}[1]
    \Statex \textbf{Inputs:} observations $Y_t, \ldots, Y_{t - n + 1}$, lead time $h$, quantile level $p \in (0,1)$
    \State Compute an estimate $\widehat{\phi}$ of $\phi$ using, e.g., LAD estimation in \eqref{e:LAD-estimator}.
    \State Set $\widehat{\Phi} = \begin{pmatrix} \widehat{\phi} & e_1 & \cdots & e_{d - 1} \end{pmatrix}$, where $\{e_1, \ldots, e_d\}$ is the standard basis of $\mathbb{R}^d$.
    \State Set $\widehat{\phi}(h) = \widehat{\Phi}^h e_1$.
    \State Compute $\approxAROptPred{t}{h}{\phi}, \ldots, \approxAROptPred{t - n + d}{h}{\phi}$ and their $p$th sample quantile $\widehat{F}_{\approxAROptPred{t}{h}{\phi}}^{\leftarrow}(p)$.
    \State Predict that $Y_{t + h} \ge \Fi{Y}(p)$ iff $\approxAROptPred{t}{h}{\phi} \ge \widehat{F}_{\approxAROptPred{t}{h}{\phi}}^{\leftarrow}(p)$.
\end{algorithmic}
\end{algorithm}

\begin{remark}
    Given an estimator $\widehat{\phi}$ of $\phi$, the matrix $\Phi$ appearing in Lemma~\ref{lem:AR_opt_pred_recursion1} can be approximated by the matrix $\widehat{\Phi} = (\begin{matrix} \widehat{\phi} & e_1 & \cdots & e_{d - 1} \end{matrix})$. This yields a plug-in approximation $\widehat{\phi}(h) = \widehat{\Phi}^h e_1$ to $\phi(h)$ and a plug-in approximation $\approxAROptPred{t}{h}{\phi}$ to $\AROptPred{t}{h}{\phi}$. The coefficient vector $\phi$ can be estimated using methods like LAD or least squares estimation (if applicable) from a training 
    set consisting of $n$ pairs $(Y_t, Y_{(t - 1):(t - d)}), \ldots, (Y_{t - n + 1}, Y_{(t - n):(t - n - d + 1)})$, where the response in the $k$th pair is $Y_{t - k + 1}$ and the covariate vector is $Y_{(t - k):(t - k - d + 1)}$. Note that we are assuming that the order $d$ is known. The approximations $\widehat{\phi}(h)$ and $\approxAROptPred{t}{h}{\phi}$ consistently approximate $\phi(h)$ and $\AROptPred{t}{h}{\phi}$, respectively; see Lemma~\ref{lem:AR_consistency1}.
\end{remark}

To make it easier to derive the properties of the predictor \eqref{eq:AR_practical_pred}, we express it in a simpler form. By Equation (2.6) in \cite{resnick2007heav},
\begin{equation}\label{eq:rank_transformed_pred_non_U}
    \I(\approxAROptPred{t}{h}{\phi} \ge \widehat{F}_{\approxAROptPred{t}{h}{\phi}}^{\leftarrow}(p))
    = \I(\widehat{F}_{\approxAROptPred{t}{h}{\phi}}(\approxAROptPred{t}{h}{\phi}) \ge p).
\end{equation}
Define $H_{\varphi,\,\tau} := \{y \in \R^d : \varphi^{\top}y \le \tau\}$, which is a half-space in $\R^d$. Write
\[
\wh U_{t+h} (\varphi,\tau) := \frac{1}{n-d+1}\sum_{k=0}^{n-d} \I(Z_{t - k} \in H_{\varphi,\tau})
\ \ \mbox{ and } \ \
\wh U_{t+h} := \wh U_{t+h} (\widehat{\phi}(h),\,\approxAROptPred{t}{h}{\phi})
\]
Since, by \eqref{eq:ecdf},
\[
\widehat{F}_{\approxAROptPred{t}{h}{\phi}}(\approxAROptPred{t}{h}{\phi})
= \frac{1}{n - d + 1}\sum_{k = 0}^{n - d} \I(\widehat{\phi}(h)^{\top}Z_{t - k} \le \approxAROptPred{t}{h}{\phi}),
\]
it follows from \eqref{eq:rank_transformed_pred_non_U} that
\begin{equation}\label{eq:rank_transformed_pred_U}
    \I(\approxAROptPred{t}{h}{\phi} \ge \widehat{F}_{\approxAROptPred{t}{h}{\phi}}^{\leftarrow}(p))
    = \I(\widehat{U}_{t + h} \ge p).
\end{equation}

We introduce some additional notation that will be useful in the next section. Let $Z$ be a random vector that has the same distribution as $Y_{t:(t - d + 1)}$ and is independent of $(Y_{t:(t-d+1)}, \widehat{\phi})$. Define
\begin{equation}\label{e:U-varphi,tau}
    U_{t+h}(\varphi,\tau) := \P_Z(H_{\varphi,\tau} ) = \E[ \I(Z\in H_{\varphi,\tau})]
    \ \mbox{ and } \
    U_{t+h}^* :=  U_{t+h} (\phi(h),\,\AROptPred{t}{h}{\phi}).
\end{equation}

\subsubsection{Asymptotic Optimality of the Estimated Predictor}\label{subsect:fixed_p_consistency}

The following result shows that the predictor in \eqref{eq:rank_transformed_pred_U} is asymptotically calibrated and optimal, under mild general conditions that essentially amount to
requiring that the coefficient vector $\phi$ be estimated consistently.

\begin{theorem}\label{thm:fixed-p}
    Suppose that $\widehat{\phi} \stackrel{\P}{\to} \phi$ and that $\AROptPred{t}{h}{\phi}$ has a density with respect to Lebesgue measure. Also suppose that the condition in Assumption~\ref{assump:F_Y_cont_at_quantile} holds. Then the predictor in \eqref{eq:rank_transformed_pred_U} is:
    
    (i) Asymptotically calibrated, i.e.,
    \begin{equation}\label{e:thm:fixed-p-i}
        \P[\widehat{U}_{t + h} \ge p] \to 1 - p = \P[Y_{t + h} \ge \Fi{Y}(p)] ,\ \ \mbox{ as }n\to\infty.
    \end{equation}
    
    (ii) Asymptotically optimal, i.e.,
    \begin{equation}\label{e:thm:fixed-p-ii}
        \Prec_{Y_{t+h} \ge \Fi{Y}(p)}(\widehat{U}_{t + h} \ge p) \to \Prec_{Y_{t+h} \ge \Fi{Y}(p)}(U_{t+h}^* \ge p),\ \ \mbox{ as }n\to\infty.
    \end{equation}
\end{theorem}
To prove this theorem, we will need the following results.
\begin{proposition}\label{prop:U_convergence}
    Suppose that $\widehat{\phi} \xrightarrow{\P} \phi$ and that $\AROptPred{t}{h}{\phi}$ has a density with respect to Lebesgue measure. 
    Then $\widehat{U}_{t + h} \xrightarrow{\P} U_{t + h}^*$, as $n \to \infty$,
    where $U_{t+h}^*$ is as in \eqref{e:U-varphi,tau}.
\end{proposition}
\begin{proof}
    See Section~\ref{subsubsect:U_convergence_prop_pf}.
\end{proof}

\begin{lemma}\label{lem:stable-conv}
    Let $\xi_n\stackrel{\P}{\to}\xi$ and also $x_n\to x$, as $n\to\infty$. If  $\P[\xi = x] = 0$,
    then, for all events $A$, as $n\to\infty$, we have
    \begin{equation}\label{e:l:stable-conv}
        \P[ A\cap\{ \xi_n \ge x_n \} ] \longrightarrow \P[ A\cap\{ \xi \ge x\}].
    \end{equation}
\end{lemma}
\begin{proof}
    Recall that $\xi_n\stackrel{\P}\to \xi$, if and only if, for every $n_k\to\infty$, there is a further subsequence 
    $\{m_k\}\subset\{n_k\}$, such that $\xi_{m_k}\to \xi$, almost surely, as $m_k\to\infty$ \citep[see Theorem 6.3.1(b) in][]{resnick2014apro}.
    
    Applying this characterization and using the fact that, by assumption, $\P[\xi = x] = 0$, it follows that
    \[
    I_A I_{\{\xi_{m_k} \ge x_{m_k}\}} \stackrel{a.s.}{\longrightarrow} I_A I_{\{\xi \ge x\}},\ \ \mbox{ as }m_k\to \infty.
    \]
    This, in view of the dominated convergence theorem, implies that
    \begin{equation}\label{e:l:stable-conv-1}
        \P[A \cap\{\xi_{m_k} \ge x_{m_k}\}] \to \P[A \cap\{\xi \ge x\}].
    \end{equation}
    
    We have thus shown that every infinite sequence $n_k\to\infty$, has a further sub-sequence  
    $m_k\to\infty$, along which \eqref{e:l:stable-conv-1} holds.  This implies \eqref{e:l:stable-conv} and completes the
    proof of the lemma.
\end{proof}
We are now in a position to prove Theorem~\ref{thm:fixed-p}.
\begin{proof}[Proof of Theorem~\ref{thm:fixed-p}]
    Introduce the random variables $\xi_n:= \widehat{U}_{t + h}$ and $\xi:= U_{t + h}^*$, and let $x_n:= x = p$. Because $\widehat{\phi} \stackrel{\P}{\to} \phi$, we get from Proposition~\ref{prop:U_convergence} that $\xi_n \xrightarrow{\P} \xi$. By Lemma \ref{lem:stable-conv} applied with $A$ being the entire sample space, we obtain \eqref{e:thm:fixed-p-i}. Note that $\P(\xi \ge p) = 1 - p = \P(Y_{t + h} \ge \Fi{Y}(p))$ since $U_{t + h}^* = F_{\AROptPred{t}{h}{\phi}}(\AROptPred{t}{h}{\phi})$ is standard uniform and $F_Y$ is continuous at $\Fi{Y}(p)$. This establishes part {\em (i)}, i.e., the asymptotic calibration of the predictor.\\
    
    Now, to prove part {\em (ii)}, observe that
    \begin{align*}
    \P(Y_{t+h} \ge \Fi{Y}(p) \mid \widehat{U}_{t + h} \ge p) &= \frac{\P[A \cap\{\xi_n \ge x_n\} ]}{ \P[\widehat{U}_{t + h} \ge p]},\ \ \ \mbox{ and }\\
     \P(Y_{t+h} \ge \Fi{Y}(p) \mid U_{t + h}^* \ge p) &= \frac{\P[A \cap\{\xi \ge x\} ]}{ \P[U_{t + h}^* \ge p]},
    \end{align*}
    where $A:=\{ Y_{t+h} \ge \Fi{Y}(p)\}$ is a fixed event that does not depend on $n$.  Thus, in view of the already established convergence \eqref{e:thm:fixed-p-i}, to prove \eqref{e:thm:fixed-p-ii}, it is enough to show that
    $\P[A \cap\{\xi_n \ge x_n\} ]\to \P[A \cap\{\xi \ge x\} ]$, which however follows from Lemma \ref{lem:stable-conv}.
\end{proof}

\begin{remark}
    See Section~\ref{subsect:ar_mod_sim_study} for the results of a simulation study of the performance of the optimal AR predictor $\I(\AROptPred{t}{h}{\phi} > \Fi{\AROptPred{t}{h}{\phi}}(p))$ and two approximations to it of the form $\I(\approxAROptPred{t}{h}{\phi} \ge \widehat{F}_{\approxAROptPred{t}{h}{\phi}}^{\leftarrow}(p))$, where the quantile estimator $\widehat{F}_{\approxAROptPred{t}{h}{\phi}}^{\leftarrow}(p)$ is either a sample quantile as in this section or an estimator based on extreme value theory.
\end{remark}

\subsection{The Optimal Extremal Precision for MA($\infty$) Models}\label{subsect:MA_infty_mod_opt_pred_asymp_aspects}

In this section we derive the optimal extremal precision (see Section~\ref{subsect:opt_pred_asymp_aspects}) for the MA($\infty$) model
\begin{equation}\label{eq:ma_infty_mod}
    Y_t = \sum_{j = 0}^{\infty} a_j\epsilon_{t - j}.
\end{equation}
Theorem~\ref{thm:MA_infty_opt_asymp_precision} will give the optimal extremal precision as the precision as $p \uparrow 1$ of the optimal predictor. That result will rely on certain assumptions about the sequence of innovations $\{\epsilon_t\}$ and the sequence of coefficients $\{a_j\}$, assumptions that need to be made in order for $Y_t$ to be a random variable that is finite almost surely. Our assumption about the innovations uses the concept of regular variation, whose definition for random variables we recall here. We use the terminology of \cite{kulik2020heav}. In particular, we denote the class of functions that are regularly varying with tail index $\rho$ by $\RV_{\rho}$. See Section~\ref{subsect:rv_info} for additional regular variation definitions and notation.

\begin{definition}\label{def:rv_rand_var2}
    A random variable $X$ is regularly varying with tail index $\alpha > 0$ and extremal skewness coefficient
    $p_X \in [0, 1]$ if
    \[
    \overline{F}_{|X|}(x) := \P(|X| > x) \in \RV_{-\alpha}\ \ 
    \mbox{ and }\ \ \lim_{x \to \infty} \frac{\overline{F}_X(x)}{\overline{F}_{|X|}(x)} = p_X.
    \]
\end{definition}

We now state our assumption about $\{\epsilon_t\}$.

\begin{assumption}\label{assump:ma_infty_innovs}
    The innovations $\epsilon_t,\ t\in\Z$, are iid and regularly varying with tail exponent $\alpha>0$ and extremal skewness coefficient $p_{\epsilon}$. Furthermore, if $\alpha > 1$, we shall assume that $\E(\epsilon_t) = 0$.
\end{assumption}

We will consider various infinite linear combinations of the innovations $\epsilon_j$'s. To this end, it
will be convenient to use the following notation:
\begin{equation}\label{e:xi-a-def}
  \series(b) = \sum_{j = 0}^{\infty} b_j \epsilon_j,
\end{equation}
where $b = \{b_j\}$ is a sequence of real coefficients. Corollary 4.2.1 of \cite{kulik2020heav} implies 
that under Assumption \ref{assump:ma_infty_innovs}, the series in  \eqref{e:xi-a-def} converges almost surely, provided the sequence $b$ belongs to the class of infinite sequences:
\begin{equation}\label{e:C_alpha}
    \seqSet:=\Big\{ \{b_j\}\, :\, \sum_{j = 0}^{\infty} |b_j|^{\delta} < \infty,\ 
    \mbox{ for some }\delta\in(0,\alpha)\cap (0,2] \Big\}.
\end{equation}
Introduce the following notation:
\begin{align}
\multiplier{\pm}{b_j} &:= p_\epsilon \I(\pm b_j>0) + (1-p_\epsilon) \I(\pm b_j<0) \label{def:multiplier} \\
\mmultiplier{\pm}{\pm}{b_j}{b'_j} &:= p_\epsilon \I(\pm b_j>0, \pm b'_j>0) + (1-p_\epsilon) \I(\pm b_j<0, \pm b'_j<0) \label{def:mmultiplier} \\
\normConst{\pm}{b}{h} &:= \sum_{j = h}^{\infty} \multiplier{\pm}{b_j}|b_j|^\alpha. \label{def:norm_const}
\end{align}

As a consequence of Corollary 4.2.1 in \cite{kulik2020heav}, we have the following result:
\begin{lemma}\label{lem:xi_a_facts}  For all $\alpha>0$ and sequences $b \in \seqSet$, the series \eqref{e:xi-a-def} converges almost surely.
Moreover, $\series(b)$ is regularly varying with tail index $\alpha$ and extremal skewness coefficient
    \[
    p_{\series(b)} = \frac{\normConst{+}{b}{0}}{\sum_{j = 0}^{\infty} |b_j|^{\alpha}}.
    \]
\end{lemma}

\begin{proof}
    See Section~\ref{subsect:xi_a_facts_lem_pf}.
\end{proof}


Recall the tail dependence coefficient $\lambda$ defined in Definition~\ref{e:lambda-xi-eta}.
\begin{lemma}\label{lem:two_series_tdc}
    Let $a, b \in \seqSet$. Assume that $p_{\series(a)} > 0$ and $p_{\series(b)} > 0$. Then
    \[
    \lambda\left(\series(a), \series(b)\right)
    = \sum_{j = 0}^{\infty} \mmultiplier{+}{+}{a_j}{b_j}\left(\frac{|a_j|^{\alpha}}{\normConst{+}{a}{0}} \bigwedge \frac{|b_j|^{\alpha}}{\normConst{+}{b}{0}}\right).
    \]
\end{lemma}
\begin{proof}
    See Section~\ref{subsect:lem:two_series_tdc_lem_pf}.
\end{proof}


We now derive the optimal extremal precision for the model in \eqref{eq:ma_infty_mod}. As in Theorem~\ref{thm:MA_infty_opt_pred}, we will need to assume that $\{Y_t\}$ is invertible. For ease of reference, we state this assumption here.
\begin{assumption}\label{assump:ma_infty_invert}
    The process defined by \eqref{eq:ma_infty_mod} is invertible, i.e., there exists a real sequence $\{b_j\}_{j = 0}^{\infty}$ such that for all $t$, $\epsilon_t = \sum_{j = 0}^{\infty} b_j Y_{t - j}$ almost surely.
\end{assumption}
\begin{theorem}\label{thm:MA_infty_opt_asymp_precision}
    Let $\{\epsilon_t\}$ be a sequence of random variables satisfying Assumption~\ref{assump:ma_infty_innovs}. Consider a sequence of coefficients $a \in \seqSet$, and define
    \begin{equation}\label{e:Yt-MA}
    Y_t = \sum_{j = 0}^{\infty} a_j\epsilon_{t - j}, \ t \in \Z.
    \end{equation}
    Suppose moreover that Assumption~\ref{assump:ma_infty_invert} holds. Then for every $h \ge 1$, the extremal optimal precision for prediction of $\I(Y_{t + h} > y_0)$ is
    \begin{equation}\label{eq:MA_infty_asymp_opt_pred}
    \lambda_h^{\rm (opt)} := \lambda(Y_{t + h}, \pred{opt}) = \frac{\normConst{+}{a}{h}}{\normConst{+}{a}{0}},    
    \end{equation}
    where $\normConst{+}{a}{h}$ and $\normConst{+}{a}{0}$ are defined according to \eqref{def:norm_const}, assuming that $\normConst{+}{a}{h} > 0$.
\end{theorem}
\begin{proof}
    Using the expression for the optimal predictor from Theorem~\ref{thm:MA_infty_opt_pred}, we have
    \[
    \lambda_h^{\rm (opt)}
    = \lambda\left(Y_{t + h}, \sum_{j = 0}^{\infty} a_{j + h}\epsilon_{t - j}\right)
    = \lambda\left(\sum_{j = 0}^{\infty} a_j\epsilon_{t + h - j}, \sum_{j = 0}^{\infty} a_{j + h}\epsilon_{t - j}\right)
    = \lambda\left(\sum_{j = -h}^{\infty} a_{j + h}\epsilon_{t - j}, \sum_{j = 0}^{\infty} a_{j + h}\epsilon_{t - j}\right).
    \]
    Applying Lemma~\ref{lem:two_series_tdc}, we get
    \begin{equation}\label{eq:lambda_opt_expr}
        \lambda_h^{\rm (opt)}
        = \sum_{j = 0}^{\infty} \mmultiplier{+}{+}{a_{j + h}}{a_{j + h}}\left(\frac{|a_{j + h}|^{\alpha}}{\normConst{+}{a}{0}} \bigwedge \frac{|a_{j + h}|^{\alpha}}{\normConst{+}{a}{h}}\right)
        = \sum_{j = 0}^{\infty} \kappa_+(a_{j + h})\frac{|a_{j + h}|^{\alpha}}{\normConst{+}{a}{0}}
        = \frac{\normConst{+}{a}{h}}{\normConst{+}{a}{0}}.
    \end{equation}
\end{proof}

At time $t$, the data available for prediction of $Y_{t + h}$ consists of at least $Y_t$. One may wonder how a predictor based solely on $Y_t$ compares to the optimal predictor. Proposition~\ref{prop:oracle} describes circumstances under which these two predictors have the same extremal precision, or the extremal optimal precision can at least be expressed or bounded in terms of the extremal precision of $\pm Y_t$.
\begin{proposition}\label{prop:oracle}
Let $\{\epsilon_t\}$ be a sequence of random variables satisfying Assumption~\ref{assump:ma_infty_innovs} and let $a \in \seqSet$ such that Assumption~\ref{assump:ma_infty_invert} holds. Suppose that for some $h\in{\mathbb N}$, $\normConst{+}{a}{h}>0$ and the sequence $\{a_j\}$ is lag-$h$ absolute decreasing, i.e., 
\begin{equation}\label{e:aj-absolute-monotone}
 |a_{j+h}| \le |a_j|,\ \ \mbox{ for all }j\ge 0.
\end{equation}

Then: 

{\em (i)} If either $\normConst{-}{a}{0} = 0$ or $a_{j+h} a_j\ge 0$ for all $j\ge 0$, we have
\begin{equation}\label{e:p:oracle-i}
 \lambda(Y_{t+h},Y_t) = \lambda(Y_{t + h}, \pred{opt}).
\end{equation}

{\em (ii)} If $\normConst{-}{a}{0}>0$, we have:
\begin{equation}\label{e:p:oracle}
 \lambda(Y_{t+h},Y_t) + \lambda(Y_{t+h},-Y_t) \le\  \lambda(Y_{t + h}, \pred{opt}) \ 
 \le \lambda(Y_{t+h},Y_t) + \Big( 1\bigvee \frac{\normConst{-}{a}{0}}{\normConst{+}{a}{0}}\Big)\lambda(Y_{t+h},-Y_t).
\end{equation}

{\em (iii)} If $p_\epsilon = 1/2$, we have $\normConst{-}{a}{0} = \normConst{+}{a}{0}$ and $\lambda(Y_{t+h},Y_t) + \lambda(Y_{t+h},-Y_t) = \lambda(Y_{t + h}, \pred{opt})$.
\end{proposition}
\begin{proof}
    See Section~\ref{subsubsect:oracle_prop_pf}.
\end{proof}

\begin{remark}
    We conjecture that when $Y_t$ has a symmetric distribution,
    \[
    \lambda(Y_{t+h},Y_t) + \lambda(Y_{t+h},-Y_t) = \lambda(Y_{t+h}, |Y_t|).
    \]
\end{remark}

\begin{example}
Using the MA$(\infty)$ representation of AR$(d)$ models and Relation 
\eqref{eq:MA_infty_asymp_opt_pred} 
one can numerically compute the optimal extremal precision 
$\lambda_h^{\rm (opt)}$ in terms of the vector of autoregressive coefficients 
$\phi$ and the extremal skewness coefficient $p_\epsilon$ of the innovations.  In the simple AR$(1)$ setting, however,
a closed-form expression is attainable.  Indeed, in this case,
\eqref{e:Yt-MA} holds with $a_j=\phi^j$ when $|\phi|<1$, and 
thus, if $p_\epsilon>0$, from \eqref{eq:MA_infty_asymp_opt_pred}
and \eqref{def:multiplier}--\eqref{def:norm_const}, we obtain:
\[
\lambda_h^{\rm (opt)} = \lambda(Y_{t+h},\pred{opt}) 
= \frac{\eta_+(a,h)}{\eta_+(a,0)}= 
\left\{ \begin{array}{ll} 
 |\phi|^{h\alpha} &,\ \mbox{ if }\phi\ge 0\mbox{ or $h$ is even }\\
 |\phi|^{h\alpha} \frac{(p_\epsilon|\phi|^\alpha+1-p_\epsilon)}{p_\epsilon + (1-p_\epsilon)|\phi|^\alpha} &, \mbox{ if $\phi<0$ and $h$ is odd. }
\end{array}
\right. 
\]
In particular, for symmetric innovations ($p_{\epsilon} = 1 / 2$), the above asymptotic optimal precision simplifies to $\lambda_h^{\rm (opt)} = |\phi|^{h\alpha}$.
\end{example}

\section{Data Analysis}\label{sect:data_analysis}

We illustrate the use of our methodological and theoretical results with an analysis of soft X-ray flux data. As stated in Section~\ref{sect:intro}, the peak value of the flux during a flare determines the flare's strength. We fit several models to a time series of flux values and used them to predict exceedances of high thresholds at multiple lead times. Section~\ref{subsect:data_processing} describes how the raw data was processed into the data used to fit models. FARIMA models were one of the types of models we fit; Section~\ref{subsect:farima_background} gives background on them, while Section~\ref{subsect:flux_farima_mods} describes how they were fit and discusses the fitted models. The results for the various prediction scenarios are presented in Section~\ref{subsect:pred_results}.

\subsection{Data Processing}\label{subsect:data_processing}

There are five classes of solar flares, A, B, C, M, and X, with the A class containing the weakest flares and the X class containing the strongest \citep[][]{fletcher2011anob}. Solar flare forecasting methods are typically used to forecast M- and X-class flares \citep[][]{chen2023edit}, i.e., those with a peak flux above $\SI{e-5}{W.m^{-2}}$ and $\SI{e-4}{W.m^{-2}}$, respectively \citep[][]{mothersbaugh2023read}. Flux values are recorded by the Geostationary Operational Environmental Satellites (GOES) of the U.S. National Oceanic and Atmospheric Administration (NOAA) and are published in several forms; we downloaded data files containing averages over one-minute intervals. Instead of using the raw data from the satellites, we used science-quality data, which contains a number of corrections \citep[][]{mothersbaugh2023read}.

\begin{figure}[!htb]
    \centering
    \includegraphics[width=0.7\linewidth]{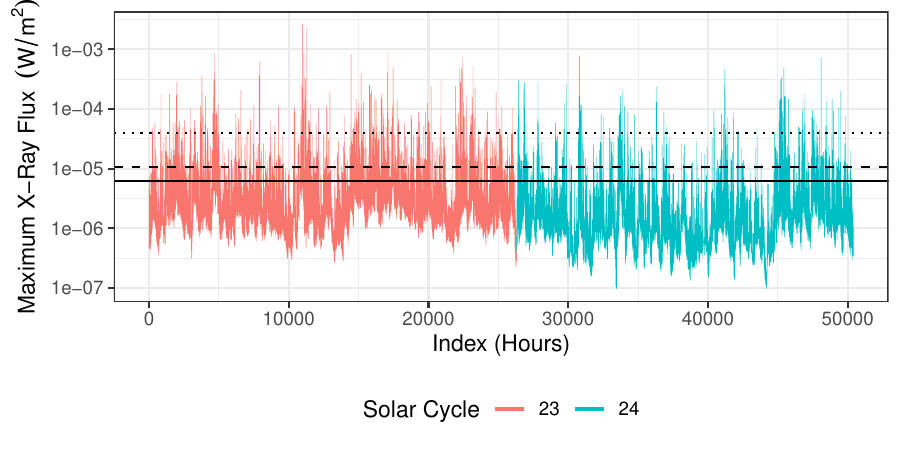}
    \caption{{\small A plot of the combined soft X-ray flux time series from the maxima of solar cycles 
    23 and 24, which correspond to the date ranges 1 January 2000 to 31 December 2002 and 1 September 2011 to 31 May 2014, respectively. The horizontal lines correspond to the 0.9, 0.95, and 0.99 sample quantiles of the combined dataset. Note that the time series is plotted on a logarithmic scale.}}
    \label{fig:goes_flux_time_series}
\end{figure}

We processed the downloaded data as follows. First, we only retained measurements recorded between January 2000 and December 2002 or between September 2011 and May 2014. These two periods roughly represent the maxima of solar cycles 23 and 24, respectively; the flux and flaring activity are elevated during the maximum of a solar cycle. These periods were identified by looking for the peak periods in the sunspot time series displayed in the U.S. Space Weather Prediction Center's Solar Cycle Progression product. To smooth the data, we took the maximum flux over each hour of the day. This produced a time series of length 50,400. Next, we used linear interpolation to impute a small number of missing values, which constituted a little less than 1\% of the values. The final time series is displayed in Figure~\ref{fig:goes_flux_time_series}. Hourly maxima are colored based on which solar cycle they belong to. Note that since we have concatenated two time series separated by several years, observation indices instead of times are shown on the horizontal axis. Figure~\ref{fig:goes_flux_time_series} suggests that it is reasonable to treat the concatenated series as a single time series.

We divided the final time series into windows of size 4,320, which is the number of hours in 180 days. Consecutive windows differed in their starting times by twelve observations, so there were 3,840 windows. Over each window, we fit three different kinds of models, and we predicted exceedances of several thresholds $h$ hours after the last hour in the window for several values of $h$. The first kind of model, which we shall call the baseline model, predicts that a future maximum flux will exceed a threshold if and only if the most recent maximum flux exceeds it. The second kind of model is an AR($d$) model, which we fit using least squares after centering the flux maxima in the training window; the maxima were centered because the flux is positive. We initially tried fitting AR models using LAD estimation, but found that the fitted models were not calibrated as well as models fit using least squares; see Table~\ref{tab:expanded_data_analysis_results} and Figure~\ref{fig:lad_vs_ols} for more details. Predictions were made using the methodology described in Algorithm~\ref{alg:AR_pred} in Section~\ref{subsubsect:AR_approx_opt_pred}. The third kind of model is a FARIMA$(0, d, 0)$ model with symmetric $\alpha$-stable innovations, which was also fit to centered data. We next give necessary background on FARIMA models.

\subsection{Background on FARIMA Models}\label{subsect:farima_background}

The idea to consider fractional autoregressive integrated moving average (FARIMA) models came from our reading of \cite{stanislavsky2009fari}, where such models were applied to the flux data. A FARIMA($p$, $d$, $q$) model is defined by the equation
\begin{equation}\label{eq:farima_p_d_q_mod}
    \phi(B)Y_t = \theta(B)(1 - B)^{-d}\epsilon_t, \quad t \in \Z,
\end{equation}
where $\phi(B)$ and $\theta(B)$ are polynomials of degrees $p$ and $q$, 
respectively, in 
the backshift operator $B$, and $d \in \R$. For simplicity, 
following \cite{kokoszka1995frac}, we shall consider 
iid standard symmetric $\alpha$-stable (S$\alpha$S) innovations with 
$\alpha\in (0,2)$. That is, their characteristic function is 
$\E[ e^{ i \epsilon_t u}] = e^{-|u|^\alpha},\ u\in\R$.
It can be shown that the $\epsilon_t$'s are then heavy-tailed with exponent 
$\alpha$:
\begin{equation}\label{e:epsilon_t-heavy-tailed}
\P[ \epsilon_t<-x] = \P[\epsilon_t>x] \sim c x^{-\alpha},\ \ \mbox{ as }x\to\infty,
\end{equation}
for some $c>0$ \citep[][]{samorodnitsky:taqqu:1994book,nolan2020univ}.  Therefore,
the innovations as well as the resulting FARIMA time series will have infinite variances
and in fact infinite means for $0<\alpha\le 1$.  One can instead consider FARIMA
models with asymmetric regularly varying innovations in the domain of attraction of 
$\alpha$-stable laws \citep[see, e.g.,][]{kokoszka1996para}.

For the sake of simplicity, we consider FARIMA$(0,d,0)$ models, for which \eqref{eq:farima_p_d_q_mod} becomes
\begin{equation}\label{eq:farima_0_d_0_mod}
    Y_t = (1 - B)^{-d}\epsilon_t, \quad t \in \Z.
\end{equation}
\cite{stanislavsky2009fari} found that the flux time series can be modeled reasonably well using such a model. It follows from Theorem 2.1 in \cite{kokoszka1995frac} that if $d$ is non-integral, $\alpha \in (1, 2)$, and $d < 1 - 1 / \alpha$, then the unique solution to \eqref{eq:farima_0_d_0_mod} is $Y_t = \sum_{j = 0}^{\infty} a_j\epsilon_{t - j}$, where the terms of $a = \{a_j\}$ are given by
\begin{equation}\label{eq:farima_causal_coefs}
    a_0 = 1, \quad a_j = \frac{\Gamma(j + d)}{\Gamma(d)\Gamma(j + 1)},
\end{equation}
$\Gamma$ being the gamma function. It is well-known that for $d > 0$,
\[
a_j = \frac{\Gamma(j + d)}{\Gamma(d)\Gamma(j + 1)} \sim \frac{j^{d - 1}}{\Gamma(d)}, \ \ \ \mbox{as $j \to \infty$}
\]
\citep[see, e.g., Lemma 3.1 in][]{kokoszka1995frac}. This implies that $a = \{a_j\}$ 
belongs to $\seqSet$ provided that $d < 1 - 1 / \alpha$. It follows from \eqref{e:epsilon_t-heavy-tailed} 
and Lemma~\ref{lem:xi_a_facts} that the solution $\{Y_t\}$ to \eqref{eq:farima_0_d_0_mod} is finite almost surely; Theorem 2.3 of 
\cite{kokoszka1995frac} implies that the solution is invertible.

\begin{remark}\label{rem:farima_asymp_opt_precs}
    Under the assumptions about $d$ and $\alpha$ given above, for any $h \ge 1$, the extremal optimal precision $\lambda_h^{\rm (opt)}$ can be computed as the ratio $\normConst{+}{a}{h} / \normConst{+}{a}{0}$, according to Theorem~\ref{thm:MA_infty_opt_asymp_precision}. Given $d$, $\alpha$, and the extremal skewness $p_{\epsilon}$ of the innovations, $\lambda_h^{\rm (opt)}$ can be approximated by truncating the series defining $\normConst{+}{a}{h}$ and $\normConst{+}{a}{0}$ using some large number of initial terms. We used this approach to investigate how the extremal optimal precision for a FARIMA$(0, d, 0)$ model varies with $d$ and $\alpha$; see Section~\ref{subsubsect:farima_asymp_opt_precs} for the results.
\end{remark}

\subsection{FARIMA Modeling of the Soft X-Ray Flux}\label{subsect:flux_farima_mods}

Recall from Section~\ref{subsect:data_processing} that the maximum flux was taken over each hour. Let $Y_t$ represent the centered maximum flux in hour $t$; we center the maxima because the flux is positive. We model $\{Y_t\}$ using a FARIMA$(0, d, 0)$ model with S$\alpha$S innovations, i.e., we write
\[
Y_t = (1 - B)^{-d}\epsilon_t, \quad t \in \Z,
\]
where the $\epsilon_t$'s are iid S$\alpha$S random variables. As stated in Section~\ref{subsect:farima_background}, under the assumptions that $d$ is non-integral, $\alpha \in (1, 2)$, and $d \in (0, 1 - 1 / \alpha)$, there is a unique solution, and it satisfies $Y_t = \sum_{j = 0}^{\infty} a_j\epsilon_{t - j}$ for all $t \in \Z$, where the $a_j$'s are defined by the recurrence relation in \eqref{eq:farima_causal_coefs}. Because $\{Y_t\}$ is stationary, the $Y_t$'s are identically distributed. Also, by Lemma~\ref{lem:xi_a_facts}, the $Y_t$'s are regularly varying with tail index $\alpha$. Note that since $\alpha > 1$, it follows from Proposition 1.4.6 in \cite{kulik2020heav} that $\E Y_t$ is finite.

\begin{figure}[!htb]
    \centering
    \includegraphics[width=0.32\textwidth]{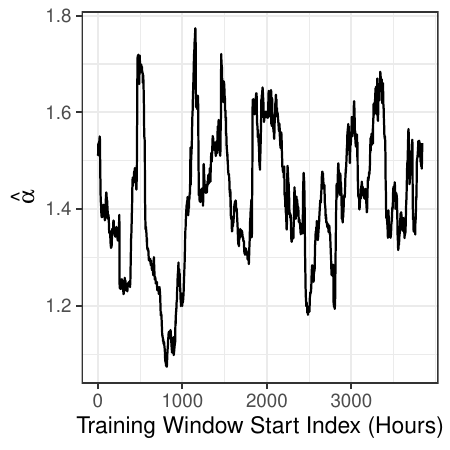}
    \includegraphics[width=0.32\textwidth]{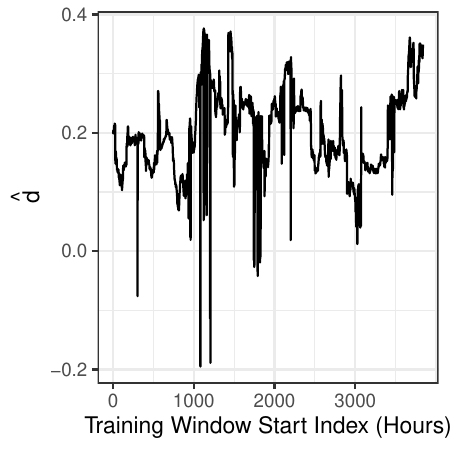}
     \includegraphics[width=0.32\textwidth]{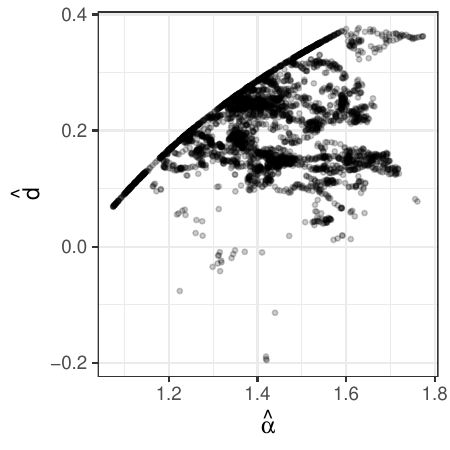}
    \caption{{\small \textbf{Left}: estimates of $\alpha$ for the FARIMA$(0, d, 0)$ models with symmetric $\alpha$-stable innovations that were fit over the training windows. \textbf{Center}: estimates of $d$ for the various models. \textbf{Right}: a scatterplot showing the relationship between $\widehat{\alpha}$ and $\widehat{d}$; for these models, $d$ must be less than $1 - 1 / \alpha$, hence the curve in the upper left.}}
    \label{fig:farima_param_estims}
\end{figure}

Since the uncentered observations are maxima over blocks of a relatively large number of values, sixty, we can approximate the marginal distribution of the centered observations by a GEV (generalized extreme value) distribution. Because of the celebrated 
Fisher-Tippett-Gnedenko theorem and its extensions, a 
GEV distribution is a common choice of model for 
block maxima \citep[][]{leadbetter:lindgren:rootzen:1983,embrechts:kluppelberg:mikosch:1997,beirlant:goegebeur:teugels:segers:2004}. 
The GEV distribution with location 
parameter $\mu$, scale parameter $\sigma > 0$, and shape 
parameter $\xi$ has distribution function 
\[
G(z) = \exp\left\{ - \left[1 + \xi\left(\frac{z - \mu}{\sigma}\right)\right]_+^{-1 / \xi}\right\},
\]
where $0^{-1/\xi} =0,$ for $\xi<0$, $0^{-1/\xi}=\infty$, for $\xi>0$, and 
$G(z) = \exp\{ - e^{-(z-\mu)/\sigma}\}$, for $\xi = 0$ \citep[][]{embrechts:kluppelberg:mikosch:1997,coles2001anin}. Assuming that the observations in a training set are approximately GEV, $\alpha$ can be estimated by fitting a GEV distribution to them and taking $\widehat{\alpha} = 1 / \widehat{\xi}$. We implemented this approach using the R package \texttt{extRemes} of \cite{gilleland2016extr}. In the left panel of Figure~\ref{fig:farima_param_estims}, the estimates of $\alpha$ for the various training windows are plotted against the index of the first observation in the window. There is no apparent trend; the estimates oscillate around roughly 1.4. The variability in $\widehat{\alpha}$ appears to be larger for earlier windows, which mostly consisted of observations from solar cycle 23, than for later windows, which mostly consisted of observations from solar cycle 24.

We estimated $d$ using the method for fitting a FARIMA$(p, d, q)$ model outlined in Appendix D of \cite{burnecki2014algo}; it is based on the procedures described in \cite{kokoszka1996para} and \cite{burnecki2013esti}. A FARIMA$(0, d, 0)$ model can be fit to observations $Y_t, \ldots, Y_{t - n + 1}$ by doing the following:
\begin{enumerate}
    \item Calculate the normalized periodogram
    \begin{equation}
        I_n(\lambda) = \left|\sum_{u = t - n + 1}^t Y_u e^{-i\lambda(u - t + n)}\right|^2, \ -\pi \le \lambda \le \pi.
    \end{equation}
    \item Estimate $d$ as
    \begin{equation}
    \widehat{d} = \mathop{\rm argmin}_{d \in (-1 / 2, 1 - 1 / \alpha)} \int_{1 / n}^{\pi} (2 - 2\cos\lambda)^d I_n(\lambda)\,d\lambda.
    \end{equation}
\end{enumerate}
The estimates of $d$ are shown in the center panel of Figure~\ref{fig:farima_param_estims}. As in the left panel, there is no trend, and the variability looks larger for earlier windows than for later windows.

The right panel of Figure~\ref{fig:farima_param_estims} 
reveals some interesting aspects of the 
relationship between $\widehat{\alpha}$ and $\widehat{d}$. Recall that
$d < 1 - 1 / \alpha$; for many training windows,
$\widehat{d}$ is as large as it is permitted to be, given $\widehat{\alpha}$. This 
suggests that over many windows, the time series exhibits very strong long-range 
dependence. However, this could also be partly due to the presence of complex 
non-stationary trends in the data that can often be confounded with long-range dependence
\citep[][]{stoev:taqqu:park:marron:2005}.

Suppose that we have fit a FARIMA$(0, d, 0)$ model to centered training values $Y_t, \ldots, Y_{t - n + 1}$, with $\widehat{d}$ being the estimate of $d$. Also suppose that we are interested in forecasting exceedance of a threshold at time $t + h$ using the most recent $\ell$ observations, $Y_t, \ldots, Y_{t - \ell + 1}$. The training set size $n$ should be considerably larger than $\ell$. By Remark~\ref{rem:finite_history_pred}, we can do this by thresholding appropriately the quantity $\sum_{r = 0}^{\ell - 1} c_r Y_{t - r}$, where
\begin{equation}\label{eq:farima_opt_pred_coefs}
    c_r = \sum_{s = 0}^r a_{s + h}b_{r - s}, \ r = 0, \ldots, \ell - 1,
\end{equation}
the $a_j$'s are given by \eqref{eq:farima_causal_coefs}, and the $b_j$'s are the coefficients in the invertible representation $\epsilon_t = \sum_{j = 0}^{\infty} b_j Y_{t - j}$, which exists by Theorem 2.3 of \cite{kokoszka1995frac}. The $c_j$'s can be estimated using estimates of the $a_j$'s and the $b_j$'s. The $a_j$'s can be estimated by plugging $\widehat{d}$ into the recurrence relation in \eqref{eq:farima_causal_coefs}, while the $b_j$'s can be estimated using the recursion
\begin{equation}\label{eq:farima_invert_coefs}
    b_j =
    \begin{cases}
      1 & \text{if $j = 0$} \\
      -\sum_{k = 1}^j a_k b_{j - k} & \text{if $j \ge 1$},
    \end{cases}
\end{equation}
with $\widehat{a}_k$ in the place of $a_k$. From \eqref{eq:farima_opt_pred_coefs} and \eqref{eq:farima_invert_coefs}, to estimate $c_0, \ldots, c_{\ell - 1}$, we need only estimate $a_0, \ldots, a_{\ell - 1 + h}$ and $b_0, \ldots, b_{\ell - 1}$. Suppose that we want to predict whether the uncentered observation $X_{t + h}$ will exceed its $p$th marginal quantile. Once $\widehat{c}_0, \ldots, \widehat{c}_{\ell - 1}$ have been computed, we can compute $\sum_{r = 0}^{\ell - 1} \widehat{c}_r Y_{t - r}, \ldots, \sum_{r = 0}^{\ell - 1} \widehat{c}_r Y_{t - n + \ell - r}$; let $q$ be the $p$th sample quantile of these sums. We predict that an exceedance will occur if and only if $\sum_{r = 0}^{\ell - 1} \widehat{c}_r Y_{t - r} \ge q$. However, it is not hard to show that the $a_j$'s form a decreasing sequence. It then follows from Proposition~\ref{prop:oracle} that if the model is correctly specified, then the optimal predictor has the same asymptotic precision as a predictor that merely thresholds the most recent observation. Hence, if the model is reasonable, then we can expect the baseline predictor to perform at least as well as the FARIMA predictor when the threshold is a high quantile.

\subsection{Prediction Results}\label{subsect:pred_results}

Table~\ref{tab:data_analysis_results} summarizes the results from our prediction 
experiments. We sought to predict exceedances of the $p$th marginal quantile $h$ steps 
ahead, for $p = 0.90, 0.95, 0.99$ and $h = 1$ and $6$ (see also 
Table~\ref{tab:expanded_data_analysis_results} for $h=12$ and $18$). 
The number of observations used to compute predictions from FARIMA models was 168, the number of hours in one week; this value was also used as the order of the AR models that we fit. The True Skill Statistic (TSS) is commonly used to evaluate solar flare forecasting methods; see Remark~\ref{rem:metric_equivs} for more information about it. Since the patterns in the TSS section of Table~\ref{tab:data_analysis_results} are very similar to those in the precision section, we shall only comment on the latter. In each row, in the section for each metric, the largest number is in bold.

\begin{table}[ht]
\centering
\begin{tabular}{cc|ccc|ccc}
  \hline
  & & \multicolumn{3}{c|}{Precision} & \multicolumn{3}{c}{TSS} \\
 $h$ & $p$ & Baseline & FARIMA & AR(168) & Baseline & FARIMA & AR(168) \\
 \hline
 1 & 0.90 & {\bf 0.492} & 0.448 & 0.354 & {\bf 0.435} & 0.392 & 0.276 \\ 
   1 & 0.95 & {\bf 0.379} & 0.368 & 0.258 & {\bf 0.334} & 0.320 & 0.225 \\ 
   1 & 0.99 & {\bf 0.244} & 0.238 & 0.104 & {\bf 0.262} & {\bf 0.262} & 0.124 \\ 
   \hline
 6 & 0.90 & 0.308 & {\bf 0.339} & 0.284 & 0.256 & {\bf 0.298} & 0.215 \\ 
   6 & 0.95 & 0.189 & {\bf 0.213} & 0.150 & 0.173 & {\bf 0.198} & 0.138 \\ 
   6 & 0.99 & 0.049 & 0.119 & {\bf 0.125} & 0.040 & 0.115 & {\bf 0.162} \\ 
   \hline
\end{tabular}
\caption{A summary of the performance of various predictors. The AR models were fit using OLS. For each pair of $h$ and $p$, the goal was to predict exceedance of the $p$th marginal quantile $h$ steps ahead. In each row, for each metric, the largest value is in boldface. For additional results, including results for AR models fit using LAD, see Table \ref{tab:expanded_data_analysis_results}.}
\label{tab:data_analysis_results}
\end{table}

It is certain that complex solar flaring phenomena are
not governed by a simple model such as a FARIMA$(0,d,0)$ or 
an AR$(d)$.  Nevertheless, we shall discuss the prediction performance of
these models from the perspective of our theoretical results.  Specifically,
the qualitative results in Proposition~\ref{prop:oracle} suggest that for a large
class of MA$(\infty)$ models with absolutely monotone coefficients, the baseline
predictor will be asymptotically optimal for large $p$.

The baseline predictor, which is based on just the most recent observation, outperforms a high-order AR$(168)$ model in all of the cases except $h=6$ and $p=0.99$; it {\em significantly} outperforms the AR model when $h = 1$. The baseline also outperforms FARIMA$(0,d,0)$ for $h=1$, though it is roughly comparable in precision, while FARIMA$(0,d,0)$ outperforms slightly the baseline for $h=6$ and $p=0.9,0.95$. 
Since the moving average representation of FARIMA$(0,d,0)$ involves
a monotone sequence of coefficients, the baseline and FARIMA predictors 
should have had similar performance, especially at extreme thresholds. The discrepancy between their precisions is relatively small, which is somewhat in line with the qualitative results in Proposition~\ref{prop:oracle}, though it suggests that simple absolutely
monotone MA$(\infty)$ models may not be appropriate for prediction over longer lags and extreme levels.

The FARIMA and baseline predictors outperform the AR$(168)$ predictor in the majority of the cases, except
for $h=6$ and $p=0.99$.  This remains true over much longer time horizons, 
$h=12$ and $18$,
as seen in the extended Table \ref{tab:expanded_data_analysis_results}. 
These results suggest that the baseline and long-range dependent models are 
most powerful for prediction of high fluxes over moderate to extreme 
levels $p=0.9$ and $p=0.95$. On the other hand, higher order short-range 
dependent models prevail over intermediate horizons $(h=6)$ at the most extreme levels 
$p=0.99$, but AR is significantly outperformed by baseline or FARIMA when $p = 0.99$ 
and $h = 12, 18$. One can speculate
that this is due to fundamentally different physical phenomena governing the 
very extreme flares $(p=0.99)$, which although persisting over several hours
are likely to be more isolated in time than the moderate flaring activity 
($p\le 0.95$) that could be the result of longer term processes persisting
and generating flares over several days.  This qualitative hypothesis is to be checked
with extensive analysis involving additional (image) covariates of the Sun.

One quantitative contribution of the above analysis is a benchmark 
lower bound for the optimal prediction precision of high fluxes from historical 
GOES time series data.  It shows that prediction of moderately extreme fluxes
with a reasonable precision of roughly 0.2-0.3 is possible over time horizons on the 
order of several hours to a day.  We anticipate that this precision can be improved, but
only by incorporating more informative covariates, which we plan to do in our future
work.

\section{Supporting Information}

Additional supporting information can be found online in the supporting information tab for this article.

\section{Data Availability Statement}

Code for the simulations and data analysis and the data for the latter can be accessed at \url{https://github.com/victorverma/opt_extreme_event_pred}. The Solar Cycle Progression product mentioned in Section~\ref{sect:data_analysis} is located at \url{https://www.swpc.noaa.gov/products/solar-cycle-progression}.

\section{Acknowledgements}

We thank Dr.~Janet Machol and Prof.\ Ward Manchester for their advice on processing the GOES data. We acknowledge NOAA and the NCEI for making that data available. SS was partially supported by the NSF grant 2319592. YC and VV were partially supported by NSF DMS 2113397, NSF PHY 2027555, NASA Federal Award No. 80NSSC23M0192 and No. 80NSSC23M0191.

\bibliographystyle{agsm}
\bibliography{bibliography}


\appendix
\newpage
\renewcommand{\thesection}{A}
\setcounter{section}{0}

\renewcommand{\theHsection}{A\arabic{section}}

\section{Appendix}\label{sect:appendix}

\subsection{Additional Information for Section~\ref{sect:foundations}}

\subsubsection{Proof of Lemma \ref{lem:calibration_levels}}\label{subsubsect:calibration_levels_lem_pf}

\begin{proof}
    {\em Part (i):} Let $p \in [0, 1]$, and suppose that the pre-image $F_{\xi}^{-1}(\{p\})$ contains two distinct points $x$ and $x'$, with $x < x'$. Then for all $x'' \in (x, x')$,
    \[
    p = F_{\xi}(x) \le F_{\xi}(x'') \le F_{\xi}(x') = p,
    \]
    so $x'' \in F_{\xi}^{-1}(\{p\})$, which means that $F_{\xi}^{-1}(\{p\})$ is an interval. When the cardinality of $F_{\xi}^{-1}(\{p\})$ is bigger than one, we shall call $F_{\xi}^{-1}(\{p\})$ a \textit{flat spot} since $F$ is flat over it. Since flat spots corresponding to distinct values of $p$ must be disjoint, we can uniquely identify a flat spot by a rational number inside it. Hence, the collection of flat spots is countable; let $\{F_{\xi}^{-1}(\{p_i\}) : i \in I\}$ be an enumeration of the flat spots, with $I \subseteq \N$.

    Suppose that for some $x \in \R$, $\Fi{\xi}(F_{\xi}(x)) \ne x$. Then $\Fi{\xi}(F_{\xi}(x)) < x$. For any $x' \in (\Fi{\xi}(F_{\xi}(x)), x)$, we must have $F_{\xi}(x') = F_{\xi}(x)$, so $F_{\xi}^{-1}(\{F_{\xi}(x)\})$ is a flat spot, and $x \in F_{\xi}^{-1}(\{F_{\xi}(x)\}) \setminus \{\inf F_{\xi}^{-1}(\{F_{\xi}(x)\})\}$. Thus,
    \begin{equation}\label{eq:x_in_union}
        x \in \bigcup_{i \in I} F_{\xi}^{-1}(\{p_i\}) \setminus \{\inf F_{\xi}^{-1}(\{p_i\})\}.
    \end{equation}
    Conversely, if \eqref{eq:x_in_union} holds, then for some $i$, $x \in F_{\xi}^{-1}(\{p_i\}) \setminus \{\inf F_{\xi}^{-1}(\{p_i\})\}$, so $\Fi{\xi}(F_{\xi}(x)) = \Fi{\xi}(p_i) < x$. We deduce that
    \[
    \{\Fi{\xi}(F_{\xi}(X)) \ne X\} = \left\{X \in \bigcup_{i \in I} F_{\xi}^{-1}(\{p_i\}) \setminus \{\inf F_{\xi}^{-1}(\{p_i\})\}\right\}.
    \]
    The probability of the latter event is bounded above by
    \[
    \sum_{i \in I} \P(X \in F_{\xi}^{-1}(\{p_i\}) \setminus \{\inf F_{\xi}^{-1}(\{p_i\})\}).
    \]
    The claim will follow if we can show that each probability in this sum equals zero.

    Observe that by the right continuity of $F_{\xi}$, $F_{\xi}(\inf F_{\xi}^{-1}(\{p_i\})) = p_i$. If $F_{\xi}^{-1}(\{p_i\})$ is closed on the right, then
    \[
    \P(X \in F_{\xi}^{-1}(\{p_i\}) \setminus \{\inf F_{\xi}^{-1}(\{p_i\})\})
    = F_{\xi}(\sup F_{\xi}^{-1}(\{p_i\})) - F_{\xi}(\inf F_{\xi}^{-1}(\{p_i\}))
    = p_i - p_i
    = 0.
    \]
    If $F_{\xi}^{-1}(\{p_i\})$ is open on the right, then
    \begin{align*}
        \P(X \in F_{\xi}^{-1}(\{p_i\}) \setminus \{\inf F_{\xi}^{-1}(\{p_i\})\})
        &= \P(X < \sup F_{\xi}^{-1}(\{p_i\})) - F_{\xi}(\inf F_{\xi}^{-1}(\{p_i\})) \\
        &= \lim_{x\,\uparrow\,\sup F_{\xi}^{-1}(\{p_i\})} F_{\xi}(x) - p_i \\
        &= p_i - p_i \\
        &= 0.
    \end{align*}
    The claim now follows.
    
    {\em Part (ii):} Let $y_n$ be such that $y_n\downarrow F_{\xi}^\leftarrow(q)$.  By the right continuity of $F_{\xi}$, we get 
    \begin{equation}\label{e:lem:calibration_levels-1}
    F_{\xi}( F_{\xi}^\leftarrow(q) ) = \lim_{n\to\infty} F_{\xi}(y_n) \ge q.
    \end{equation}
    Since $q\in {\rm Range}(F_\xi)$, there exists a $\tau$ such that
    $F_\xi(\tau) = q$. Thus,  $\tau\in \{y : F_{\xi}(y) \ge q\}$ and in view of
    \eqref{e:gen-inverse}, we have $\tau \ge F_{\xi}^\leftarrow(q)$. This implies that
     $F_{\xi}( F_{\xi}^\leftarrow(q) )\le  F_{\xi}(\tau) =q$ and combined with \eqref{e:lem:calibration_levels-1},
     yields $\P[ \xi > F_{\xi}^\leftarrow(q) ] = 1-q,$
     completing the proof of \eqref{e:exact-calibration}.
\end{proof}

\subsubsection{More Examples of Explicit Optimal Predictors}\label{subsect:more_example_opt_preds}

For two additional models, we derive an expression for an optimal extreme event predictor that is not based on the density ratio $r$ appearing in Theorem~\ref{thm:base_thm}.

\begin{proposition}\label{prop:non-linear}
    Let $Y = \varphi_{g(X)}(\epsilon)$, where $X \ind \epsilon$, $X$ has a density $f_X$ with respect to a $\sigma$-finite Borel measure $\mu$ on $\R$, and $\epsilon$ has a density $f_\epsilon$ with respect to a $\sigma$-finite Borel measure $\nu$ on $\R$. Assume that for all $u \in \R$, the function $\varphi_u(\cdot)$ is strictly increasing and smooth and that $\lim_{v \to \pm\infty} \varphi_u(v) = \pm\infty$. Assume moreover that for all $v \in \R$, $\varphi_{u_1}(v) \le \varphi_{u_2}(v)$, for all $u_1 \le u_2$, and that $u \mapsto \varphi_u^{-1}(v)$ is smooth, where, for a fixed $u$, $\varphi_u^{-1}$ is the inverse of $\varphi_u(\cdot)$. Finally, assume that $\nu \circ \varphi_{g(x)}^{-1}$ is absolutely continuous with respect to $\nu$ for all $x \in \R$. Let $p \in (0, 1)$ and $y_0 \in \R$ such that $\P(Y \le y_0) = p$. Then $\I(g(X) > \tau)$ is an optimal predictor of $\I(Y > y_0)$, where $\tau$ is a constant chosen to calibrate the predictor.
\end{proposition}
\begin{proof}
    Note that for a fixed $u$, $\varphi_u^{-1}$ exists and is strictly increasing. Furthermore, for any $t$, $v > \varphi_u(t)$ and $v < \varphi_u(t)$ imply that $\varphi_u^{-1}(v) > t$ and $\varphi_u^{-1}(v) < t$, respectively; these in turn imply that $\lim_{v \to \pm\infty} \varphi_u^{-1}(v) = \pm\infty$.

    For all Borel subsets $A$ and $B$ of $\R$,
    \begin{align}
        \P((X, Y) \in A \times B)
        &= \E[\E[\I(X \in A, Y \in B) \mid X]] \nonumber \\
        &= \int f_X(x)\E[\I(X \in A)\I(\varphi_{g(X)}(\epsilon) \in B) \mid X = x]\,\mu(dx) \nonumber \\
        &= \int_A f_X(x)\P(\varphi_{g(x)}(\epsilon) \in B)\,\mu(dx) \nonumber \\
        &= \int_A f_X(x)\P(\epsilon \in \varphi_{g(x)}^{-1}(B))\,\mu(dx) \nonumber \\
        &= \int_A f_X(x)\left[\int_{\varphi_{g(x)}^{-1}(B)} f_{\epsilon}(v)\,\nu(dv)\right]\,\mu(dx) \nonumber \\
        &= \int_A f_X(x)\left[\int_B f_{\epsilon}(\varphi_{g(x)}^{-1}(y))\,(\nu \circ \varphi_{g(x)}^{-1})(dy)\right]\,\mu(dx) \label{eq:iterated_integral}
    \end{align}
    Note that for each $x$, $\nu \circ \varphi_{g(x)}^{-1}$ is $\sigma$-finite. By assumption, $\nu$ is $\sigma$-finite, so there exist Borel sets $B_1, B_2, B_3, \ldots$ that cover $\R$ and have finite $\nu$-measure. For each $i$, $(\nu \circ \varphi_{g(x)}^{-1})(\varphi_{g(x)}(B_i)) = \nu(B_i) < \infty$, so $\varphi_{g(x)}(B_1), \varphi_{g(x)}(B_2), \varphi_{g(x)}(B_3), \ldots$ have finite measure with respect to $\nu \circ \varphi_{g(x)}^{-1}$. Furthermore, since $\varphi_{g(x)}$ is surjective, those images cover $\R$. It follows that $\nu \circ \varphi_{g(x)}^{-1}$ is $\sigma$-finite. Let $\delta_x$ be the Radon-Nikodym derivative of $\nu \circ \varphi_{g(x)}^{-1}$ with respect to $\nu$. Then \eqref{eq:iterated_integral} becomes
    \begin{align}
        \P((X, Y) \in A \times B)
        &= \int_A f_X(x)\left[\int_B f_{\epsilon}(\varphi_{g(x)}^{-1}(y))\delta_x(y)\,\nu(dy)\right]\,\mu(dx) \nonumber \\
        &= \int_{A \times B} f_X(x)f_{\epsilon}(\varphi_{g(x)}^{-1}(y))\delta_x(y)\,(\mu \times \nu)(dx, dy) \label{eq:integrated_density}.
    \end{align}
    It follows that the integrand in \eqref{eq:integrated_density} is the joint density of $X$ and $Y$ with respect to $\mu \times \nu$. Similar calculations show that $f_0$, the conditional density of $X$ given that $Y \le y_0$, and $f_1$, the conditional density of $X$ given that $Y > y_0$, are given by
    \begin{align*}
        f_0(x) &= \frac{1}{p}\int_{(-\infty, y_0]} f(x, y)\,\nu(dy) = \frac{f_X(x)}{p}F_{\epsilon}(\varphi_{g(x)}^{-1}(y_0)) \\
        f_1(x) &= \frac{1}{1 - p}\int_{(y_0, \infty)} f(x, y)\,\nu(dy) = \frac{f_X(x)}{1 - p}\overline{F}_{\epsilon}(\varphi_{g(x)}^{-1}(y_0)),
    \end{align*}
    where $F_{\epsilon}$ is the distribution function of $\epsilon$ and $\overline{F}_{\epsilon} = 1 - F_{\epsilon}$ is the tail function of $\epsilon$. We compute the ratio $r$ of $f_1$ to $f_0$:
    \begin{equation}\label{eq:non-linear_density_ratio}
        r(x) = \left(\frac{p}{1 - p}\right)\frac{\overline{F}_{\epsilon}(\varphi_{g(x)}^{-1}(y_0))}{F_{\epsilon}(\varphi_{g(x)}^{-1}(y_0))},
    \end{equation}
    By Theorem~\ref{thm:base_thm}, for any $q \in \Range(F_{r(X)}) \cap (0, 1)$, $\I(r(X) > \Fi{r(X)}(q))$ is an optimal predictor of $\I(Y > y_0)$ calibrated at level $q$. We have $r(X) > \Fi{r(X)}(q)$ if and only if
    \begin{equation}\label{e:prop-non-linear}
        \left(\frac{p}{1 - p}\right)\frac{\overline{F}_{\epsilon}(\varphi_{g(X)}^{-1}(y_0))}{F_{\epsilon}(\varphi_{g(X)}^{-1}(y_0))}
        > \Fi{r(X)}(q).
    \end{equation}
    Let $U = \varphi_{g(X)}^{-1}(y_0)$ and let $\tilde{r} : \R \to \R$ be the function that maps $u$ to $\left(\frac{p}{1 - p}\right)\frac{\overline{F}_{\epsilon}(u)}{F_{\epsilon}(u)}$. Then the optimal prediction condition \eqref{e:prop-non-linear} can be written as $\tilde{r}(U) > \Fi{r(X)}(q)$. Since $\tilde{r}$ is decreasing and right-continuous, it follows from part (ii) of Lemma~\ref{lem:opt_pred_cond_simplifier} that this condition is equivalent to the condition $\varphi_{g(X)}^{-1}(y_0) = U < \rho$ for some real number $\rho$.

    Now, when $u_1 \le u_2$,
    \[
    y_0 = \varphi_{u_1}( \varphi_{u_1}^{-1} (y_0)  )\le \varphi_{u_2}( \varphi_{u_1}^{-1} (y_0)  ).
    \]
    On the other hand, $\varphi_{u_2}( \varphi_{u_2}^{-1} (y_0)  ) =y_0$, which by the fact that $\varphi_{u_2}(\cdot)$ is
    strictly increasing implies that $\varphi_{u_1}^{-1} (y_0) \ge \varphi_{u_2}^{-1} (y_0)$. That is, for a fixed
    $y_0$, we have shown that the function $\psi : \R \to \R$ defined by $\psi(u) = \varphi_u^{-1}(y_0)$ is decreasing. The condition $\varphi_{g(X)}^{-1}(y_0) < \rho$ can be written as $\psi(g(X)) < \rho$, or $-\psi(g(X)) > -\rho$. By assumption, $\psi$ is smooth. It follows from part (i) of Lemma~\ref{lem:opt_pred_cond_simplifier} that $-\psi(g(X)) > -\rho$ is equivalent to $g(X) > \tau$ for some $\tau$, which means that $\I(g(X) > \tau)$ is an optimal predictor of $\I(Y > y_0)$.
\end{proof}

\begin{remark}
    Suppose that $Y = g(X) + \epsilon$, where $X \ind \epsilon$, $X$ has a density with respect to some $\sigma$-finite Borel measure on $\R$, and $\epsilon$ has a density with respect to the Lebesgue measure on $\R$. Then $Y = \varphi_{g(X)}(\epsilon)$, where $\varphi_u(v) = u + v$. All of the assumptions of Proposition~\ref{prop:non-linear} hold, so given a threshold $y_0$, $\I(g(X) > \tau)$ is an optimal predictor of $\I(Y > y_0)$ for some appropriately chosen $\tau$. Note that this result is similar to the corollary of Theorem~\ref{thm:add_err_mod_thm} given in Remark~\ref{rem:homoscedastic_noise}, the difference being that here $X$ is a random variable instead of a random vector.
\end{remark}

\begin{remark} Another example illustrating Proposition \ref{prop:non-linear} is when
$\phi_u(v) = G(u, v)$ is a joint distribution function with strictly positive joint density $g(u,v) = \partial^2_{u,v} G(u,v)$.
\end{remark}

\begin{remark}
    Recall that random variables $X$ and $Y$ are positively dependent if and only if $\Cov(g(X), h(Y)) \ge 0$ for all $g, h \in L^2(X)$ that are increasing. We have not derived the form of the optimal predictor when $X$ and $Y$ are positively dependent. However, we conjecture that the positive dependence of $X$ and $Y$ implies that $Y = k(X, \epsilon)$ for some function $k$ and some random variable $\epsilon$ that is independent of $X$, in which case Proposition~\ref{prop:non-linear} could be used to derive an optimal predictor of $\I(Y > y_0)$ in terms of $X$, under some additional assumptions that make the proposition applicable.
\end{remark}

The next proposition can be used to derive an optimal predictor for a class of state space models.

\begin{proposition}\label{prop:monotone-link}
    Let $Y = g(X + \delta) + \epsilon$, where $X$, $\delta$, and $\epsilon$ are independent random variables. Suppose that $X$ has a density $f_X$ with respect to the Lebesgue measure, that $\delta$ has a density $f_{\delta}$ with respect to a $\sigma$-finite Borel measure $\mu$, and that $\epsilon$ has a density $f_{\epsilon}$ with respect to Lebesgue measure. Let $p \in (0, 1)$ and $y_0 \in \R$ such that $\P(Y \le y_0) = p$. If $g$ is a increasing function, then an optimal predictor of $\I(Y > y_0)$ can be written in the form $\I(X > x_0)$ for some constant $x_0$ that calibrates the predictor.
\end{proposition}
\begin{proof}
    We first compute the joint density of $X$ and $Y$. Let $\lambda$ and $\nu$ represent Lebesgue measure. For all Borel subsets $A$ and $B$ of $\R$,
    \begin{align}
        \P((X, Y) \in A \times B)
        &= \E[\E[\I(X \in A, Y \in B) \mid X, \delta]] \nonumber \\
        &= \int_{\R}\int_{\R} f_X(x)f_{\delta}(u)\E[\I(X \in A, Y \in B) \mid X = x, \delta = u]\,\mu(du)\lambda(dx) \nonumber \\
        &= \int_{\R}\int_{\R} f_X(x)f_{\delta}(u)\I(x \in A)\P(\epsilon \in B - g(x + u))\,\mu(du)\lambda(dx) \nonumber \\
        &= \int_A\int_{\R} f_X(x)f_{\delta}(u)\left[\int_{B - g(x + u)} f_{\epsilon}(v)\,\nu(dv)\right]\,\mu(du)\lambda(dx) \label{eq:triple_integral1} \\
        &= \int_A\int_{\R} f_X(x)f_{\delta}(u)\left[\int_B f_{\epsilon}(y - g(x + u))\,\nu(dy)\right]\,\mu(du)\lambda(dx) \label{eq:triple_integral2} \\
        &= \int_{A \times B}\left[\int_{\R} f_X(x)f_{\epsilon}(y - g(x + u))f_{\delta}(u)\,\mu(du)\right](\lambda \times \nu)(dx, dy), \nonumber
    \end{align}
    where the transition from \eqref{eq:triple_integral1} to \eqref{eq:triple_integral2} is based on the translation invariance of $\nu$ and a change of variable. It follows that the joint density $f$ of $X$ and $Y$ with respect to $\lambda \times \nu$ is the expression in brackets on the last line. The conditional densities $f_0$ and $f_1$ of $X$ given $Y \le y_0$ and $Y > y_0$,
    respectively, can be obtained from $f$:
    \begin{align*}
        f_0(x)
        &= \frac{1}{p}\int_{(-\infty, y_0]} f(x, y)\,\nu(dy)
        = \frac{f_X(x)}{p}\int_{(-\infty, y_0]}\int_{\R} f_{\epsilon}(y - g(x + u))f_{\delta}(u)\,\mu(du)\nu(dy) \\
        f_1(x)
        &= \frac{1}{1 - p}\int_{(y_0, \infty)} f(x, y)\,\nu(dy)
        = \frac{f_X(x)}{1 - p}\int_{(y_0, \infty)}\int_{\R} f_{\epsilon}(y - g(x + u))f_{\delta}(u)\,\mu(du)\nu(dy)
    \end{align*}
    The ratio of conditional densities $r(x)$ is
    \begin{align}
        r(x)
        &= \frac{f_1(x)}{f_0(x)} \nonumber \\
        &= \left(\frac{p}{1 - p}\right)\frac{\int_{(y_0, \infty)} \int_{\R}   f_\epsilon(y-g(x+u))f_\delta(u)\,\mu(du)\nu(dy)} 
        {\int_{(-\infty, y_0]} \int_{\R} f_\epsilon(y-g(x+u))f_\delta(u)\,\mu(du)\nu(dy)} \nonumber \\
        &= \left(\frac{p}{1 - p}\right)\frac{\int_{\R} \overline F_\epsilon (y_0 - g(x+u)) f_\delta(u)\,\mu(du)}{ \int_{\R} F_\epsilon (y_0 - g(x+u))f_\delta(u)\,\mu(du)} \label{eq:state_space_r_expr},
    \end{align}
    with $\overline{F}_\epsilon = 1 - F_{\epsilon}$ being the tail function of $\epsilon$. It follows from Theorem~\ref{thm:base_thm} and Lemma~\ref{lem:calibration_levels} that an optimal predictor of $\I(Y > y_0)$ at level $q$ is $\I(r(X) > \Fi{r(X)}(q))$, for any $q \in \Range(F_{r(X)}) \cap (0, 1)$.
    
    Notice that because $g$ is non-decreasing, for every $u \in \R$, we have that the integrand in the numerator in \eqref{eq:state_space_r_expr} is a non-decreasing function of $x$ and at the same time the integrand in the denominator
    is a non-increasing function of $x$. This implies that $r$ is non-decreasing. Let
    \[
    x_0 = \sup \{x \in \R : r(x) \le \Fi{r(X)}(q)\}.
    \]
    By reasoning similar to that in the proof of Lemma~\ref{lem:opt_pred_cond_simplifier},
    \[
    \{X > x_0\} \subseteq \{r(X) > \Fi{r(X)}(q)\} \subseteq \{X \ge x_0\}.
    \]
    Unlike in that proof, we cannot say that $\{r(X) > \Fi{r(X)}(q)\} = \{X > x_0\}$, because we do not know whether $r(x_0) \le \Fi{r(X)}(q)$. However, since $X$ has a density with respect to Lebesgue measure, $\P(X = x_0) = 0$. It follows that $\I(X > x_0)$ is an optimal predictor of $\I(Y > y_0)$ at level $q$.
\end{proof}

\begin{remark}
    Consider the state space model:
    \begin{align}\label{e:state-space-X}
     X_t  &= \sum_{j=1}^q A_j X_{t-j} + \delta_t\\
     Y_t  &= g \left( \sum_{j=0}^{r - 1} \beta_j^\top X_{t-j} \right) + \epsilon_t, \label{e:state-space-Y}
    \end{align}
    where the $X_t$'s and $\beta_j$'s take values in $\R^d$ and the $A_j$'s are $d\times d$ coefficient matrices.  
    The function $g :\R \to\R$ is assumed to be increasing. We assume that the $\delta_t$'s are iid, that the $\epsilon_t$'s are iid, and that the $\delta_t$'s and the $\epsilon_t$'s are independent of each other. We shall view $X_t$ as the covariate or the state variable, and $Y_t$ as the target response variable. The goal is to predict $\I(Y_{t+h} > y_0)$ via $X_s$, $s\le t$, where $y_0$ satisfies $\P(Y_{t + h} \le y_0) = p \in (0, 1)$.
    
    It can be assumed that $q=r$. If $q < r$, then the sum in \eqref{e:state-space-X} could be padded with extra terms $A_j X_{t - j}$ such that $A_j = 0$. If $r < q$, then the sum in \eqref{e:state-space-Y} could be padded in a similar fashion. Furthermore, we may assume that $q=r=1$, though at the expense of increasing the dimension of the state space. Consider the state space model
    \begin{equation}\label{e:state-space-simplified}
    \tilde X_t = A \tilde X_{t-1} + \tilde \delta_t,
    \quad
    Y_t = g(\beta^{\top}\tilde X_t) + \epsilon_t,
    \end{equation}
    where $\tilde X_t^\top := (X_t^\top, \ldots, X_{t-q+1}^\top)^\top$ takes values in $\R^{q d}$.  Defining 
    $$
    A := \left(\begin{array}{lllll}
        A_1 & A_2 & \cdots & A_{q-1} & A_q\\
        I_d & 0 & \cdots & 0 & 0 \\
        0  & I_d & \cdots & 0& 0\\
        \vdots & \vdots &  \ddots &\vdots & \vdots\\
        0 & 0 & \cdots & I_d & 0
        \end{array}\right),\ \ \tilde \delta_t := \left(\begin{array}{l} \delta_t \\ 0\\ \vdots \\ 0
        \end{array}\right),  \ \ \mbox{ and }\ \ \beta := \left(\begin{array}{l} \beta_1 \\ \beta_2\\ \vdots \\ \beta_q
        \end{array}\right) 
       $$
    we obtain that a solution $\{Y_t\}$ to the simplified state space model \eqref{e:state-space-simplified} of dimensions $q=r=1$ yields a solution to the model in \eqref{e:state-space-X} and \eqref{e:state-space-Y}.
    
    Relation \eqref{e:state-space-simplified} implies
    \begin{align}\label{e:Y-state-space-decomposition}
        Y_{t + h} & = g \Big( \beta^\top ( A \tilde X_{t+h-1} + \tilde  \delta_{t+h})\Big ) + \epsilon_{t+h} \nonumber \\
        & = g \Big( \beta^\top A \tilde X_{t + h-1} +  \beta^\top \tilde  \delta_{t+h}  \Big) + \epsilon_{t+h} \nonumber \\
        & = \cdots \nonumber \\
        & = g \Big( \beta^\top A^h \tilde X_t + \sum_{j=1}^h \beta^\top A^{h - j}\tilde  \delta_{t+j}  \Big) +  \epsilon_{t+h} \nonumber\\
        & = g \Big( \beta^\top A^h \tilde X_t + \delta_{t,h} \Big) + \epsilon_{t+h},
    \end{align}
    where 
    \begin{equation}\label{e:eta_k,h}
        \delta_{t,h}:= \sum_{j=1}^h \beta^\top A^{h - j}\tilde  \delta_{t+j}.
    \end{equation}
    
    Assuming that the innovations $\tilde{\delta}_s,\ s\ge t+1$ and $\epsilon_s,\ s\ge t+1$ are independent from $\tilde  X_u,\ u \le t$, and assuming that $\beta^\top A^h \tilde X_t$, $\delta_{t,h}$ and $\epsilon_{t+h}$ have densities with respect to the appropriate measures, then it follows from Proposition~\ref{prop:monotone-link} that $\I(\beta^\top A^h \tilde X_t > x_0)$ is an optimal predictor of $\I(Y > y_0)$ for some $x_0$.
\end{remark}

\subsubsection{ROC Dominance Property}\label{subsect:roc_domin_prop}

We continue with the setup of Section~\ref{subsect:general_opt_preds}, where the objective is to predict $\I(Y > y_0)$, with $\P(Y > y_0) = 1 - p$ for some $p \in (0, 1)$, using a predictor $\I(h(X) > \tau)$ that is calibrated at a level $q \in \Range(F_{h(X)}) \cap (0, 1)$, i.e., one for which $\P(h(X) > \tau) = 1 - q$. In this section, we define the concept of a family of predictors and show that a family of optimal predictors has a certain dominance property that is based on ROC curves. This property motivates an optimality condition for families of predictors; under this condition, a family of optimal predictors can be viewed as optimal.

Since an ROC curve shows how the true positive rate varies as the false positive rate changes \citep[][]{krzanowski2009rocc}, we first define these quantities. We then express them in terms of $p$, $q$, and the precision $\lambda = \P(Y > y_0 \mid h(X) > \tau)$.
\begin{definition}[True and False Positive Rates]\label{def:tpr_fpr}
    Let $\I(h(X) > \tau)$ be a predictor of $\I(Y > y_0)$.
    
    (i) The true positive rate (TPR) of $\I(h(X) > \tau)$ is $\P(h(X) > \tau \mid Y > y_0)$.

    (ii) The false positive rate (FPR) of $\I(h(X) > \tau)$ is $\P(h(X) > \tau \mid Y \le y_0)$.
\end{definition}

\begin{remark}\label{rem:sample_metrics}
    The performance of a binary classifier on a dataset can be summarized by a {\em confusion matrix} like that in Table~\ref{tab:conf_mat}; the two classes are represented by + and --, and $\TP$, $\FP$, $\FN$, and $\TN$ are the counts of true positives, false positives, false negatives, and true negatives, respectively. For example, $\TP$ is the number of observations belonging to the positive class that were predicted to belong to that class.
    \begin{table}[h]
        \centering
        \begin{tabular}{@{}cc|cc@{}}
            \multicolumn{1}{c}{} &\multicolumn{1}{c}{} &\multicolumn{2}{c}{Observed} \\ 
            \multicolumn{1}{c}{} & 
            \multicolumn{1}{c|}{} & 
            \multicolumn{1}{c}{+} & 
            \multicolumn{1}{c}{--} \\ 
            \cline{2-4}
            \multirow[c]{2}{*}{Predicted} 
            & + & $\TP$ & $\FP$   \\[1.5ex]
            & -- & $\FN$ & $\TN$ \\ 
            \cline{2-4}
            \end{tabular}
        \caption{A confusion matrix summarizing the performance of a binary classifier on a dataset; the two classes are represented by + and --, and ${\rm TP}$, ${\rm FP}$, ${\rm FN}$, and ${\rm TN}$ are the counts of true positives, false positives, false negatives, and true negatives, respectively.}
        \label{tab:conf_mat}
    \end{table}
    
    The TPR, FPR, and precision are typically defined in terms of the confusion matrix as follows \citep[][]{fawcett2006anin}:
    \begin{equation}\label{eq:sample_metrics}
        \TPR := \frac{\text{TP}}{TP + FN}, \ \FPR := \frac{FP}{FP + TN}, \ \text{precision} := \frac{TP}{TP + FP}.        
    \end{equation}
    In \eqref{eq:sample_metrics}, the TPR is the sample proportion of positive observations predicted to be positive, the FPR is the sample proportion of negative observations predicted to be positive, and the precision is the sample proportion of observations predicted to be positive that actually are positive. The probabilities in Definition~\ref{def:tpr_fpr} are the population analogues of the TPR and FPR in \eqref{eq:sample_metrics} when an observation $(X, Y)$ belongs to the positive class if $Y > y_0$ and the negative class otherwise. Similarly, the precision as defined in Definition~\ref{def:pred_properties} is the population analogue of the sample precision in \eqref{eq:sample_metrics}.
\end{remark}

It is a simple exercise in algebra to show that if
\begin{equation}\label{eq:rates_and_precision}
    \P(Y > y_0) = 1 - p, \ 
    \P(h(X) > \tau) = 1 - q, \ \text{and} \
    \P(Y > y_0 \mid h(X) > \tau) = \lambda,
\end{equation}
then we have the joint probabilities in Table~\ref{tab:joint_probs}, where, e.g., the entry in the upper left indicates that $\P(Y > y_0, h(X) > \tau) = (1 - q)\lambda$.
\begin{table}[!htb]
    \centering
    \begin{tabular}{@{}cc|cc@{}}
        \multicolumn{1}{c}{} &\multicolumn{1}{c}{} &\multicolumn{2}{c}{Observed Event} \\ 
        \multicolumn{1}{c}{} & 
        \multicolumn{1}{c|}{} & 
        \multicolumn{1}{c}{$\{Y > y_0\}$} & 
        \multicolumn{1}{c}{$\{Y \le y_0\}$} \\ 
        \cline{2-4}
        \multirow[c]{2}{*}{Predicted Event}
        & $\{h(X) > \tau\}$ & $(1 - q)\lambda$ & $(1 - q)(1 - \lambda)$ \\[1.5ex]
        & $\{h(X) \le \tau\}$ & $1 - p - (1 - q)\lambda$ & $p - (1 - q)(1 - \lambda)$ \\ 
        \cline{2-4}
    \end{tabular}
    \caption{Joint probabilities under the conditions in \eqref{eq:rates_and_precision}. The upper left entry indicates that $\P(Y > y_0, h(X) > \tau) = (1 - q)\lambda$, etc.}
    \label{tab:joint_probs}
\end{table}

It follows from Table~\ref{tab:joint_probs} that
\begin{align}
    \TPR &= \P(h(X) > \tau \mid Y > y_0)
    = \frac{\P(h(X) > \tau, Y > y_0)}{\P(Y > y_0)}
    = \frac{(1 - q)\lambda}{1 - p} \label{eq:tpr_expr} \\
    \FPR &= \P(h(X) > \tau \mid Y \le y_0)
    = \frac{\P(h(X) > \tau, Y \le y_0)}{\P(Y \le y_0)}
    = \frac{(1 - q)(1 - \lambda)}{p} \label{eq:fpr_expr}
\end{align}

\begin{remark}
    It can be shown that the two probabilities in the bottom row of Table~\ref{tab:joint_probs} are valid, i.e., in $[0, 1]$, if and only if
    \begin{equation}\label{ineq:lambda_bounds}
        \left(\frac{1 - p - q}{1 - q}\right) \vee 0
        \le \lambda
        \le \left(\frac{1 - p}{1 - q}\right) \wedge 1.
    \end{equation}
    At first glance, these inequalities seem to offer bounds on $\lambda$ that could be informative in practice. However, consideration of the values of $p$ and $q$ that would typically be used in practice makes it clear that this is not the case. We view the event $\{Y > y_0\}$ as extreme when $p \approx 1$. The natural choice for the calibration level $q$ is $p$. If $q = p$, then the lower bound in \eqref{ineq:lambda_bounds} is zero, and the upper bound is one. If $q \approx p$, then the lower bound would still be zero; the upper bound may be smaller than one, but it would not be that much smaller.
\end{remark}

\begin{remark}\label{rem:metric_equivs}
    In addition to the precision, there are other metrics that are often used to evaluate the performance of binary predictors. Some important examples are the True Skill Statistic (TSS), the Heidke Skill Score (HSS), and the $F_1$ score. We can use the probabilities in Table~\ref{tab:joint_probs} to express these as functions of $p$, $q$, and $\lambda$ that are strictly increasing in $\lambda$. Note that here we use the population definitions of these metrics; we have adapted these from the sample definitions in the references cited below. Because each metric is a strictly increasing function of $\lambda$, it follows that maximizing the precision while keeping the alarm rate equal to $1 - q$ is equivalent to maximizing any of the alternative metrics subject to the same constraint.
    
    The TSS is defined as $\TPR - \FPR$ \citep[][]{wilks2019stat}. Thus,
    \[
    \TSS
    := \TPR - \FPR
    = \frac{(1 - q)\lambda}{1 - p} - \frac{(1 - q)(1 - \lambda)}{p}
    = \frac{(1 - q)(\lambda + p - 1)}{p(1 - p)}.
    \]

    The definition of the HSS is more complicated; our explanation of the definition is based on the explanation in \cite{wilks2019stat}. The accuracy of $\I(h(X) > \tau)$ is defined as the probability that its prediction is correct; using Table~\ref{tab:joint_probs}, we can write this as
    \[
    \P(Y > y_0, h(X) > \tau) + \P(Y \le y_0, h(X) \le \tau) = (1 - q)\lambda + p - (1 - q)(1 - \lambda).
    \]
    If $X$ were independent of $Y$ and thus contained no information about it, then the accuracy would be
    \[
    \P(Y > y_0)\P(h(X) > \tau) + \P(Y \le y_0)\P(h(X) \le \tau) = (1 - p)(1 - q) + p q.
    \]
    This value can be viewed as a baseline that a predictor ought to improve upon. The maximum possible improvement is $1 - (1 - p)(1 - q) - p q$. The HSS is defined as the amount by which the predictor's accuracy improves upon the baseline, relative to the maximum possible improvement:
    \begin{align*}
        \HSS
        &:= \frac{(1 - q)\lambda + p - (1 - q)(1 - \lambda) - (1 - p)(1 - q) - p q}{1 - (1 - p)(1 - q) - p q} \\
        &= \frac{(1 - q)(2\lambda - 1) + p(1 - q) - (1 - p)(1 - q)}{1 - (1 - p - q + p q) - p q} \\
        &= \frac{(1 - q)[2\lambda - 1 + p - (1 - p)]}{p + q - 2 p q} \\
        &= \frac{2(1 - q)(\lambda + p - 1)}{p + q - 2 p q}.
    \end{align*}
    
    Using the definition of the $F_1$ score \citep[][]{murphy2022prob},
    \[
    F_1
    := \frac{2}{1 / \lambda + 1 / \TPR}
    = \frac{2}{1 / \lambda + (1 - p) / (1 - q)\lambda}
    = \frac{2(1 - q)\lambda}{2 - p - q}.
    \]

    An alternative to the preceding metrics is the Extremal Dependence Index (EDI), which was introduced in \cite{ferro2011extr}. The EDI is intended for evaluation of predictors of rare events and has several desirable properties, one of which is that as the event rate $1 - p$ goes to zero, the EDI does not tend to a trivial limit, unlike other metrics. The EDI is defined as 
    \[
    \frac{\log\FPR - \log\TPR}{\log\FPR + \log\TPR}.
    \]
    We show below that the derivative of the EDI with respect to $\lambda$ is positive on $(0, 1)$, so the EDI is also a strictly increasing function of $\lambda$, and maximizing the precision is equivalent to maximizing the EDI, assuming the same constraint is used for both.
    
    By \eqref{eq:tpr_expr} and \eqref{eq:fpr_expr}, we have
    \begin{equation}\label{eq:rate_deriv_signs}
        \frac{d\TPR}{d\lambda} = \frac{1 - q}{1 - p} > 0 \ \  \text{and} \ \ 
        \frac{d\FPR}{d\lambda} = -\frac{1 - q}{p} < 0,
    \end{equation}
    assuming that $0 < p, q < 1$. By the chain rule, we can write the derivative of the EDI with respect to $\lambda$ as
    \begin{equation}\label{eq:EDI_chain_rule}
        \frac{d\EDI}{d\lambda}
        = \frac{\partial\EDI}{\partial\TPR}\frac{d\TPR}{d\lambda} + \frac{\partial\EDI}{\partial\FPR}\frac{d\FPR}{d\lambda}.
    \end{equation}
    Let $D$ be the denominator of the EDI. We have
    \begin{align*}
        \frac{\partial\EDI}{\partial\TPR}
        &= \frac{1}{D^2}\left[-\frac{1}{\TPR}\left(\log\FPR + \log\TPR\right) - \frac{1}{\TPR}\left(\log\FPR - \log\TPR\right)\right] \\
        &= \frac{1}{D^2}\left(-\frac{2}{\TPR}\log\FPR\right); \\
        \frac{\partial\EDI}{\partial\FPR}
        &= \frac{1}{D^2}\left[\frac{1}{\FPR}\left(\log\FPR + \log\TPR\right) - \frac{1}{\FPR}\left(\log\FPR - \log\TPR\right)\right] \\
        &= \frac{1}{D^2}\left(\frac{2}{\FPR}\log\TPR\right).
    \end{align*}
    Plugging the expressions for the two partial derivatives into \eqref{eq:EDI_chain_rule} yields
    \begin{equation}\label{eq:EDI_deriv_expr}
        \frac{d\EDI}{d\lambda}
        = \frac{1}{D^2}\left[-\frac{2}{\TPR}\log\FPR \cdot \frac{d\TPR}{d\lambda} + \frac{2}{\FPR}\log\TPR \cdot \frac{d\FPR}{d\lambda}\right].
    \end{equation}
    Under our assumption that $0 < p, q < 1$, it follows from \eqref{eq:tpr_expr} and \eqref{eq:fpr_expr} that 
    $0 < \TPR, \FPR < 1$ for $\lambda \in (0, 1)$. Combining this observation with \eqref{eq:rate_deriv_signs}, it follows that the expression for the derivative in \eqref{eq:EDI_deriv_expr} is positive.
\end{remark}

We now turn to the notions of a family of predictors and an ROC curve based on a family.
\begin{definition}
    A family $\mathcal{F}(h)$ of predictors of $\I(Y > y_0)$ is a set of the form $\{\I(h(X) > \tau) : \tau \in \R\}$. We denote the precision, the TPR, and the FPR for a member of this family by $\lambda(\tau)$, $\TPR(\tau)$, and $\FPR(\tau)$, respectively.
\end{definition}
\begin{definition}
    The suboptimality region of the family $\mathcal{F}(h)$ is defined as
    \[
    D(h) := \bigcup_{\tau \in \R} \{(\FPR(\tau), \TPR) : 0 \le \TPR \le \TPR(\tau)\}.
    \]
\end{definition}
\begin{definition}
    Given two families $\mathcal{F}(h)$ and $\mathcal{F}(k)$ of predictors of $\I(Y > y_0)$, we say that $\mathcal{F}(h)$ dominates $\mathcal{F}(k)$ in ROC order if $D(h) \supseteq D(k)$. In this case, we write $\mathcal{F}(h) \ge_{ROC} \mathcal{F}(k)$.
\end{definition}

It turns out that the family based on an optimal predictor, i.e., based on the density ratio in Theorem~\ref{thm:base_thm}, dominates any other family in ROC order. To prove this, we first need to show that the optimization problem solved by Theorem~\ref{thm:base_thm}, the problem of maximizing the precision subject to a constraint on the alarm rate, is equivalent to an optimization problem that involves maximizing the TPR subject to a different constraint.

Let $E = \{Y > y_0\}$ and let $A = \{h(X) > \tau\}$. Abbreviate the alarm rate as AR and the false negative rate as FNR. Then
\[
\AR = \P(A), \
\text{precision} = \P(E \mid A), \
\TPR = \P(A \mid E), \
\FPR = \P(A \mid E^c), \
\FNR = \P(A^c \mid E).
\]
These conditions are equivalent:
\begin{align}
    \AR &= 1 - q \label{eq:AR_cond} \\
    (1 - p)\TPR + p\FPR &= 1 - q \label{eq:TPR_FPR_cond}
\end{align}
Condition \eqref{eq:AR_cond} implies condition \eqref{eq:TPR_FPR_cond} because
\begin{align*}
    (1 - p)\TPR + p\FPR &= (1 - p)\frac{\P(E \mid A)\P(A)}{\P(E)} + p\frac{\P(E^c \mid A)\P(A)}{\P(E^c)} \\
    &= \P(E \mid A)(1 - q) + \P(E^c \mid A)(1 - q) \\
    &= 1 - q,
\end{align*}
while condition \eqref{eq:TPR_FPR_cond} implies \eqref{eq:AR_cond} because
\begin{align*}
    1 - q &= (1 - p)\TPR + p\FPR \\
    &= (1 - p)\frac{\P(E \mid A)\P(A)}{\P(E)} + p\frac{\P(E^c \mid A)\P(A)}{\P(E^c)} \\
    &= \P(E \mid A)\P(A) + \P(E^c \mid A)\P(A) \\
    &= \AR.
\end{align*}

If conditions \eqref{eq:AR_cond} and \eqref{eq:TPR_FPR_cond} hold, then
\[
\text{precision}
= \P(E \mid A)
= \frac{\P(A \mid E)\P(E)}{\P(A)}
= \frac{1 - p}{1 - q}\P(A \mid E)
= \frac{1 - p}{1 - q}\TPR
= \frac{1 - p}{1 - q}(1 - \FNR).
\]
Thus, the following two problems are equivalent:
\begin{problem}[Solved by Theorem~\ref{thm:base_thm}]
    Maximize the precision subject to the constraint that $\AR = 1 - q$.
\end{problem}
\begin{problem}\label{prob:TPR_FPR_prob}
    Maximize the $\TPR$ (or minimize the $\FNR$) subject to the constraint that
    \[
    (1 - p)\TPR + p\FPR = 1 - q.
    \]
\end{problem}

\begin{remark}\label{rem:NP_comparison}
    The restatement in Problem~\ref{prob:TPR_FPR_prob} of the optimization problem solved by Theorem~\ref{thm:base_thm} clarifies the relationship between our approach and the Neyman-Pearson paradigm in classification. In classification, a covariate vector $X$ and a binary response $Y$ are given, and we want to construct a function of $X$, a classifier, that accurately predicts the value of $Y$. In the Neyman-Pearson paradigm, the goal is to find a classifier that minimizes the type II error rate subject to a constraint that bounds above the type I error rate, or vice versa \citep[][]{tong2016asur}. The Neyman-Pearson lemma can be used to show that the optimal classifier under this paradigm predicts that $Y = 1$ if and only if the density ratio $f_1(X) / f_0(X)$ lies above a threshold, where $f_i$ is the conditional density of $X$ given that $Y = i$ for $i = 0, 1$.

    Thus, if one wants to predict $\I(Y > y_0)$, then the Neyman-Pearson approach says to choose the predictor $\I(h(X) > \tau)$ that maximizes the TPR subject to the constraint $\FPR \le \alpha$ for a pre-specified $\alpha \in [0, 1]$. The objective function to be maximized is the same as in Problem~\ref{prob:TPR_FPR_prob}, but the constraint is different. However, we are interested in extreme scenarios, i.e., those in which $p \approx 1$. Under this condition, the constraint in \eqref{eq:TPR_FPR_cond} says that the FPR should be approximately $1 - q$, so it essentially controls the type I error rate, just like the constraint for the Neyman-Pearson problem.
\end{remark}

The proposition below says that a family based on an optimal predictor is dominant in ROC order under certain conditions.
\begin{proposition}\label{prop:opt_ROC_curve}
    Let $p \in (0, 1)$ and $y_0 \in \R$ such that $\P(Y \le y_0) = p$. Let $\mathcal{F}^* = \mathcal{F}(h)$ and $\mathcal{F} = \mathcal{F}(k)$ be two families of predictors of $\I(Y > y_0)$, with the former consisting of predictors that are optimal in the sense of Theorem~\ref{thm:base_thm}. Let $TPR^*(\tau)$, $FPR^*(\tau)$, $TPR(\tau)$, and $FPR(\tau)$ be the true and false positive rates for the two families.

    (i) $TPR^*$, $FPR^*$, $TPR$, and $FPR$ are decreasing functions of $\tau$.

    (ii) If $TPR$ and $FPR$ are continuous functions of $\tau$, then $\mathcal{F}^* \ge_{ROC} \mathcal{F}$.
\end{proposition}
\begin{proof}
    (i) Let $\tau_1, \tau_2 \in \R$ such that $\tau_1 \le \tau_2$. Since tail functions are decreasing,
    \begin{align*}
        \TPR^*(\tau_2) &= \P(k(X) > \tau_2 \mid Y > y_0)
        \le \P(k(X) > \tau_1 \mid Y > y_0)
        = \TPR^*(\tau_1); \\
        \FPR^*(\tau_2) &= \P(k(X) > \tau_2 \mid Y \le y_0)
        \le \P(k(X) > \tau_1 \mid Y \le y_0)
        = \FPR^*(\tau_1).
    \end{align*}
    Thus, $\TPR^*(\cdot)$ and $\FPR^*(\cdot)$ are decreasing. In a similar way, $\TPR(\cdot)$ and $\FPR(\cdot)$ can be shown to be decreasing.

    (ii) Let $(\FPR^*(\tau^*), \TPR^*(\tau^*))$ be a point on the ROC curve for $\mathcal{F}^*$. Consider the line with slope $-p / (1 - p)$ that passes through this point. The equation of this line must be
    \begin{equation}\label{eq:constraint_line}
        \TPR = -\frac{p}{1 - p}\FPR + \frac{1 - q}{1 - p},
    \end{equation}
    where $1 - q = \P(h(X) > \tau^*)$. Since $(0, 0)$ and $(1, 1)$ lie on opposite sides of the line, and the ROC curve for $\mathcal{F}$ is a continuous curve that connects those two points, the line intersects that curve as well. Let $(\FPR(\tau), \TPR(\tau))$ be the point of intersection. First assume that $(\FPR^*(\tau^*), \TPR^*(\tau^*))$ and $(\FPR(\tau), \TPR(\tau))$ are distinct.
    
    Since both points satisfy \eqref{eq:constraint_line}, both $\I(h(X) > \tau^*)$ and $\I(k(X) > \tau)$ satisfy the constraint in Problem~\ref{prob:TPR_FPR_prob}. Because the latter predictor is optimal, $\TPR^*(\tau^*) \ge \TPR(\tau)$. It follows that $\FPR^*(\tau^*) \le \FPR(\tau)$ as both points are on a line with negative slope. The fact that these points are distinct implies that these inequalities are actually strict, i.e., $\TPR^*(\tau^*) > \TPR(\tau)$ and $\FPR^*(\tau^*) < \FPR(\tau)$. Let $\tilde{\tau} \in \R$ such that $\FPR(\tilde{\tau}) = \FPR^*(\tau^*)$; $\tilde{\tau}$ exists because $\FPR(\cdot)$ is continuous. Then $\FPR(\tilde{\tau}) < \FPR(\tau)$, so $\tilde{\tau} \ge \tau$, and $\TPR(\tilde{\tau}) \le \TPR(\tau) \le \TPR^*(\tau^*)$. Thus, $(\FPR(\tilde{\tau}), \TPR(\tilde{\tau}))$ has the same first coordinate as $(\FPR^*(\tau^*), \TPR^*(\tau^*))$, but a second coordinate that is no larger.
    
    If $(\FPR(\tau), \TPR(\tau))$ and $(\FPR^*(\tau^*), \TPR^*(\tau^*))$ are the same point, then we can take $\tilde{\tau} = \tau$ above and derive the same conclusion. It follows that $\mathcal{F}^* \ge_{ROC} \mathcal{F}$.
\end{proof}

\begin{remark}
    A result similar to part (ii) of Proposition~\ref{prop:opt_ROC_curve} can be proven under the Neyman-Pearson paradigm described in Remark~\ref{rem:NP_comparison}; the optimal Neyman-Pearson predictor gives rise to an ROC curve that is optimal in that it lies above any other ROC curve at all values of the FPR. See implication (ii) of Result 4.4 in \cite{pepe2004thes} and Section 2.3.2 in \cite{krzanowski2009rocc}.
\end{remark}

\begin{remark}
    Some researchers argue that when there is class imbalance, it is better to consider the precision-recall (PR) curve instead of the ROC curve. The term ``recall'' is a synonym for the TPR. The argument is that a large increase in the number of false positives may not be apparent from the change in the FPR because of the large number of true negatives, but it would be apparent from the change in the precision \citep[][]{davis2006ther}. We can prove a result like part (ii) of Proposition~\ref{prop:opt_ROC_curve} that says that an optimal predictor gives rise to an optimal PR curve, a curve that lies above any other PR curve. The argument below is based on the proof of Theorem 3.2 in \cite{davis2006ther}.

    We use the same definitions and notation as in Proposition~\ref{prop:opt_ROC_curve}. Let $\lambda^*(\tau)$ and $\lambda(\tau)$ be the precisions and $\AR^*(\tau)$ and $\AR(\tau)$ the alarm rates for the two families. Suppose that the PR curve for $\mathcal{F}^*$ does not lie above that for $\mathcal{F}$. Then there exists a point $(\TPR(\tau), \lambda(\tau))$ and a point $(\TPR^*(\tau^*), \lambda^*(\tau^*))$ such that $\TPR(\tau) = \TPR^*(\tau^*)$, but $\lambda(\tau) > \lambda^*(\tau^*)$. Adapting \eqref{eq:tpr_expr},
    \[
    \TPR(\tau) = \frac{\AR(\tau)\lambda(\tau)}{1 - p} \text{ and }
    \TPR^*(\tau^*) = \frac{\AR^*(\tau^*)\lambda^*(\tau^*)}{1 - p},
    \]
    so we have $\AR(\tau)\lambda(\tau) = \AR^*(\tau^*)\lambda^*(\tau^*)$. Since $\lambda(\tau) > \lambda^*(\tau^*)$, it must be that $\AR(\tau) < \AR^*(\tau^*)$. Adapting \eqref{eq:fpr_expr}, we get
    \[
    \FPR(\tau)
    = \frac{\AR(\tau)(1 - \lambda(\tau))}{p}
    < \frac{\AR^*(\tau^*)(1 - \lambda^*(\tau^*))}{p}
    = \FPR^*(\tau^*).
    \]
    However, this yields a contradiction, as the dominance of the ROC curve for $\mathcal{F}^*$ over the ROC curve for $\mathcal{F}$ implies that $\FPR(\tau) \ge \FPR^*(\tau^*)$. Hence,  $\mathcal{F}^*$ must also have a dominant PR curve.
\end{remark}

\subsection{Additional Information for Section~\ref{sect:lin_ts_mod_opt_preds}}

\subsubsection{Proof of Lemma~\ref{lem:AR_opt_pred_recursion1}}\label{subsect:AR_opt_pred_recursion_lem_pf}

Lemma~\ref{lem:AR_opt_pred_recursion1} follows from the lemma below.
\begin{lemma}\label{lem:AR_opt_pred_recursion2}
    Let $\{Y_t\}$ be an AR($d$) process satisfying \eqref{eq:AR_mod}, with the innovations $\epsilon_t$ satisyfing Assumption~\ref{assump:innov_density}.

    (i) Define sequences of vectors $\{\phi(h)\}$ and $\{\psi(h)\}$ in $\R^d$ and $\R^h$, respectively, as follows. When $-(d - 1) \le h \le 0$, set
    \begin{align}
        \phi_j(h) &= \delta_{|h|, j}, \quad 0 \le j \le d - 1 \label{eq:phi_recur_rel1} \\
        \psi_j(h) &= 0, \quad 1 \le j \le h \nonumber
    \end{align}
    where $\delta$ is the Kronecker delta.

    For $h \ge 1$, define $\phi(h)$ and $\psi(h)$ recursively by
    \begin{align}
        \phi_j(h) &= \sum_{i = 1}^d \phi_i\phi_j(h - i), \quad 0 \le j \le d - 1 \label{eq:phi_recur_rel2} \\
        \psi_j(h) &=
        \begin{cases}
            \sum_{i = 1}^{(h - j) \wedge d} \phi_i\psi_j(h - i), & 1 \le j \le h - 1 \\
            1, & j = h
        \end{cases} \nonumber
    \end{align}
    
    Then for all $h \ge -(d - 1)$,
    \begin{equation}\label{eq:y_t_plus_h_expr}
        Y_{t + h} = \sum_{j = 0}^{d - 1} \phi_j(h)Y_{t - j} + \sum_{j = 1}^h \psi_j(h)\epsilon_{t + j},
    \end{equation}
    with the last sum having a density with respect to the Lebesgue measure.

    (ii) Let $\Phi$ be the $d \times d$ matrix $(\begin{matrix} \phi & e_1 & \cdots & e_{d - 1} \end{matrix})$, where $e_1, \ldots, e_{d - 1}$ are the first $d - 1$ standard basis vectors of $\R^d$. Then $\phi(h) = \Phi^h e_1$ for all $h \ge 1$.
\end{lemma}
\begin{proof}
    (i) Note that the claim is true when $-(d - 1) \le h \le 0$ because of how $\phi(h)$ and $\psi(h)$ are defined when $h$ is in that range. For $h \ge 1$, we prove the claim by induction on $h$. The base case is $h = 1$. We have
    \begin{align*}
        Y_{t + 1} &= \sum_{i = 1}^d \phi_i Y_{t + 1 - i} + \epsilon_{t + 1} \\
        &= \sum_{i = 1}^d \phi_i\left[\sum_{j = 0}^{d - 1} \phi_j(1 - i)Y_{t - j} + \sum_{j = 1}^{1 - i} \psi_j(1 - i)\epsilon_{t + j}\right] + \epsilon_{t + 1} \\
        &= \sum_{i = 1}^d \phi_i\sum_{j = 0}^{d - 1} \phi_j(1 - i)Y_{t - j} + \epsilon_{t + 1} \\
        &= \sum_{j = 0}^{d - 1} \phi_j(1)Y_{t - j} + \sum_{j = 1}^1 \psi_j(1)\epsilon_{t + j}.
    \end{align*}
    Suppose that the claim is true for integers smaller than some $h \ge 2$. Then
    \begin{align*}
        Y_{t + h} &= \sum_{i = 1}^d \phi_i Y_{t + h - i} + \epsilon_{t + h} \\
        &= \sum_{i = 1}^d \phi_i\left[\sum_{j = 0}^{d - 1} \phi_j(h - i)Y_{t - j} + \sum_{j = 1}^{h - i} \psi_j(h - i)\epsilon_{t + j}\right] + \epsilon_{t + h} \\
        &= \sum_{i = 1}^d \sum_{j = 0}^{d - 1} \phi_i \phi_j(h - i)Y_{t - j} + \sum_{i = 1}^d \sum_{j = 1}^{h - i} \phi_i\psi_j(h - i)\epsilon_{t + j} + \epsilon_{t + h} \\
        &= \sum_{j = 0}^{d - 1} \sum_{i = 1}^d \phi_i \phi_j(h - i)Y_{t - j} + \sum_{i = 1}^{(h - 1) \wedge d} \sum_{j = 1}^{h - i} \phi_i\psi_j(h - i)\epsilon_{t + j} + \epsilon_{t + h} \\
        &= \sum_{j = 0}^{d - 1} \phi_j(h)Y_{t - j} + \sum_{j = 1}^{h - 1} \sum_{i = 1}^{(h - j) \wedge [(h - 1) \wedge d]} \phi_i\psi_j(h - i)\epsilon_{t + j} + \epsilon_{t + h} \\
        &= \sum_{j = 0}^{d - 1} \phi_j(h)Y_{t - j} + \sum_{j = 1}^{h - 1} \sum_{i = 1}^{(h - j) \wedge d} \phi_i\psi_j(h - i)\epsilon_{t + j} + \epsilon_{t + h} \\
        &= \sum_{j = 0}^{d - 1} \phi_j(h)Y_{t - j} + \sum_{j = 1}^h \psi_j(h)\epsilon_{t + j},
    \end{align*}
    which completes the induction step. Under Assumption~\ref{assump:innov_density}, $\epsilon_{t + 1}, \ldots, \epsilon_{t + h}$ are independent and have densities with respect to the Lebesgue measure. Since $\psi_h(h) = 1$, it follows that the last sum also has a density with respect to the Lebesgue measure. 

    (ii) It follows from \eqref{eq:phi_recur_rel2} that for $h \ge 1$,
    \begin{equation}\label{eq:phi_recur_rel3}
        \phi(h) = (\begin{matrix} \phi(h - 1) & \cdots & \phi(h - d) \end{matrix})\phi.
    \end{equation}
    In particular, using \eqref{eq:phi_recur_rel1}, we get that
    \[
    \phi(1)
    = (\begin{matrix} \phi(0) & \cdots & \phi(-(d - 1)) \end{matrix})\phi
    = (\begin{matrix} e_1 & \cdots & e_d \end{matrix})\phi
    = \phi,
    \]
    $\{e_1, \ldots, e_d\}$ being the standard basis of $\R^d$. We show by induction that for every $h \ge 1$, $\Phi^h = (\begin{matrix} \phi(h) & \cdots & \phi(h - d + 1) \end{matrix})$. When $h = 1$, the claim is that
    \[
    \Phi
    = (\begin{matrix} \phi(1) & \cdots & \phi(-(d - 2)) \end{matrix})
    = (\begin{matrix} \phi & e_1 & \cdots & e_{d - 1} \end{matrix}),
    \]
    which is true by the definition of $\Phi$. Suppose that the claim is true for some $h \ge 1$. Then $\Phi^h\phi = (\begin{matrix} \phi(h) & \cdots & \phi(h - d + 1) \end{matrix})\phi = \phi(h + 1)$, by \eqref{eq:phi_recur_rel3}, and
    \[
    \Phi^{h + 1}
    = \Phi^h\Phi
    = (\begin{matrix} \Phi^h\phi & \Phi^h e_1 & \cdots & \Phi^h e_{d - 1} \end{matrix})
    = (\begin{matrix} \phi(h + 1) & \phi(h) & \cdots & \phi(h - d + 2) \end{matrix}).
    \]
    Hence, the claim is true for all $h \ge 1$, and $\phi(h) = \Phi^h e_1$.
\end{proof}

\subsubsection{The Optimal Predictor for Finite-Variance Autoregressive Models}\label{subsubsect:opt_pred_finite_var_ar}

Remark~\ref{rem:opt_pred_finite_var_ar} asserts that when the process $\{Y_t\}$ has finite variances, the random variable $\AROptPred{t}{h}{\phi}$ given in \eqref{eq:AR_opt_pred} is the same as the linear combination of $Y_t, \ldots, Y_{t - d + 1}$ that approximates $Y_{t + h}$ with minimal mean squared error (MSE), i.e.,
\begin{equation}\label{eq:AR_MSE_proof}
    \phi(h) = \argmin_{\varphi \in \R^d} \E[(Y_{t + h} - \varphi^{\top}Y_{t:(t - d + 1)})^2].
\end{equation}
This assertion can be proved as follows. Let $\gamma(\cdot)$ be the autocovariance function of $\{Y_t\}$. Using the notation of \cite{brockwell1991time}, let $\Gamma_d = (\gamma(i - j))_{i, j = 1, \ldots, d}$ and $\gamma_d^{(h)} = (\gamma(h), \ldots, \gamma(h + d - 1))^{\top}$. By equation (5.1.9) in \cite{brockwell1991time}, for any $\varphi \in \R^d$, the linear combination $\varphi^{\top}Y_{t:(t - d + 1)}$ minimizes the MSE in \eqref{eq:AR_MSE_proof} if and only if
\begin{equation}\label{eq:best_lin_pred_sys_of_eqs}
    \Gamma_d\varphi = \gamma_d^{(h)}.
\end{equation}
We shall show that \eqref{eq:best_lin_pred_sys_of_eqs} is satisfied when $\varphi = \phi(h)$. This will follow if for all $h \ge 0$,
\begin{equation}\label{eq:Gamma_d_Phi_relation}
    \Gamma_d \Phi^h = (\begin{matrix} \gamma_d^{(h)} & \cdots & \gamma_d^{(h - d + 1)} \end{matrix}),
\end{equation}
where $\Phi$ is as in Lemma~\ref{lem:AR_opt_pred_recursion1}, because \eqref{eq:Gamma_d_Phi_relation} implies that
\[
\Gamma_d \phi(h)
= \Gamma_d \Phi^h e_1
= \gamma_d^{(h)}.
\]
We prove \eqref{eq:Gamma_d_Phi_relation} by induction on $h$. The base case is $h = 0$; we have
\[
(\begin{matrix} \gamma_d^{(0)} & \cdots & \gamma_d^{(0 - d + 1)} \end{matrix})
= \left(\begin{matrix} \gamma(0) & \cdots & \gamma(-(d - 1)) \\ \vdots & \ddots & \vdots \\ \gamma(d - 1) & \cdots & \gamma(0) \end{matrix}\right).
\]
The $(i, j)$-entry of this matrix is $\gamma(1 - j + i - 1) = \gamma(i - j)$, so the matrix is the same as $\Gamma_d$. Now assume that the claim holds for some $h \ge 0$. Then
\[
\Gamma_d\Phi^{h + 1}
= \Gamma_d\Phi^h\Phi
= (\begin{matrix} \gamma_d^{(h)} & \cdots & \gamma_d^{(h - d + 1)} \end{matrix})\Phi
= (\begin{matrix} \phi_1\gamma_d^{(h)} + \cdots + \phi_d\gamma_d^{(h - d + 1)} & \gamma_d^{(h)} & \cdots & \gamma_d^{(h - d + 2)} \end{matrix}).
\]
The induction step will be complete if $\gamma_d^{(h + 1)} = \phi_1\gamma_d^{(h)} + \cdots + \phi_d\gamma_d^{(h - d + 1)}$. However, this follows from equations (3.3.8) and (3.3.9) in \cite{brockwell1991time}.

\subsubsection{A Consistency Lemma for Autoregressive Models}

\begin{lemma}\label{lem:AR_consistency1}
    Let $\phi$ be the vector of coefficients in model \eqref{eq:AR_mod}, and let $\widehat{\phi}(n) := \widehat{\phi}$ be an estimator of $\phi$, with $n$ being the size of the training set. Suppose that $\widehat{\phi} \xrightarrow{\P} \phi$ as $n \to \infty$. Then $\widehat{\phi}(h) \xrightarrow{\P} \phi(h)$ and $\approxAROptPred{t}{h}{\phi} \xrightarrow{\P} \AROptPred{t}{h}{\phi}$.
\end{lemma}
\begin{proof}
Let $\|\cdot\|_F$ be the Frobenius norm. Observe that
  \[
  \|\widehat{\Phi} - \Phi\|_F = \|(\begin{matrix} \widehat{\phi} - \phi & 0 & \cdots & 0 \end{matrix})\|_F = \sqrt{(\widehat{\phi} - \phi)^{\top}(\widehat{\phi} - \phi)} = \|\widehat{\phi} - \phi\|_2,
  \]
  so $\widehat{\Phi} \xrightarrow{\P} \Phi$. By the continuous mapping theorem, $\widehat{\Phi}^h \xrightarrow{\P} \Phi^h$. Because
  \[
  \|\widehat{\phi}(h) - \phi(h)\|_2 = \|(\widehat{\Phi}^h - \Phi^h)e_1\|_2 \le \|\widehat{\Phi}^h - \Phi^h\|_2,
  \]
  we may conclude that $\widehat{\phi}(h) \xrightarrow{\P} \phi(h)$. We have
  \[
  (\widehat{\phi}(h)^{\top}, Y_{t:(t - d + 1)}^{\top})^{\top} \xrightarrow{\P} (\phi(h)^{\top}, Y_{t:(t - d + 1)}^{\top})^{\top}
  \]
  since $\|(\widehat{\phi}(h)^{\top}, Y_{t:(t - d + 1)}^{\top})^{\top} - (\phi(h)^{\top}, Y_{t:(t - d + 1)}^{\top})^{\top}\|_2 = \|\widehat{\phi}(h) - \phi(h)\|_2$. Since the function $g : \R^{2 d} \to \R$ defined by
  \[
  g(u_0, \ldots, u_{d - 1}, v_0, \ldots, v_{d - 1}) = \sum_{k = 0}^{d - 1} u_k v_k
  \]
  is continuous, $\sum_{k = 0}^{d - 1} (\widehat{\phi}(h))_k Y_{t - k} \xrightarrow{p} \sum_{k = 0}^{d - 1} (\phi(h))_k Y_{t - k}$, or $\approxAROptPred{t}{h}{\phi} \xrightarrow{p} \AROptPred{t}{h}{\phi}$.
\end{proof}

\subsubsection{Proof of Proposition~\ref{prop:U_convergence}}\label{subsubsect:U_convergence_prop_pf}

\begin{proof}
    Set $\widetilde U_{t+h} :=  U_{t+h} (\widehat{\phi}(h),\,\approxAROptPred{t}{h}{\phi})$. Observe that by the triangle inequality,
    \begin{equation}\label{ineq:U_triangle_ineq}
        |\wh U_{t+h} - U_{t+h}^*| \le |\wh U_{t+h} - \widetilde U_{t+h}| + |\widetilde U_{t+h} - U_{t+h}^*|.
    \end{equation}
    We will show that the two terms on the right side converge in probability to zero.
    
    The first term is handled using Theorem 3.1 from \cite{adams2012unif}, which gives a uniform law of large numbers for stationary ergodic processes. That $\{Z_t = Y_{t:(t - d + 1)},\ t \in \mathbb Z\}$ is a stationary ergodic $\R^d$-valued time series readily follows from its moving average representation.  Namely, observe that without loss of generality, the iid sequence $\epsilon_t,\ t\in \Z$, supported on
    a probability space $(\Omega,{\cal F},\P)$ can always be written as $\epsilon_t(\omega):= \epsilon_0\circ T^t(\omega)$, for some invertible, 
    bi-measurable, measure-preserving, and ergodic transformation $T:\Omega \to \Omega$.  Therefore,
    $$
     Z_t = f(\epsilon_t,\epsilon_{t-1},\cdots) = Z_0\circ T^t,\ \ t\in \Z,
    $$
    where $f:\R^{\N} \to \R^d$ is ${\cal B}(\R^\N) \vert {\cal B}(\R^d)$ measurable, with ${\cal B}(\R^\N)$ denoting the minimal $\sigma$-field generated
    by all finite-dimensional cylinder sets with Borel faces.  Indeed, letting $f = (f_j)_{j=0}^{d - 1}$ in the AR$(d)$ case one can write $f$ explicitly as:
    $$
     f_j(x_0,x_1,\cdots) := \sum_{k=0}^\infty a_k(\phi) x_{j+k},  
    $$
    where $Y_t = \sum_{k=0}^\infty a_k(\phi) \epsilon_{t-k}$ is the causal moving average representation of the one-dimensional AR$(d)$ time-series.
    Now, the ergodicity of the map $T$ entails the ergodicity of the process $\{Z_t = Z_0\circ T^t \}$. In order to use Theorem 3.1 from \cite{adams2012unif}, we must also have that the collection of all half-spaces $H_{\varphi,\,\tau}$ is a separable 
    Vapnik-Chervonenkis class; this is shown in Lemma~\ref{lem:half_spaces_VC_class}. The latter uniform law of large numbers theorem now yields
    \begin{align*}
      | \wh U_{t+h} - \widetilde U_{t+h} |& = |\wh U_{t+h} (\varphi,\tau) 
     - U_{t+h} ( \varphi,\tau )| \Big{\vert}_{(\varphi,\tau)=(\widehat{\phi}(h),\,\approxAROptPred{t}{h}{\phi})} \\
     &\le \sup_{\varphi,\tau} |  \wh U_{t+h} (\varphi,\tau) -
       U_{t+h} (\varphi,\tau)| \stackrel{a.s.}{\longrightarrow }0,\ \ 
       \mbox{ as }n\to\infty,
    \end{align*}
    which shows that the first term on the right side in \eqref{ineq:U_triangle_ineq} vanishes.
    
    Now we handle the second term on the right side in \eqref{ineq:U_triangle_ineq}. 
    Recall that $Z$ in \eqref{e:U-varphi,tau} is chosen to be independent from 
    $(Y_{t:(t-d+1)}, \wh \phi)$.  Thus, $Z$ is independent from 
    $\wh Y_{t+h}(\wh \phi) = \wh\phi(h)^\top Y_{t:(t-d+1)}$ and $\wh \phi(h)$ and in 
    view of \eqref{e:U-varphi,tau},  we can write:
    \begin{equation}\label{eq:U_tilde_minus_U}
       \widetilde U_{t+h} - U_{t+h}^* = \E \left[ \I( \wh \phi(h)^\top Z \le \wh Y_{t+h}(\wh \phi) )
       - \I( \phi(h)^\top Z \le \wh Y_{t+h}(\phi) ) \Big\vert \widehat{\phi}, Y_{t:(t - d + 1)} \right] \\
    \end{equation}
   Introduce
   \[
   \xi_n := \wh \phi(h)^\top Z - \wh Y_{t+h}(\wh \phi) \ \text{and} \ \xi := \phi(h)^\top Z - \wh Y_{t+h}(\phi).
   \]
   In view of \eqref{eq:U_tilde_minus_U},
   \begin{align*}
       \E|\widetilde U_{t+h} - U_{t+h}^*| & = \E\left|\E\left[\I[\xi_n \le 0] - \I[\xi\le 0] \Big\vert \widehat{\phi}, Y_{t:(t - d + 1)} \right]\right| \\
       &\le \E|\I[\xi_n \le 0] - \I[\xi\le 0]|
    \end{align*}
    We will argue that the right-hand side vanishes. Indeed, since $\wh\phi\stackrel{\P}{\to} \phi$, it follows from Lemma~\ref{lem:AR_consistency1} that $\xi_n\stackrel{\P}{\to}\xi$. Because $\phi(h)^\top Z \overset{d}{=} \AROptPred{t}{h}{\phi}$, $\phi(h)^\top Z$ has a density with respect to Lebesgue measure. This together with the fact that $\phi(h)^{\top}Z \ind \AROptPred{t}{h}{\phi}$ implies that $\P[\xi=0] = 0$.  Therefore, by the continuous mapping theorem $\I(\xi_n\le 0) \stackrel{\P}{\to}\I(\xi\le 0)$, since the function $x\mapsto \I(x\le 0)$ is a.s.~continuous with respect to the law of $\xi$. (See also the proof of Lemma \ref{lem:stable-conv} for a more detailed argument). The fact that the latter difference of indicators is also 
    bounded entails 
    \[
    \E|\I[\xi_n \le 0] - \I[\xi\le 0]| \to 0,\ \ \mbox{ as $n\to\infty$. }
    \]
    This completes the proof.
\end{proof}

\subsubsection{Half-Spaces in $\R^d$ Form a Separable Vapnik-Chervonenkis Class}

The proof of Proposition~\ref{prop:U_convergence} requires the fact that the collection of all half-spaces in $\R^d$, i.e., sets of the form $\{y \in \R^d : \varphi^{\top}y \le \tau\}$ for some $\varphi \in \R^d$ and $\tau \in \R$, is a separable Vapnik-Chervonenkis (VC) class. In this section, we define the notions of VC class and separability and prove the fact, which is stated as Lemma~\ref{lem:half_spaces_VC_class}. The definitions below are adapted from those in \cite{adams2012unif}. In each definition, $\mathcal{C}$ denotes a collection of subsets of $\R^d$.

\begin{definition}
    The collection $\mathcal{C}$ shatters a finite set $F$ in $\R^d$ if every subset of $F$ can be expressed as $F \cap C$ for some $C \in \mathcal{C}$.
\end{definition}

\begin{definition}
    The collection $\mathcal{C}$ has VC dimension $d$ if $\mathcal{C}$ shatters some finite set of cardinality $d$, but does not shatter any finite set of cardinality $d + 1$. If $\mathcal{C}$ shatters sets of every finite cardinality, then $\mathcal{C}$ is said to have VC dimension $\infty$. The collection $\mathcal{C}$ is a VC class if it has finite VC dimension.
\end{definition}

\begin{definition}
    The collection $\mathcal{C}$ is separable if it has a countable subcollection $\mathcal{C}_0$ such that every $C \in \mathcal{C}$ can be written as $\lim_{k \to \infty} C_k$, where $C_k \in \mathcal{C}_0$ for each $k$.
\end{definition}

\begin{lemma}\label{lem:half_spaces_VC_class}
    For $\varphi \in \R^d$ and $\tau \in \R$, define $H_{\varphi,\,\tau} := \{y \in \R^d : \varphi^{\top}y \le \tau\}$. Then $\mathcal{H} = \{H_{\varphi,\,\tau} : \varphi \in \R^d, \, \tau \in \R\}$ is a separable VC class.
\end{lemma}
\begin{proof}
    That $\mathcal{H}$ is a VC class follows from Theorem 9.3 in \cite{shalev-shwartz2014unde}, which implies that $\mathcal{H}$ has VC dimension $d + 1$. It remains to show separability.

    Choose $\{\varphi^{(k)}\} \subset \Q^d$ so that $\varphi^{(k)} \to \varphi$ as $k \to \infty$ and $\varphi^{(k)}$ is not a multiple of $\varphi$ for any $k$. Fix a reference point $y^*$ in $\partial H_{\varphi,\,\tau}$, the hyperplane that forms the boundary of $H_{\varphi,\,\tau}$. Construct a sequence $\{\tau_k\} \subset \Q$ by choosing $\tau_k$ so that
    \begin{equation}\label{ineq:tau_k_bounds}
        \left(1 + \frac{1}{k}\right)\|\varphi^{(k)} - \varphi\|^{1 / 2}
        \ge \tau_k - (\varphi^{(k)})^{\top}y^*
        \ge \|\varphi^{(k)} - \varphi\|^{1 / 2}.
    \end{equation}
    The bounds in \eqref{ineq:tau_k_bounds} have several implications. One is that $y^*$ is an interior point of $H_{\varphi^{(k)},\,\tau_k}$ for every $k$. Another is that $\tau_k - (\varphi^{(k)})^{\top}y^* \to 0$ as $k \to \infty$. Since $(\varphi^{(k)})^{\top}y^* \to \varphi^{\top}y^* = \tau$, we must have $\tau_k \to \tau$ as well. We will prove that $H_{\varphi^{(k)},\,\tau_k} \to H_{\varphi,\,\tau}$ by showing that for any $y \in \R^d$,
    \begin{equation}\label{eq:half_space_lim}
        \lim_{k \to \infty} \I(y \in H_{\varphi^{(k)},\,\tau_k}) = \I(y \in H_{\varphi,\,\tau}).
    \end{equation}

    Suppose that $y$ is in the interior of $H_{\varphi,\,\tau}$; then $\varphi^{\top}y < \tau$. We have
    \[
    (\varphi^{(k)})^{\top}y - \tau_k
    = (\varphi^{(k)} - \varphi)^{\top}y - (\tau_k - \tau) + \varphi^{\top}y - \tau.
    \]
    It follows from this that for sufficiently large $k$, $(\varphi^{(k)})^{\top}y - \tau_k < 0$, i.e., $y$ is in the interior of $H_{\varphi^{(k)},\,\tau_k}$. Relation~\eqref{eq:half_space_lim} must then hold for $y$ since the indicator on the right side equals one and the indicator on the left side eventually equals one. The case $y \in H_{\varphi,\,\tau}^c$ can be handled similarly.

    Now suppose that $y \in \partial H_{\varphi,\,\tau}$. Let $k$ be such that $y \notin H_{\varphi^{(k)},\,\tau_k}$. The line segment $L = \{\lambda y + (1 - \lambda)y^* : 0 \le \lambda \le 1\}$ must lie inside $\partial H_{\varphi,\,\tau}$ since the boundary is convex. The mapping $\lambda \mapsto (\varphi^{(k)})^{\top}[\lambda y + (1 - \lambda)y^*] - \tau_k$ defined on $[0, 1]$ is continuous; it takes a negative value at $\lambda = 0$ and a positive value at $\lambda = 1$. Thus, there exists $\widetilde{\lambda} \in (0, 1)$ such that $\widetilde{y} = \widetilde{\lambda}y + (1 - \widetilde{\lambda})y^*$ satisfies $(\varphi^{(k)})^{\top}\widetilde{y} - \tau_k = 0$, i.e., $\widetilde{y} \in \partial H_{\varphi^{(k)},\,\tau_k}$. The point $\widetilde{y}$ is also closer to the reference point $y^*$ than is $y$:
    \begin{equation}\label{eq:y_tilde_y_star_dist}
        \|\widetilde{y} - y^*\|
        = \|\widetilde{\lambda}y + (1 - \widetilde{\lambda})y^* - y^*\|
        = \widetilde{\lambda}\|y - y^*\|.
    \end{equation}
    Since $\widetilde{y}$ is in both $\partial H_{\varphi,\,\tau}$ and $\partial H_{\varphi^{(k)},\,\tau_k}$, it is also in $\partial H_{\varphi^{(k)} - \varphi,\,\tau_k - \tau}$:
    \begin{equation}\label{eq:y_tilde_in_sum_boundary}
        \tau_k - \tau
        = (\varphi^{(k)})^{\top}\widetilde{y} - \varphi^{\top}\widetilde{y}
        = [\varphi^{(k)} - \varphi]^{\top}\widetilde{y}.
    \end{equation}
    It follows from \eqref{eq:y_tilde_y_star_dist} and \eqref{eq:y_tilde_in_sum_boundary} that $\|y - y^*\|
    \ge \delta(y^*, \partial H_{\varphi^{(k)} - \varphi,\,\tau_k - \tau})$, the distance between $y^*$ and $\partial H_{\varphi^{(k)} - \varphi,\,\tau_k - \tau}$. Choose any $z \in \partial H_{\varphi^{(k)} - \varphi,\,\tau_k - \tau}$. The distance can be computed as the length of the projection of $y^* - z$ along $\varphi^{(k)} - \varphi$, the normal vector to the hyperplane that is $\partial H_{\varphi^{(k)} - \varphi,\,\tau_k - \tau}$. Thus,
    \begin{align*}
        \|y - y^*\|
        &\ge \left\|\frac{(\varphi^{(k)} - \varphi)^{\top}(y^* - z)}{\|\varphi^{(k)} - \varphi\|^2}(\varphi^{(k)} - \varphi)\right\| \\
        &\ge \frac{|(\varphi^{(k)})^{\top}y^* - \varphi^{\top}y^* - (\varphi^{(k)} - \varphi)^{\top}z|}{\|\varphi^{(k)} - \varphi\|} \\
        &\ge \frac{\tau_k - (\varphi^{(k)})^{\top}y^*}{\|\varphi^{(k)} - \varphi\|} \\
        &\ge \|\varphi^{(k)} - \varphi\|^{-1 / 2},
    \end{align*}
    where the last two steps follow from \eqref{ineq:tau_k_bounds}. Since the final expression goes to $\infty$ as $k \to \infty$, there can only be finitely many $k$ for which $y \notin H_{\varphi^{(k)},\,\tau_k}$. Hence, $y$ lies in the $k$th half-space for large enough $k$, which means that \eqref{eq:half_space_lim} holds for $y$. This completes the proof that $H_{\varphi^{(k)},\,\tau_k} \to H_{\varphi,\,\tau}$.
\end{proof}

\subsubsection{A Simulation Study for an Autoregressive Model}\label{subsect:ar_mod_sim_study}

Optimal extreme event prediction for autoregressive models was discussed in Section~\ref{subsect:explicit_lin_ts_mod_opt_preds}; prediction when the AR coefficients are unknown was discussed in Section~\ref{subsect:AR_practical_pred}. Here, we use simulated data to evaluate the performance of various predictors for an AR model. We restate the model from Sections \ref{subsect:explicit_lin_ts_mod_opt_preds} and \ref{subsect:AR_practical_pred} for the sake of convenience:
\begin{equation}\label{eq:restated_ar_mod}
    Y_t = \sum_{i=1}^d \phi_i Y_{t-i} + \epsilon_t,\ \ \ t\in \Z.
\end{equation}
It was shown that for a stationary process satisfying \eqref{eq:restated_ar_mod}, for $h \ge 1$, the optimal predictor of $\I(Y_{t + h} > y_0)$ is $\I(\AROptPred{t}{h}{\phi} > \tau)$ for an appropriate choice of $\tau$ (Theorem~\ref{thm:AR_opt_pred}). In Section~\ref{subsubsect:AR_approx_opt_pred}, we discussed prediction of $\I(Y_{t + h} > \Fi{Y}(p))$, $\Fi{Y}(p)$ being the $p$th quantile of the marginal distribution of $Y_t$; we considered the use of the predictor $\I(\approxAROptPred{t}{h}{\phi} \ge \widehat{F}_{\approxAROptPred{t}{h}{\phi}}^{\leftarrow}(p))$, which can be estimated from training data, unlike the predictor given by Theorem~\ref{thm:AR_opt_pred}, $\I(\AROptPred{t}{h}{\phi} > \Fi{\AROptPred{t}{h}{\phi}}(p))$, which we shall call the oracle predictor.

The most straightforward way to compute the estimate $\widehat{F}_{\approxAROptPred{t}{h}{\phi}}^{\leftarrow}(p)$ of $\Fi{\approxAROptPred{t}{h}{\phi}}(p)$ is to use an empirical quantile, as was done in Section~\ref{subsubsect:AR_approx_opt_pred}. When $\{Y_t\}$ is a heavy-tailed process and $p$ is close to one, an alternative is to use an extreme quantile estimator based on extreme value theory 
\citep[see, e.g.,][and the references therein]{weissman1978esti,drees2003extr,chavez-demoulin2018extr}.

We conducted a simulation study that compared the performance of the oracle predictor, a predictor based on the empirical quantile estimator, and a predictor based on an extreme quantile estimator. Data was generated from the AR(5) model
\begin{equation}\label{eq:sim_ar_mod}
    Y_t = 0.3 Y_{t - 1} + 0.19 Y_{t - 2} - 0.035 Y_{t - 3} - 0.01 Y_{t - 4} + 0.0025 Y_{t - 5} + \epsilon_t,    
\end{equation}
with the $\epsilon_t$'s being iid $t$-distributed with one degree of freedom. The corresponding AR polynomial has roots $-3 \pm 1i$, $4 \pm 2i$, and $2$; since these roots lie outside the closed unit disk, it follows from Proposition 13.3.2 of \cite{brockwell1991time} that \eqref{eq:sim_ar_mod} has a unique strictly stationary, causal solution and that the solution is invertible. Set
\[
\phi = (0.3, 0.19, -0.035, -0.01, 0.0025).
\]

First, a large time series of size $10^6$ was generated from the model in \eqref{eq:sim_ar_mod}. For each $p \in \{0.90, 0.95, 0.99, 0.999\}$, $\Fi{Y}(p)$ was estimated as the $p$th sample quantile of the time series. All possible linear combinations of the form $\phi^{\top}Y_{t:(t - 4)}$ were computed using the time series; $\Fi{\AROptPred{t}{h}{\phi}}(p)$ was estimated as the $p$th sample quantile of these linear combinations.

Next, one hundred times, a time series of length $10^4 + 10^6$ was generated from the model in \eqref{eq:sim_ar_mod}; the first $10^4$ observations were reserved for training, and the last $10^6$ were reserved for testing. Inspired by \cite{davis1992mest}, LAD estimation was used on the training set to compute an estimate $\widehat{\phi}$ of $\phi$ for the non-oracle predictors; the order of the model was assumed to be known. Essentially, a median regression was carried out in which the response was the observation one step ahead and the covariates were the five most recent observations. Then, for the three predictors, one-step-ahead prediction was performed for each $p$; prediction worked the same way for each $p$.

For a test observation $Y_{t + 1}$, the oracle predictor predicted that $Y_{t + 1} > \Fi{Y}(p)$ iff 
\[
\phi^{\top}Y_{t:(t - 4)} > \Fi{\AROptPred{t}{h}{\phi}}(p).
\]
The predictor based on the empirical quantile estimator predicted an exceedance iff
\[
\widehat{\phi}^{\top}Y_{t:(t - 4)} \ge \widehat{F}_{\approxAROptPred{t}{h}{\phi}}^{\leftarrow}(p),
\]
where $\widehat{F}_{\approxAROptPred{t}{h}{\phi}}^{\leftarrow}(p)$ was computed as the $p$th sample quantile of all linear combinations of the form $\widehat{\phi}^{\top}Y_{s:(s - 4)}$ that could be computed from the training set.

The predictor based on the extreme quantile estimator predicted an exceedance iff
\[
\widehat{\phi}^{\top}Y_{t:(t - 4)} \ge \widetilde{F}_{\approxAROptPred{t}{h}{\phi}}^{\leftarrow}(p),
\]
where $\widetilde{F}_{\approxAROptPred{t}{h}{\phi}}^{\leftarrow}(p)$ was computed from the training set using the following procedure. For each of ten equally-spaced values between the third quartile of the training linear combinations and the tenth largest linear combination, a Generalized Pareto (GP) model was fit to the declustered excesses above the value. The excesses were the differences between the training observations above the value and the value, and they were declustered using the intervals method introduced in \cite{ferro2003infe}. The purpose of declustering was to extract excesses that could reasonably be viewed as independent, which was important since the GP model was fit using maximum likelihood estimation. The GP distribution with location parameter $\mu$, scale parameter $\sigma > 0$, and shape parameter $\xi$ has distribution function
\[
G(y) =
\begin{cases}
    1 - \left(1 + \xi\left(\frac{y - \mu}{\sigma}\right)\right)^{-1 / \xi}, & \xi \ne 0 \\
    1 - \exp\left(-\frac{y - \mu}{\sigma}\right), & \xi = 0,
\end{cases}
\]
with the support being $[\mu, \infty)$ when $\xi \ge 0$ and $[\mu, \mu - \sigma / \xi]$ when $\xi < 0$; GP distributions are frequently used to model excesses above a threshold \citep[][]{embrechts:kluppelberg:mikosch:1997, coles2001anin}.

The Kolmogorov-Smirnov test was used to identify the value for which the GP model fit best; the value with the smallest test statistic was used. Denote this value by $\tau$, and let $p_0$ be the proportion of linear combinations less than or equal to $\tau$. The extreme quantile estimate was computed as
\[
\widetilde{F}_{\approxAROptPred{t}{h}{\phi}}^{\leftarrow}(p) = \tau + \widehat{G}^{\leftarrow}\left(\frac{p - p_0}{1 - p_0}\right),
\]
the second term on the right being the $(p - p_0) / (1 - p_0)$ quantile of the fitted model $\widehat{G}$.

Once all predictions were made for the test set, the precision was computed as the number of exceedances that were predicted to be exceedances divided by the number of predicted exceedances. Since there were 100 training set-test set pairs, there were 100 precisions for each pair of a predictor and a value of $p$. Finally, the asymptotic precision of the optimal predictor was computed by deriving the causal representation of the solution to \eqref{eq:sim_ar_mod} and then applying Theorem~\ref{thm:MA_infty_opt_asymp_precision} using the first $10^6$ terms in each series.

\begin{figure}[!htb]
    \centering
    \includegraphics[width=0.7\textwidth]{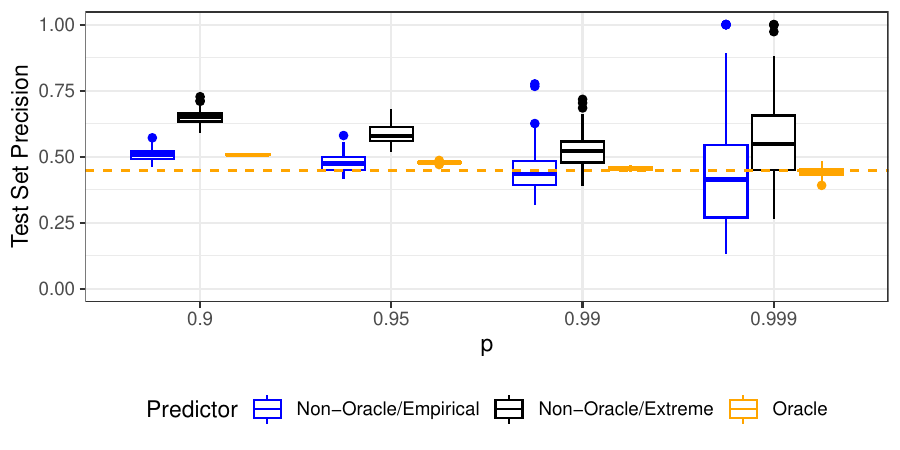}
    \caption{{\small Results from the autoregressive model simulation study. Data was generated from a particular AR(5) model. One hundred training set-test set pairs were generated; the model was fit to the training set, and the fitted model was used to predict exceedances of high quantiles in the test set. The quantile levels used are on the horizontal axis. On each of the 100 runs, the precision was calculated from the test set predictions. The dashed orange line represents the asymptotic precision of the optimal predictor.}}
    \label{fig:quantile_type_impact}
\end{figure}

The results of the simulation study are displayed in Figure~\ref{fig:quantile_type_impact}. Orange boxplots represent the precisions of the oracle predictor. As expected, as $p$ approaches one, the oracle precision converges to the asymptotic oracle precision, which is represented by the dashed orange line. The blue boxplots represent the precisions of the predictor based on the empirical quantile estimator; this predictor exhibits strong performance, as its precisions are mostly close to the corresponding oracle precisions. However, this is not the case when $p = 0.999$; most of the precisions are smaller than the corresponding oracle precisions. This is not surprising, as the training set size was $10,000$. The black boxplots represent the precisions of the predictor based on the extreme quantile estimator. Interestingly, for each $p$, all or almost all of these precisions lie well above the corresponding oracle precisions. It turned out that the extreme quantile estimates were substantially too large, so the predictor based on the extreme quantile estimator was calibrated at a higher level than the other predictors. This is not guaranteed to yield higher precisions, but it ended up doing so in this case.

\subsubsection{Regular Variation}\label{subsect:rv_info}

This section gives several definitions and auxiliary results that are used in Section~\ref{subsect:MA_infty_mod_opt_pred_asymp_aspects}. The definitions and results are based on those in \cite{kulik2020heav}.

\begin{definition}\label{def:rv_fun}
    A function $f : \R_+ \to \R_+$ is regularly varying (at $\infty$) if there exists $\alpha \in \R$ such that for all $x > 0$,
    \[
    \lim_{t \to \infty} \frac{f(t x)}{f(t)} = x^{\alpha}.
    \]
    In that case, we write $f \in \RV_{\alpha}$.
\end{definition}

\begin{definition}\label{def:sv_fun}
    A function $f : \R_+ \to \R_+$ is slowly varying if $f \in \RV_0$.
\end{definition}

\begin{definition}\label{def:rv_rand_var}
    A random variable $X$ is regularly varying with tail index $\alpha > 0$ and extremal skewness $p_X \in [0, 1]$ if
    \[
    \overline{F}_{|X|}(x) := \P(|X| > x) \in \RV_{-\alpha}
    \]
    and
    \[
    \lim_{x \to \infty} \frac{\overline{F}_X(x)}{\overline{F}_{|X|}(x)} = p_X.
    \]
\end{definition}

\begin{definition}\label{def:rv_rand_vec}
    A random vector $X \in \R^d$ is regularly varying with tail index $\alpha > 0$, scaling sequence $\{c_n\}_{n = 1}^{\infty}$, and exponent measure $\nu_X$ if
    \[
    \lim_{n \to \infty} n\P(c_n^{-1}X \in A) = \nu_X(A)
    \]
    for all Borel subsets $A$ of $\R^d \setminus \{0\}$ that are bounded away from zero and satisfy $\nu_X(\partial A) = 0$. The sequence $\{c_n\}$ is nonnegative and increasing, and the measure $\nu_X$ is a nonzero Borel measure on $\R^d \setminus \{0\}$ that is finite on sets that are bounded away from zero.
\end{definition}

\begin{remark}
    The scaling sequence-exponent measure pair isn't necessarily unique. However, if $\{c_n'\}$ and $\nu_X'$ is another pair, then as observed in \cite{kulik2020heav}, there exists $\gamma \in \R$ such that $\nu_X = \gamma \nu_X'$ and $c_n' \sim \gamma^{1 / \alpha}c_n$ as $n \to \infty$.
\end{remark}

\begin{remark}\label{rem:rv_rand_var_as_vec}
    A regularly varying random variable $X$ with tail index $\alpha$ and extremal skewness $p_X$ is also a regularly varying random vector with tail index $\alpha$. The scaling sequence can be chosen to be $\{c_n\}_{n = 1}^{\infty} = \{\overline{F}_{|X|}^{\leftarrow}(1 - 1 / n)\}_{n = 1}^{\infty}$, and the exponent measure $\nu_X$ can be defined by
    \[
    \nu_X(dx) = [p_X \I(x > 0) + (1 - p_X)\I(x < 0)]\alpha|x|^{-\alpha - 1}\,dx.
    \]
    Thus, for any Borel subset $A$ of $\R$ that is bounded away from zero and satisfies $\nu_X(\partial A) = 0$,
    \[
    \lim_{n \to \infty} n\P(c_n^{-1}X \in A) = \nu_X(A)
    \]
    It can be shown that $c_n \to \infty$ as $n \to \infty$. Note that
    \begin{equation}\label{eq:p_X_nu_X_rel}
        p_X
        = \lim_{x \to \infty} \frac{\P(X > x)}{\P(|X| > x)}
        = \lim_{n \to \infty} \frac{n\P(X > c_n)}{n\P(|X| > c_n)}
        = \lim_{n \to \infty} \frac{n\P(c_n^{-1}X \in (1, \infty))}{n\P(c_n^{-1}X \in [-1, 1]^c)}
        = \frac{\nu_X((1, \infty))}{\nu_X([-1, 1]^c)}.
    \end{equation}
\end{remark}

\begin{lemma}\label{lem:rv_lin_map}
    Let $X \in \R^d$ be a regularly varying random vector with tail index $\alpha$, scaling sequence $\{c_n\}_{n = 1}^{\infty}$, and exponent measure $\nu_X$. Let $A$ be an $m \times d$ matrix. Then $A X$ is regularly varying with tail index $\alpha$, scaling sequence $\{c_n\}_{n = 1}^{\infty}$, and exponent measure
    \[
    \nu_{A X} = \nu_X \circ A^{-1},
    \]
    as long as this measure is nonzero.
\end{lemma}
\begin{proof}
    Apply Corollary~2.1.14 of \cite{kulik2020heav}. 
    
\end{proof}

\subsubsection{Proof of Lemma~\ref{lem:xi_a_facts}}\label{subsect:xi_a_facts_lem_pf}

\begin{proof}
    Because Assumption~\ref{assump:ma_infty_innovs} is satisfied and $b \in \seqSet$, we can apply Corollary~4.2.1 of \cite{kulik2020heav} to conclude that $\xi(b)$ converges almost surely. We can also apply the corollary to obtain that $\xi(b)$ is regularly varying with tail index $\alpha$. As discussed in Remark~\ref{rem:rv_rand_var_as_vec}, a regularly varying random variable can be viewed as a regularly varying random vector. There exists an exponent measure $\nu_{b_j\epsilon_j}$ of $b_j\epsilon_j$ that equals $\nu_{\epsilon_j} \circ b_j^{-1}$, by Lemma~\ref{lem:rv_lin_map}. By Corollary~4.2.1 of \cite{kulik2020heav}, $\xi(b)$ has an exponent measure $\nu_{\xi(b)}$ such that
    \[
    \nu_{\xi(b)} = \sum_{j = 0}^{\infty} \nu_{b_j\epsilon_j}.
    \]
    Using \eqref{eq:p_X_nu_X_rel},
    \begin{equation}\label{eq:extr_skew_frac}
        p_{\xi(b)} = \frac{\nu_{\xi(b)}((1, \infty))}{\nu_{\xi(b)}([-1, 1]^c)} = \frac{\sum_{j = 0}^{\infty} \nu_{b_j\epsilon_j}((1, \infty))}{\sum_{j = 0}^{\infty} \nu_{b_j\epsilon_j}([-1, 1]^c)}
    \end{equation}
    If $a_j = 0$, then $\nu_{b_j\epsilon_j}$ is concentrated on $\{0\}$, so $\nu_{b_j\epsilon_j}((1, \infty)) = 0$. Otherwise,
    \begin{align}
        \nu_{a_j\epsilon_j}((1, \infty)) &= \nu_{\epsilon_j} \circ b_j^{-1}((1, \infty)) \nonumber \\
        &= \begin{cases}
            \nu_{\epsilon_j}((1 / b_j, \infty)), & b_j > 0 \\
            \nu_{\epsilon_j}((-\infty, 1 / b_j)), & b_j < 0
        \end{cases} \nonumber \\
        &= \begin{cases}
            |b_j|^{\alpha}\nu_{\epsilon_j}((1, \infty)), & b_j > 0 \\
            |b_j|^{\alpha}\nu_{\epsilon_j}((-\infty, -1)), & b_j < 0
        \end{cases} \nonumber \\
        &= \begin{cases}
            |b_j|^{\alpha}p_{\epsilon}, & b_j > 0 \\
            |b_j|^{\alpha}(1 - p_{\epsilon}), & b_j < 0
        \end{cases} \nonumber \\
        &= \kappa_+(b_j)|b_j|^{\alpha}, \label{eq:numer_expr}
    \end{align}
    which also gives the correct value when $b_j = 0$. 
    
    If $b_j = 0$, then $\nu_{b_j\epsilon_j}([-1, 1]^c) = 0$; if $b_j \ne 0$, then
    \begin{align}
        \nu_{b_j\epsilon_j}([-1, 1]^c) &= \nu_{\epsilon_j} \circ b_j^{-1}([-1, 1]^c) \nonumber \\
        &= \begin{cases}
            \nu_{\epsilon_j}([-1 / b_j, 1 / b_j]^c), & b_j > 0 \\
            \nu_{\epsilon_j}([1 / b_j, -1 / b_j]^c), & b_j < 0
        \end{cases} \nonumber \\
        &= \nu_{\epsilon_j}([-1 / |b_j|, 1 / |b_j|]^c) \nonumber \\
        &= |b_j|^{\alpha}\nu_{\epsilon_j}([-1, 1]^c) \nonumber \\
        &= |b_j|^{\alpha}, \label{eq:denom_expr}
    \end{align}
    which works when $b_j = 0$ as well.

    Substituting the expressions from \eqref{eq:numer_expr} and \eqref{eq:denom_expr} in \eqref{eq:extr_skew_frac} yields
    \[
    p_{\xi(b)}
    = \frac{\sum_{j = 0}^{\infty} \nu_{b_j\epsilon_j}((1, \infty))}{\sum_{j = 0}^{\infty} \nu_{b_j\epsilon_j}([-1, 1]^c)}
    = \frac{\normConst{+}{b}{0}}{\sum_{j = 0}^{\infty} |b_j|^{\alpha}}.
    \]
\end{proof}

\subsubsection{Proof of Lemma~\ref{lem:two_series_tdc}}\label{subsect:lem:two_series_tdc_lem_pf}

We first prove several lemmas needed in the proof of Lemma~\ref{lem:two_series_tdc}.
\begin{lemma}\label{lem:tail_fun_quant_fun_rel}
    Suppose that $X$ is regularly varying with tail index $\alpha$ and extremal skewness $p_X > 0$. Then as $p \uparrow 1$,
    \begin{equation}\label{eq:lambda_cond}
        \overline{F}_X(\Fi{X}(p)) \sim 1 - p.
    \end{equation}
\end{lemma}
\begin{proof}
    Since $\overline{F}_X \in \RV_{-\alpha}$, $1 / \overline{F}_X \in \RV_{\alpha}$; by Proposition 1.1.8 of \cite{kulik2020heav}, there exists $\ell_X \in \RV_0$ such that
    \[
    (1 / \overline{F}_X)(t) = t^{\alpha}\ell_X(t^{\alpha}) \ \ \text{and} \ \ (1 / \overline{F}_X)^{\leftarrow}(t) = t^{1 / \alpha}(\ell_X^{\#})^{1 / \alpha}(t),
    \]
    where $\ell_X^{\#}$ is the de Bruijn conjugate of $\ell_X$. We have
    \[
    \Fi{X}(p)
    = \inf\{x : F_X(x) \ge p\}
    = \inf\left\{x : (1 / \overline{F}_X)(x) \ge \frac{1}{1 - p}\right\}
    = (1 / \overline{F}_X)^{-1}\left(\frac{1}{1 - p}\right).
    \]
    Note that
    \[
    [\Fi{X}(p)]^{\alpha}
    = \left[(1 / \overline{F}_X)^{-1}\left(\frac{1}{1 - p}\right)\right]^{\alpha}
    = \left[\left(\frac{1}{1 - p}\right)^{1 / \alpha}(\ell_X^{\#})^{1 / \alpha}\left(\frac{1}{1 - p}\right)\right]^{\alpha}
    = \left(\frac{1}{1 - p}\right)\ell_X^{\#}\left(\frac{1}{1 - p}\right).
    \]
    Hence,
    \begin{align*}
        \lim_{p \uparrow 1} \frac{\overline{F}_X(\Fi{X}(p))}{1 - p} &= \lim_{p \uparrow 1} \frac{1}{(1 - p)(1 / \overline{F}_X)(\Fi{X}(p))} \\
        &= \lim_{p \uparrow 1} \frac{1}{(1 - p)[\Fi{X}(p)]^{\alpha}\ell_X([\Fi{X}(p)]^{\alpha})} \\
        &= \lim_{p \uparrow 1} \frac{1}{\ell_X^{\#}\left(\frac{1}{1 - p}\right)\ell_X\left(\frac{1}{1 - p}\ell_X^{\#}\left(\frac{1}{1 - p}\right)\right)} \\
        &= 1 \quad \text{\citep[by Theorem 1.5.13 of][]{bingham1987regu}},
    \end{align*}
    so \eqref{eq:lambda_cond} holds.
\end{proof}

\begin{lemma}\label{lem:two_seq_asymp_eq}
    Let $X$ be a regularly varying random variable with tail index $\alpha$ and extremal skewness $p_X > 0$. Let $\{a_n\}_{n = 1}^{\infty}$ and $\{b_n\}_{n = 1}^{\infty}$ be real sequences such that as $n \to \infty$, $a_n \to \infty$, $b_n \to \infty$, and $\overline{F}_X(a_n) \sim \overline{F}_X(b_n)$. Then $a_n \sim b_n$ as $n \to \infty$.
\end{lemma}
\begin{proof}
    Suppose that the claim is false. Then there exists $\epsilon > 0$ and a subsequence $\{n_k\}_{k = 1}^{\infty}$ of the sequence of natural numbers such that $a_{n_k} / b_{n_k} \notin (1 - \epsilon, 1 + \epsilon)$ for all $k$. Either $a_{n_k} / b_{n_k} \le 1 - \epsilon$ for infinitely many values of $k$ or $a_{n_k} / b_{n_k} \ge 1 + \epsilon$ for infinitely many values of $k$. Suppose that the former is true; assume without loss of generality that it is true for all $k$. Then
    \[
    \overline{F}_X(a_{n_k})
    = \overline{F}_X\left(\frac{a_{n_k}}{b_{n_k}}b_{n_k}\right)
    \ge \overline{F}_X((1 - \epsilon)b_{n_k}),
    \]
    so
    \begin{equation}\label{ineq:two_ratios}
        \frac{\overline{F}_X(a_{n_k})}{\overline{F}_X(b_{n_k})}
        \ge \frac{\overline{F}_X((1 - \epsilon)b_{n_k})}{\overline{F}_X(b_{n_k})}.
    \end{equation}
    Since $\overline{F}_X(a_n) / \overline{F}_X(b_n) \to 1$ as $n \to \infty$, $\overline{F}_X(a_{n_k}) / \overline{F}_X(b_{n_k}) \to 1$ as $k \to \infty$. Also, because $X$ is regularly varying with positive extremal skewness, $\overline{F}_X$ is regularly varying with index $-\alpha$ \citep[][]{kulik2020heav}. Combining these observations with \eqref{ineq:two_ratios}, we get
    \[
    1
    = \lim_{k \to \infty} \frac{\overline{F}_X(a_{n_k})}{\overline{F}_X(b_{n_k})} \\
    \ge \lim_{k \to \infty} \frac{\overline{F}_X((1 - \epsilon)b_{n_k})}{\overline{F}_X(b_{n_k})} \\
    = (1 - \epsilon)^{-\alpha},
    \]
    which is impossible. A similar contradiction arises if $a_{n_k} / b_{n_k} \ge 1 + \epsilon$ for infinitely many values of $k$. Hence, we must have $a_n \sim b_n$ as $n \to \infty$.
\end{proof}

The final lemma gives a sufficient condition for the existence of the tail dependence coefficient, which was defined in Definition~\ref{def:tdc}.
\begin{lemma}\label{lem:tdc_suff_cond}
    Let $X$ and $Y$ be random variables such that both of these conditions hold:
    
    (i) $\P(X > \Fi{X}(p)) \sim \P(Y > \Fi{Y}(p)) \sim 1 - p$ as $p \uparrow 1$.
    
    (ii) $\lim_{n \to \infty} n\P(X > \Fi{X}(1 - 1 / n), Y > \Fi{Y}(1 - 1 / n)) = \lambda \in [0, 1]$.
    
    Then the tail dependence coefficient $\lambda(X, Y)$ exists and equals $\lambda$:
    \[
    \lambda(X, Y)
    = \lim_{p \uparrow 1} \P(X > \Fi{X}(p) \mid Y > \Fi{Y}(p))
    = \lambda.
    \]
\end{lemma}
\begin{proof}
    Define $g : (0, 1) \to \R$ by $g(p) = \P(X > \Fi{X}(p), Y > \Fi{Y}(p))$; note that $g$ is increasing. Let $\{p_k\}$ be a sequence in $(0, 1)$ such that $p_k \uparrow 1$ as $k \to \infty$. For each $k$, set $n_k = \lfloor\frac{1}{1 - p_k}\rfloor$; as $k \to \infty$, $n_k \to \infty$. We have
    \begin{equation}\label{eq:tdc_expr1}
        \lambda = \lim_{k \to \infty} n_k\P(X > \Fi{X}(1 - 1 / n_k), Y > \Fi{Y}(1 - 1 / n_k)) = \lim_{k \to \infty} \frac{g\left(1 - 1 / n_k\right)}{1 / n_k}.
    \end{equation}
    We also have
    \begin{gather}
        n_k \le \frac{1}{1 - p_k} < n_k + 1 \nonumber \\
        \frac{1}{n_k + 1} < 1 - p_k \le \frac{1}{n_k}, \label{eq:1_minus_p_k}
    \end{gather}
    which implies that $1 - p_k \sim \frac{1}{n_k}$ and $1 - p_k \sim \frac{1}{n_k + 1}$ as $k \to \infty$. It also follows from \eqref{eq:1_minus_p_k} that
    \begin{gather*}
        1 - \frac{1}{n_k} \le p_k < 1 - \frac{1}{n_k + 1} \\
        \frac{g\left(1 - \frac{1}{n_k}\right)}{1 - p_k} \le \frac{g(p_k)}{1 - p_k} \le \frac{g\left(1 - \frac{1}{n_k + 1}\right)}{1 - p_k} \\
        \frac{g\left(1 - \frac{1}{n_k}\right)}{\frac{1}{n_k}} \cdot \frac{\frac{1}{n_k}}{1 - p_k} \le \frac{g(p_k)}{1 - p_k} \le \frac{g\left(1 - \frac{1}{n_k + 1}\right)}{\frac{1}{n_k + 1}} \cdot \frac{\frac{1}{n_k + 1}}{1 - p_k}.
    \end{gather*}
    Letting $k \to \infty$ on the last line and applying \eqref{eq:tdc_expr1} yields
    \begin{equation}\label{eq:tdc_expr2}
        \lim_{k \to \infty} \frac{g(p_k)}{1 - p_k} = \lambda.
    \end{equation}
    Hence,
    \begin{align*}
        \lambda(X, Y)
        &= \lim_{k \to \infty} \P(X > \Fi{X}(p_k) \mid Y > \Fi{Y}(p_k)) \\
        &= \lim_{k \to \infty} \frac{\P(X > \Fi{X}(p_k), Y > \Fi{Y}(p_k))}{\P(Y > \Fi{Y}(p_k))} \\
        &= \lim_{k \to \infty} \left[\frac{g(p_k)}{1 - p_k} \cdot \frac{1 - p_k}{\P(Y > \Fi{Y}(p_k))}\right] \\
        &= \lambda,
    \end{align*}
    where the last line follows from \eqref{eq:tdc_expr2} and Lemma~\ref{lem:tail_fun_quant_fun_rel}.
\end{proof}

We now prove Lemma~\ref{lem:two_series_tdc}.
\begin{proof}[Proof of Lemma~\ref{lem:two_series_tdc}]
    Let
    \[
    Z
    = \left(\begin{matrix} X \\ Y \end{matrix}\right)
    = \left(\begin{matrix} \series(a) \\ \series(b) \end{matrix}\right)
    = \sum_{j = 0}^{\infty} A_j\epsilon_j
    = \sum_{j = 0}^{\infty}
    \left(\begin{matrix} a_j \\ b_j \end{matrix}\right)\epsilon_j.
    \]
    Since $a, b \in \seqSet$, there exist $\delta_a, \delta_b \in (0, \alpha) \cap (0, 2]$ such that
    \[
    \sum_{j = 0}^{\infty} |a_j|^{\delta_a} <\infty \ \ \mbox{ and }\ \ 
    \sum_{j = 0}^{\infty} |b_j|^{\delta_b} < \infty.
    \]
    Set $\delta = \delta_a \vee \delta_b$. Since $\delta / 2 \le 1$,
    \[
    \|A_j\|_2^{\delta}
    = (a_j^2 + b_j^2)^{\delta / 2}
    \le (a_j^2)^{\delta / 2} + (b_j^2)^{\delta / 2}
    = |a_j|^{\delta} + |b_j|^{\delta}
    \]
    It follows that
    \[
    \sum_{j = 0}^{\infty} \|A_j\|_2^{\delta}
    \le \sum_{j = 0}^{\infty} |a_j|^{\delta} + \sum_{j = 0}^{\infty} |b_j|^{\delta}
    < \infty.
    \]
    Thus, we can apply Corollary 4.2.1 of \cite{kulik2020heav} to get that $Z$ is regularly varying with tail index $\alpha$ and exponent measure
    \[
    \nu_Z = \sum_{j = 0}^{\infty} \nu_{A_j\epsilon_j} =  \sum_{j = 0}^{\infty} \nu_{\epsilon} \circ A_j^{-1},
    \]
    where $\nu_{\epsilon}$ is defined by
    \[
    \nu_{\epsilon}(dt) = \kappa_+(t)\alpha|t|^{-\alpha - 1}\,dt.
    \]
    and $A_j^{-1}(S)$ is the pre-image of the set $S$ under the mapping $t \mapsto A_j t$. Note that if $S \cap \spn\{A_j\} = \emptyset$, then $\nu_{A_j \epsilon_j}(S) = 0$. Let $\{c_n\}$ be the scaling sequence associated with $\nu_Z$. We will show that
    \[
    \lambda(X, Y) = \lim_{n \to \infty} n\P(X > x_0 c_n, Y > y_0 c_n) = \nu_Z((x_0, \infty) \times (y_0, \infty))
    \]
    for some $x_0, y_0 > 0$ such that $x_0 c_n \sim \Fi{X}(1 - 1 / n)$ and $y_0 c_n \sim \Fi{Y}(1 - 1 / n)$ as $n \to \infty$. Examining the geometry of the support of $\nu_{A_j\epsilon_j}$ and considering cases will allow us to compute $\nu_{A_j\epsilon_j}((x_0, \infty) \times (y_0, \infty))$, which will then make it possible to compute $\nu_Z((x_0, \infty) \times (y_0, \infty)) = \sum_{j = 0}^{\infty} \nu_{A_j\epsilon_j}((x_0, \infty) \times (y_0, \infty))$.

    Since $X$ is regularly varying and $p_X = p_{\series(a)} > 0$, $\overline{F}_X(\Fi{X}(p)) \sim 1 - p$ as $p \uparrow 1$ by Lemma~\ref{lem:tail_fun_quant_fun_rel}. This implies that $\lim_{n \to \infty} n\P(X > \Fi{X}(1 - 1 / n)) = 1$.
    Suppose we could find an $x_0 > 0$ such that $\lim_{n \to \infty} n\P(X > x_0 c_n) = 1$. It would then follow from Lemma~\ref{lem:two_seq_asymp_eq} that as $n \to \infty$, $\Fi{X}(1 - 1 / n) \sim x_0 c_n$.

    For any $x_0 > 0$,
    \begin{align}
        \lim_{n \to \infty} n\P(X > x_0 c_n) &= \lim_{n \to \infty} n\P(X > x_0 c_n, Y \in \R) \nonumber \\
        &= \lim_{n \to \infty} n\P(c_n^{-1} Z \in (x_0, \infty) \times \R) \nonumber \\
        &= \nu_Z((x_0, \infty) \times \R) \nonumber \\
        &= x_0^{-\alpha}\nu_Z((1, \infty) \times \R) \nonumber \\
        &= x_0^{-\alpha}\sum_{j = 0}^{\infty} \nu_{A_j\epsilon_j}((1, \infty) \times \R) \label{eq:X_lim}
    \end{align}
    
    We have
    \[
    \nu_{A_j\epsilon_j}((1, \infty) \times \R) = \nu_{A_j\epsilon_j}([(1, \infty) \times \R] \cap \spn\{A_j\}).
    \]
    A vector is in the intersection above if and only if it equals $t A_j = t\left(\begin{matrix} a_j & b_j \end{matrix}\right)^{\top}$ for some $t \in \R$ and $t a_j > 1$. If $a_j = 0$, then this is impossible, so the intersection is empty. If $a_j > 0$, then we must have $t > 1 / a_j$, so
    \[
    \nu_{A_j\epsilon_j}((1, \infty) \times \R)
    = \nu_{A_j\epsilon_j}(\{t A_j : t > 1 / a_j\})
    = p_{\epsilon}\int_{1 / a_j}^{\infty} \alpha u^{-\alpha - 1}\,du
    = p_{\epsilon}a_j^{\alpha}.
    \]
    If $a_j < 0$, then we need $t < 1 / a_j$, so
    \[
    \nu_{A_j\epsilon_j}((1, \infty) \times \R)
    = \nu_{A_j\epsilon_j}(\{t A_j : t < 1 / a_j\})
    = (1 - p_{\epsilon})\int_{-\infty}^{1 / a_j} \alpha(-u)^{-\alpha - 1}\,du
    = (1 - p_{\epsilon})|a_j|^{\alpha}.
    \]
    Hence, whether $a_j$ is zero, positive, or negative,
    \[
    \nu_{A_j\epsilon_j}((1, \infty) \times \R) = \kappa_+(a_j)|a_j|^{\alpha}.
    \]

    Using \eqref{eq:X_lim}, in order to have $\lim_{n \to \infty} n\P(X > x_0 c_n) = 1$, we must have
    \[
    1
    = x_0^{-\alpha}\sum_{j = 0}^{\infty} \nu_{A_j\epsilon_j}((1, \infty) \times \R)
    = x_0^{-\alpha}\sum_{j = 0}^{\infty} \kappa_+(a_j)|a_j|^{\alpha} = x_0^{-\alpha}\normConst{+}{a}{0}.
    \]
    Since $p_X = p_{\series(a)} > 0$, Lemma~\ref{lem:xi_a_facts} implies that $\normConst{+}{a}{0} > 0$. Thus, we can set $x_0 = \normConst{+}{a}{0}^{1 / \alpha} > 0$. We have thus found an $x_0 > 0$ such that as $n \to \infty$, $x_0 c_n \sim \Fi{X}(1 - 1 / n)$. In a similar way, we can show that when $y_0 = \normConst{+}{b}{0}^{1 / \alpha}$, then $y_0 c_n \sim \Fi{Y}(1 - 1 / n)$ as $n \to \infty$.

    Now,
    \begin{align}
        \lim_{n \to \infty} n\P(X > x_0 c_n, Y > y_0 c_n)
        &= \lim_{n \to \infty} n\P(c_n^{-1} Z \in (x_0, \infty) \times (y_0, \infty)) \nonumber \\
        &= \nu_Z((x_0, \infty) \times (y_0, \infty)) \label{eq:nu_Z_val} \\
        &= \sum_{j = 0}^{\infty} \nu_{A_j\epsilon_j}((x_0, \infty) \times (y_0, \infty)) \nonumber \\
        &= \sum_{j = 0}^{\infty} \nu_{A_j\epsilon_j}\left([(x_0, \infty) \times (y_0, \infty)] \cap \spn\{A_j\}\right). \label{eq:tdc_sum}
    \end{align}
    A vector is in the intersection on the last line if and only if it equals $t A_j = t\left(\begin{matrix} a_j & b_j \end{matrix}\right)^{\top}$ for some $t \in \R$ and $t a_j > x_0, t b_j > y_0$. If one of $a_j, b_j$ equals zero, then this cannot happen. If $a_j, b_j > 0$, then the pair of inequalities becomes $t > x_0 / a_j, t > y_0 / b_j$, which hold if and only if $t > (x_0 / a_j) \vee (y_0 / b_j)$. Thus,
    \begin{align*}
        \nu_{A_j\epsilon_j}([(x_0, \infty) \times (y_0, \infty)] \cap \spn\{A_j\}) &= \nu_{A_j\epsilon_j}\left(\left\{t A_j : t > \frac{x_0}{a_j} \bigvee \frac{y_0}{b_j}\right\}\right) \\
        &= p_{\epsilon}\int_{(x_0 / a_j) \vee (y_0 / b_j)}^{\infty} \alpha u^{-\alpha - 1}\,du \\
        &= p_{\epsilon}\left(\frac{x_0}{a_j} \bigvee \frac{y_0}{b_j}\right)^{-\alpha} \\
        &= p_{\epsilon}\left(\frac{a_j^{\alpha}}{x_0^{\alpha}} \bigwedge \frac{b_j^{\alpha}}{y_0^{\alpha}}\right) \\
        &= p_{\epsilon}\left(\frac{a_j^{\alpha}}{\normConst{+}{a}{0}} \bigwedge \frac{b_j^{\alpha}}{\normConst{+}{b}{0}}\right)
    \end{align*}
    Similarly, if $a_j, b_j < 0$, then
    \[
    \nu_{A_j\epsilon_j}([(x_0, \infty) \times (y_0, \infty)] \cap \spn\{A_j\}) = (1 - p_{\epsilon})\left(\frac{|a_j|^{\alpha}}{\normConst{+}{a}{0}} \bigwedge \frac{|b_j|^{\alpha}}{\normConst{+}{b}{0}}\right).
    \]
    If $a_j > 0$ and $b_j < 0$, then $t a_j > x_0, t b_j > y_0$ become $t > x_0 / a_j, t < y_0 / b_j$, which cannot hold simultaneously as $x_0 / a_j > 0$ and $y_0 / b_j < 0$. If the signs of $a_j$ and $b_j$ are reversed, then we arrive at another impossibility. Hence, in all cases,
    \begin{equation}\label{eq:tdc_term}
        \nu_{A_j\epsilon_j}([(x_0, \infty) \times (y_0, \infty)] \cap \spn\{A_j\}) = \mmultiplier{+}{+}{a_j}{b_j}\left(\frac{|a_j|^{\alpha}}{\normConst{+}{a}{0}} \bigwedge \frac{|b_j|^{\alpha}}{\normConst{+}{b}{0}}\right).
    \end{equation}
    Substituting the expression on the right side of \eqref{eq:tdc_term} in \eqref{eq:tdc_sum} yields
    \begin{equation}\label{eq:tdc_expr}
        \lim_{n \to \infty} n\P(X > x_0 c_n, Y > y_0 c_n)
        = \sum_{j = 0}^{\infty} \mmultiplier{+}{+}{a_j}{b_j}\left(\frac{|a_j|^{\alpha}}{\normConst{+}{a}{0}} \bigwedge \frac{|b_j|^{\alpha}}{\normConst{+}{b}{0}}\right);
    \end{equation}
    the expression on the right side of \eqref{eq:tdc_expr} is the expression for $\lambda(\series(a), \series(b))$ in the statement of the lemma.
    
    We next show that
    \begin{equation}\label{eq:two_tdc_lims}
        \lim_{n \to \infty} n\P(X > \Fi{X}(1 - 1 / n), Y > \Fi{Y}(1 - 1 / n))
        = \lim_{n \to \infty} n\P(X > x_0 c_n, Y > y_0 c_n).
    \end{equation}
    Let $\epsilon > 0$. Since $x_0 c_n \sim \Fi{X}(1 - 1 / n)$ and $y_0 c_n \sim \Fi{Y}(1 - 1 / n)$ as $n \to \infty$, for sufficiently large $n$,
    \[
    \frac{\Fi{X}(1 - 1 / n)}{x_0 c_n}, \frac{\Fi{Y}(1 - 1 / n)}{y_0 c_n} \in (1 - \epsilon, 1 + \epsilon).
    \]
    so
    \[
    \Fi{X}(1 - 1 / n) \in ((1 - \epsilon)x_0 c_n, (1 + \epsilon)x_0 c_n)
    \ \text{and} \
    \Fi{Y}(1 - 1 / n) \in ((1 - \epsilon)y_0 c_n, (1 + \epsilon)y_0 c_n).
    \]
    Thus, for large enough $n$,
    \[
    n\P(X > \Fi{X}(1 - 1 / n), Y > \Fi{Y}(1 - 1 / n))
    \ge n\P(X > (1 + \epsilon)x_0 c_n, Y > (1 + \epsilon)y_0 c_n),
    \]
    which implies that
    \begin{align}
        \liminf_{n \to \infty} n\P(X > \Fi{X}(1 - 1 / n), Y > \Fi{Y}(1 - 1 / n))
        &\ge \liminf_{n \to \infty} n\P(X > (1 + \epsilon)x_0 c_n, Y > (1 + \epsilon)y_0 c_n) \nonumber \\
        &\ge \liminf_{n \to \infty} n\P(c_n^{-1}Z \in (1 + \epsilon)[(x_0, \infty) \times (y_0, \infty)]) \nonumber \\
        &\ge \nu_Z((1 + \epsilon)[(x_0, \infty) \times (y_0, \infty)]) \nonumber \\
        &\ge (1 + \epsilon)^{-\alpha}\nu_Z((x_0, \infty) \times (y_0, \infty)). \label{eq:liminf_lower_bound}
    \end{align}
    In a similar way, we can show that
    \begin{equation}\label{eq:limsup_upper_bound}
        \limsup_{n \to \infty} n\P(X > \Fi{X}(1 - 1 / n), Y > \Fi{Y}(1 - 1 / n))
        \le (1 - \epsilon)^{-\alpha}\nu_Z((x_0, \infty) \times (y_0, \infty)).
    \end{equation}
    Chaining together \eqref{eq:liminf_lower_bound} and \eqref{eq:limsup_upper_bound} and letting $\epsilon \downarrow 0$ yields
    \[
    \lim_{n \to \infty} n\P(X > \Fi{X}(1 - 1 / n), Y > \Fi{Y}(1 - 1 / n))
    = \nu_Z((x_0, \infty) \times (y_0, \infty)),
    \]
    which implies \eqref{eq:two_tdc_lims}, because of \eqref{eq:nu_Z_val}. Combining \eqref{eq:tdc_expr} and \eqref{eq:two_tdc_lims}, we get
    \begin{equation}\label{eq:tdc_expr3}
        \lim_{n \to \infty} n\P(X > \Fi{X}(1 - 1 / n), Y > \Fi{Y}(1 - 1 / n))
        = \sum_{j = 0}^{\infty} \mmultiplier{+}{+}{a_j}{b_j}\left(\frac{|a_j|^{\alpha}}{\normConst{+}{a}{0}} \bigwedge \frac{|b_j|^{\alpha}}{\normConst{+}{b}{0}}\right).
    \end{equation}    
    That $\lambda(\series(a), \series(b)) = \lambda(X, Y)$ exists and has the form specified in the statement of the lemma now follows from \eqref{eq:tdc_expr3} and Lemma~\ref{lem:tdc_suff_cond}, which is applicable because of Lemma~\ref{lem:tail_fun_quant_fun_rel}.
\end{proof}

\subsubsection{Proof of Proposition~\ref{prop:oracle}}\label{subsubsect:oracle_prop_pf}

\begin{proof}
Recall from \eqref{eq:lambda_opt_expr} that
\begin{equation}\label{e:p:oracle-lambda-opt}
\lambda(Y_{t + h}, \pred{opt}) = \sum_{j \ge 0} \kappa_+(a_{j + h})\frac{|a_{j + h}|^{\alpha}}{\normConst{+}{a}{0}}.
\end{equation}

On the other hand, by Lemma~\ref{lem:two_series_tdc},
\begin{align}\label{e:lambda-plus-long}
\lambda(Y_{t+h},Y_t) & =  \sum_{j\ge 0} \kappa_{+, +}(a_{j + h}, a_j)\Big(\frac{|a_{j+h}|^\alpha}{\normConst{+}{a}{0}} \bigwedge 
\frac{|a_{j}|^\alpha}{\normConst{+}{a}{0}} \Big)  \\
& = \sum_{j\ge 0} \kappa_{+, +}(a_{j + h}, a_j)\frac{|a_{j+h}|^\alpha}{\normConst{+}{a}{0}}, 
\label{e:lambda-plus}
\end{align}
where we used the lag-$h$ absolute monotonicity property $|a_{j+h}| \le |a_j|,\ j\ge 0$.

\medskip

\noindent {\em We now turn to proving (i).} If $a_{j+h}a_j\ge 0$ for all $j$, since $a_j=0$ entails $a_{j+h}=0$ 
(by \eqref{e:aj-absolute-monotone}), it follows that
\begin{align*}
    \kappa_{+, +}(a_{j + h}, a_j)
    &= p_\epsilon \I(a_{j+h}>0,a_j>0) + (1-p_\epsilon) \I(a_{j+h}<0,a_j<0) \\
    &= p_\epsilon \I(a_{j+h}>0) + (1-p_\epsilon) \I(a_{j+h}<0) \\
    &= \kappa_+(a_{j + h}),
\end{align*}
and hence the right-hand sides of \eqref{e:p:oracle-lambda-opt} and \eqref{e:lambda-plus} are identical.  This proves 
\eqref{e:p:oracle-i}. 

On the other hand, if 
$$
0 = \normConst{-}{a}{0} = \sum_{j\ge 0} \kappa_-(a_j)|a_j|^{\alpha} = \sum_{j\ge 0} (p_\epsilon \I(a_{j}<0) + 
(1-p_\epsilon) \I(a_j>0)) |a_j|^\alpha,
$$
it follows that we have either $p_\epsilon =1$ and $a_j\ge 0$ for all $j$, or $p_\epsilon =0$ and $a_j\le 0$ for all $j$. Note that our assumption that $\normConst{+}{a}{h} > 0$ rules out the possibility that $0 < p_{\epsilon} < 1$ and $a_j = 0$ for all $j$. In both cases $a_{j+h} a_j \ge 0$, for all $j\ge 0$, which as argued above, entails \eqref{e:p:oracle-i}.

\medskip
\noindent {\em Part {\em (iii)} is a special case of {\em (ii)}. We now focus on proving part {\em (ii)}.}
Applying Lemma~\ref{lem:two_series_tdc} again, we obtain
\begin{align}
 \label{e:lambda-minus}
\lambda(Y_{t+h},-Y_t) & =  \sum_{j\ge 0} \kappa_{+, -}(a_{j + h}, a_j)\Big(\frac{|a_{j+h}|^\alpha}{\normConst{+}{a}{0}} \bigwedge 
\frac{|a_{j}|^\alpha}{\normConst{-}{a}{0}} \Big),
\end{align}
For all $j\ge 0$, we have
\begin{equation}\label{e:p:oracle-2}
\Big(\frac{|a_{j+h}|^\alpha}{\normConst{+}{a}{0}} \bigwedge 
\frac{|a_{j}|^\alpha}{\normConst{+}{a}{0}} \Big)  \le 
\Big(1\bigvee \frac{\normConst{-}{a}{0}}{\normConst{+}{a}{0}} \Big)\cdot \Big(\frac{|a_{j+h}|^\alpha}{\normConst{+}{a}{0}} \bigwedge 
\frac{|a_{j}|^\alpha}{\normConst{-}{a}{0}} \Big),
\end{equation}
where we used the fact that $x\wedge ( w y) \le (1\vee w) (x \wedge y)$, for all non-negative $x, y,$ and $w$.
At the same time, for the sum of the weights in the right-hand sides of \eqref{e:lambda-plus-long} and \eqref{e:lambda-minus},
we have
\begin{align}\label{e:p:oracle-1}
\kappa_{+, +}(a_{j + h}, a_j) + \kappa_{+, -}(a_{j + h}, a_j)
&= (p_\epsilon \I( a_{j+h}>0,a_j>0) + (1-p_\epsilon) \I(a_{j+h}<0,a_j<0)) \nonumber\\
& \quad\quad + (p_\epsilon \I( a_{j+h}>0,a_j<0) + (1-p_\epsilon) \I(a_{j+h}<0,a_j>0)) \nonumber  \\
& = (p_\epsilon \I(a_{j+h}>0) + (1-p_\epsilon) \I(a_{j+h}<0)) \I(a_j\not = 0) \nonumber \\
&= \kappa_+(a_{j + h}),
\end{align}
where the last equality follows from the definition of $\kappa_+(a_{j + h})$ and lag-$h$ absolute monotonicity. In view of Relations \eqref{e:lambda-plus-long}, \eqref{e:lambda-minus}, \eqref{e:p:oracle-2}, and \eqref{e:p:oracle-1}, we obtain
$$
 \sum_{j\ge 0} \kappa_+(a_{j + h})\Big(\frac{|a_{j+h}|^\alpha}{\normConst{+}{a}{0}} \bigwedge 
\frac{|a_{j}|^\alpha}{\normConst{+}{a}{0}}\Big) \le \lambda(Y_{t+h},Y_t) + \Big( 1\bigvee \frac{\normConst{-}{a}{0}}{\normConst{+}{a}{0}}\Big)\lambda(Y_{t+h},-Y_t).
$$
By \eqref{e:aj-absolute-monotone}, however, we have $|a_{j+h}|^\alpha \wedge |a_j|^\alpha = |a_{j+h}|^\alpha$,
and in view of \eqref{e:p:oracle-lambda-opt} the left-hand side of the last relation equals $\lambda(Y_{t + h}, \pred{opt})$.  
This completes the proof of the upper bound in \eqref{e:p:oracle}.  In view of \eqref{e:p:oracle-lambda-opt},
and Relations \eqref{e:lambda-plus-long}, \eqref{e:lambda-minus}, and \eqref{e:p:oracle-1},
the lower bound in \eqref{e:p:oracle} is an immediate consequence of the inequalities:
$$
\frac{|a_{j+h}|^\alpha}{\normConst{+}{a}{0}} \bigwedge 
\frac{|a_{j}|^\alpha}{\normConst{\pm}{a}{0}} \le \frac{|a_{j+h}|^\alpha}{\normConst{+}{a}{0}}.
$$ 
\end{proof}

\subsection{Additional Information for Section~\ref{sect:data_analysis}}

\subsubsection{Extremal Optimal Precisions for FARIMA$(0, d, 0)$ Models}\label{subsubsect:farima_asymp_opt_precs}

Recall from Section~\ref{subsect:farima_background} that a FARIMA$(0, d, 0)$ model with symmetric $\alpha$-stable (S$\alpha$S) innovations is defined by the equation
\begin{equation}\label{eq:farima_0_d_0_mod2}
    Y_t = (1 - B)^{-d}\epsilon_t, \quad t \in \Z,    
\end{equation}
where $B$ is the backshift operator and the $\epsilon_t$'s are iid S$\alpha$S random variables. Also recall that when $d$ is non-integral, $\alpha \in (1, 2)$, and $d \in (0, 1 - 1 / \alpha)$, then this equation has a unique solution $Y_t = \sum_{j = 0}^{\infty} a_j\epsilon_{t - j}$ that is finite almost surely and invertible. By Theorem~\ref{thm:MA_infty_opt_asymp_precision}, the extremal optimal precision $\lambda_h^{\rm (opt)}$ for any $h \ge 1$ equals the ratio $\normConst{+}{a}{h} / \normConst{+}{a}{0}$, where $a = \{a_j\}_{j = 0}^{\infty}$. Furthermore, as stated in Remark~\ref{rem:farima_asymp_opt_precs}, $\lambda_h^{\rm (opt)}$ can be estimated by truncating the series defining the numerator and denominator using some large number of initial terms.

\begin{figure}[!htb]
    \centering
    \includegraphics[width=0.7\textwidth]{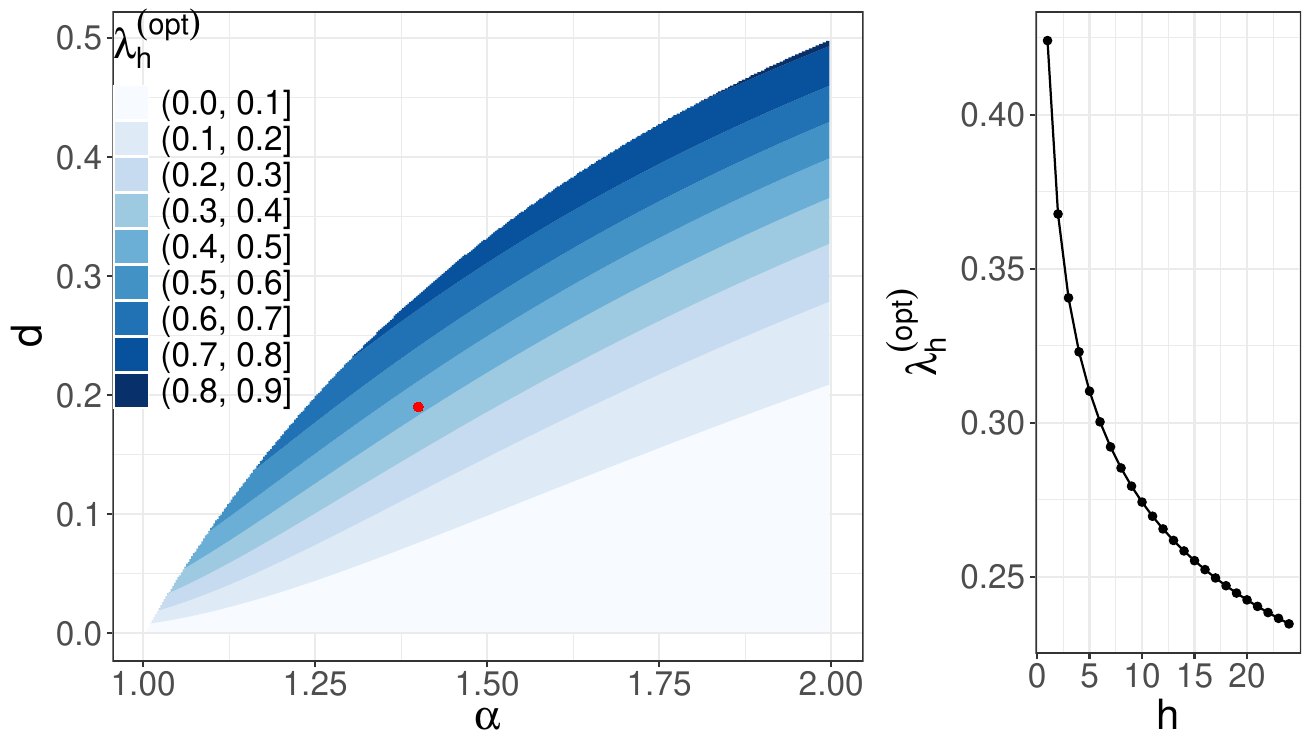}
    \caption{{\small \textbf{Left}: a contour plot showing how the one-step-ahead extremal optimal precision for the FARIMA(0, $d$, 0) model with symmetric $\alpha$-stable innovations \eqref{eq:farima_0_d_0_mod2} varies as $\alpha$ and $d$ vary. \textbf{Right}: the extremal optimal precision as a function of $h$, for $(\alpha, d) = (1.4, 0.19)$, which is represented by the red dot in the left panel.}}
    \label{fig:farima_asymp_precisions}
\end{figure}

For $h = 1$ and various pairs of $d$ and $\alpha$ satisfying $\alpha \in (1, 2)$, $d \in (0, 1 - 1 / \alpha)$, we approximated $\lambda_h^{\rm (opt)}$ by truncating the series using their first million terms. The approximations are shown in the left panel of Figure~\ref{fig:farima_asymp_precisions}. We see that for a fixed $\alpha$, as $d$ increases, $\lambda_h^{\rm (opt)}$ increases. This is intuitive as the level of long-range dependence in $\{Y_t\}$ increases as $d$ increases. For example, when $d = 0$, \eqref{eq:farima_0_d_0_mod2} simplifies to $Y_t = \epsilon_t$, $t \in \Z$, so $\{Y_t\}$ cannot be long-range dependent as it consists of independent random variables. The higher the level of long-range dependence, the more $Y_{t + h}$ depends on $Y_s$'s in the distant past, and the better the optimal predictor performs. On the other hand, for a fixed $d$, $\lambda_h^{\rm (opt)}$ decreases as $\alpha$ increases. The right panel of Figure~\ref{fig:farima_asymp_precisions} shows that for $\alpha = 1.4$ and $d = 0.19$, $\lambda_h^{\rm (opt)}$ decreases as the forecast horizon increases, as one would expect.

\subsubsection{Further Data Analysis Results}

\begin{table}[ht]
\centering
\begin{tabular}{cc|cccc|cccc}
  \hline
  & & \multicolumn{4}{c|}{Precision} & \multicolumn{4}{c}{TSS} \\
 $h$ & $p$ & Baseline & FARIMA & OLS AR & LAD AR & Baseline & FARIMA & OLS AR & LAD AR \\
 \hline
 1 & 0.90 & \textbf{0.492} & 0.448 & 0.354 & 0.442 & \textbf{0.435} & 0.392 & 0.276 & 0.419 \\
   1 & 0.95 & \textbf{0.379} & 0.368 & 0.258 & 0.332 & 0.334 & 0.320 & 0.225 & \textbf{0.352} \\
   1 & 0.99 & \textbf{0.244} & 0.238 & 0.104 & 0.148 & \textbf{0.262} & \textbf{0.262} & 0.124 & 0.230 \\
   \hline
 6 & 0.90 & 0.308 & \textbf{0.339} & 0.284 & 0.302 & 0.256 & \textbf{0.298} & 0.215 & 0.293 \\
   6 & 0.95 & 0.189 & \textbf{0.213} & 0.150 & 0.157 & 0.173 & \textbf{0.198} & 0.138 & 0.174 \\
   6 & 0.99 & 0.049 & 0.119 & \textbf{0.125} & 0.076 & 0.040 & 0.115 & \textbf{0.162} & 0.131 \\
   \hline
12 & 0.90 & 0.284 & \textbf{0.292} & 0.253 & 0.282 & 0.197 & 0.205 & 0.152 & \textbf{0.216} \\
  12 & 0.95 & \textbf{0.194} & 0.176 & 0.127 & 0.181 & 0.148 & 0.124 & 0.076 & \textbf{0.176} \\
  12 & 0.99 & \textbf{0.098} & 0.073 & 0.049 & 0.028 & \textbf{0.085} & 0.061 & 0.056 & 0.029 \\
   \hline
18 & 0.90 & 0.270 & \textbf{0.309} & 0.253 & 0.279 & 0.210 & 0.256 & 0.183 & \textbf{0.260} \\
  18 & 0.95 & 0.175 & \textbf{0.197} & 0.142 & 0.165 & 0.154 & 0.173 & 0.117 & \textbf{0.184} \\
  18 & 0.99 & 0.049 & \textbf{0.116} & 0.067 & 0.060 & 0.040 & \textbf{0.115} & 0.064 & 0.083 \\
   \hline
\end{tabular}
\caption{An expanded summary of the performance of various predictors. The order of the AR models was 168. For each pair of $h$ and $p$, the goal was to predict exceedance of the $p$th marginal quantile $h$ steps ahead. In each row, for each metric, the largest value is bolded. The LAD numbers are bigger than the corresponding OLS numbers, with the exception of the LAD precisions for $h \in \{6, 12, 18\}, p = 0.99$ and the LAD TSS values for $h \in \{6, 12\}, p = 0.99$. We believe that the LAD predictor performed better in most cases because it had a higher alarm rate in most cases; see Figure \ref{fig:lad_vs_ols}.}
\label{tab:expanded_data_analysis_results}
\end{table}

\begin{figure}[t!]
    \centering
    \includegraphics[width=0.32\textwidth]{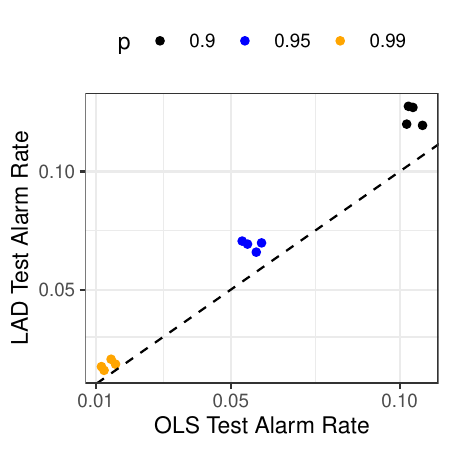}
    \caption{{\small A comparison of the test set alarm rates for the AR models fit using OLS and LAD. See Table~\ref{tab:expanded_data_analysis_results} for their precisions and TSS values. Different points represent different lead times $h$. Given $p$, the test alarm rate for the $p$th quantile is defined as the proportion of times in the test set for which an exceedance of the $p$th quantile was predicted. Along the dashed line, the OLS and LAD rates are equal. For every pair of $h$ and $p$, the LAD rate is higher than the OLS rate. A higher alarm rate corresponds to a lower threshold for predicting extremeness, which yields a higher TPR and a higher FPR. In this case, the changes led to increases in the precision and TSS for most pairs of $h$ and $p$.}}
    \label{fig:lad_vs_ols}
\end{figure}

\end{document}